\documentclass{article}
\usepackage[utf8]{inputenc}
\usepackage{hyperref}
\usepackage{xcolor}
\usepackage{bm} 
\usepackage{amssymb}
\usepackage{mathtools}
\usepackage{nicefrac}
\usepackage{cancel}
\PassOptionsToPackage{hyphens}{url}
\usepackage{hyperref}
\usepackage{amsthm}
\usepackage{enumerate}
\usepackage{multicol}
\usepackage{booktabs}
\usepackage{pgfplotstable} 
\usepackage{makecell}
\usepackage{colortbl}

\usepackage{textcomp}

\usepackage{algorithm}
\usepackage{algpseudocode}

\usepackage{siunitx}

\usepackage{subfig}
\usepackage{caption}
\usepackage[export]{adjustbox} 
\pgfplotsset{compat=newest}

\usepackage{localacronyms}

\title{Joint Reconstruction and Low-Rank Decomposition for Dynamic Inverse Problems}
\author{Simon Arridge\thanks{Department of Computer Science, University College London, London, United Kingdom  (\href{mailto:S.Arridge@cs.ucl.ac.uk}{S.Arridge@cs.ucl.ac.uk}).} \and Pascal Fernsel\thanks{Center for Industrial Mathematics, University of Bremen, Bremen, Germany 
		(\href{mailto:pfernsel@math.uni-bremen.de}{pfernsel@math.uni-bremen.de}). To whom correspondence should be addressed.} \and Andreas Hauptmann\thanks{Department of Mathematical Sciences, University of Oulu, Oulu, Finland; 
		Department of Computer Science, University College London, London, United Kingdom
  (\href{mailto:andreas.hauptmann@oulu.fi}{andreas.hauptmann@oulu.fi}).}}

\renewcommand{\vec}[1]{\bm{\mathit{#1}}}
\newcommand{\diff}{\mathop{} \! \mathrm{d}}

\DeclareMathOperator{\TV}{TV}
\DeclareMathOperator{\Var}{Var}
\DeclareMathOperator{\cov}{cov}

\theoremstyle{plain}
\newtheorem{lemma}{Lemma}[section]
\newtheorem{theorem}{Theorem}[section]
\newtheorem{corollary}{Corollary}[section]

\theoremstyle{definition}
\newtheorem{definition}{Definition}[section]

\definecolor{myGreen}{RGB}{25,142,33}

\newcommand{\R}{\mathbb{R}}

\newcommand{\N}{\mathbb{N}}

\newcommand{\scaleeq}[2]{\scalebox{#1}{\text{$\displaystyle #2$}}}

\newcommand{\ra}[1]{\renewcommand{\arraystretch}{#1}}

\begin{document}

\maketitle
\begin{abstract}
A primary interest in dynamic inverse problems is to identify the underlying temporal behaviour of the system from outside measurements. In this work we consider the case, where the target can be represented by a decomposition of spatial and temporal basis functions and hence can be efficiently represented by a low-rank decomposition. We then propose a joint reconstruction and low-rank decomposition method based on the Nonnegative Matrix Factorisation to obtain the unknown from highly undersampled dynamic measurement data. The proposed framework allows for flexible incorporation of separate regularisers for spatial and temporal features. For the special case of a stationary operator, we can effectively use the decomposition to reduce the computational complexity and obtain a substantial speed-up.
The proposed methods are evaluated for two simulated phantoms and we compare the obtained results to a separate low-rank reconstruction and subsequent decomposition approach based on the widely used principal component analysis. 
\newline

	
\end{abstract}


\section{Introduction} \label{sec:Introduction}

Several inverse problems are concerned with reconstruction of solutions in multiple physical dimensions such as space, time and frequency. Generally, such problems require very large datasets in order to satisfy conditions for accurate reconstruction, whereas in practice only subsets of such complete data can be measured. Furthermore, the information content of the solutions from such reduced data may be much less than suggested by the complete set.
In these cases, regularisation in the reconstruction process is required to compensate for the reduced information content, for instance by correlating features between auxiliary physical dimensions. 

For instance, dynamic inverse problems have gained considerable interest in recent years. This development is partly driven by the increase in computational resources and the possibility to handle large data size more efficiently, but also novel and more efficient imaging devices enabling wide areas of applications in medicine and industrial imaging. For instance in medical imaging, dynamic information is essential for accurate diagnosis of heart diseases 
or for applications in angiography to determine blood flow by injecting a contrast agent to the patient's blood stream. But also in nondestructive testing and chemical engineering, tomographic imaging has become increasingly popular to monitor dynamic processes.
The underlying problem in these imaging scenarios is often, that a fine temporal sampling, i.e. in the discrete setting a large number of channels, is only possible under considerable restrictions to sampling density at each time instance. This limitation often renders time-discrete (static) reconstructions insufficient.
Additionally, an underlying problem in many dynamic applications is given by the specific temporal dynamics of the process, which are often non-periodic and hence prevents temporal binning approaches. Thus, it is essential to include the dynamic nature of the imaging task in the reconstruction process.

With increasing computational resources, it has become more feasible to address the reconstruction task as a fully dynamic problem in a spatio-temporal setting.
In these approaches it is essential to include the dynamic information in some form into the reconstruction task. 
This could be done for instance by including a regularisation on the temporal behaviour as penalty in a variational setting \cite{Schmitt2002,Schmitt2002a}. Such approaches have been used in a wide variety of applications, such as magnetic resonance imaging \cite{feng2014golden,lustig2006kt,steeden2018real}, X-ray tomography \cite{bubba2017shearlet,niemi2015dynamic} and applications to process monitoring with electrical resistance tomography \cite{chen2018extended}. 
 
More advanced approaches aim to include a physical motion model and estimate the motion of the target from the measurements itself. This can be done for instance by incorporating an image registration step into the reconstruction algorithm and reformulate the reconstruction problem as a joint motion-estimation and reconstruction task \cite{burger2017variational,burger2018variational,djurabekova2019application,lucka2018enhancing}. Another possibility is the incorporation of an explicit motion model by methamorphsis as considered in \cite{Chen:2018aa,Gris:2019aa}.


In this work we consider another possibility to incorporate regularisation, and in particular temporal regularity, to the reconstruction task by assuming a low-dimensional representation of the unknown. This leads naturally to a low-rank description of the underlying inverse problem and is especially suitable to reduce data size in cases where we have much fewer basis functions to represent the unknown than the temporal sampling. In a continuous setting, this yields the analysis of low-rank approximations in tensor product of Hilbert spaces, for which we refer the reader to \cite{KressnerUschmajew16,Uschmajew13DoctoralThesis}. We rather focus on low-rank approximation methods in a discretised framework, which leads to the use of specific matrix factorisation approaches and their optimisation techniques.

In particular, in this work we propose a joint reconstruction and decomposition in a variational framework using non-negative matrix factorisation, which naturally represents the physical assumption of nonnegativity of the dynamic target and allows for a variety of regularising terms on spatial and temporal basis functions. Following this framework, we propose two algorithms, that either jointly recover the reconstruction and the low-rank decomposition, or alternatively recovers only the low-rank representation of the unknown without the need to construct the full spatio-temporal target in the reconstruction process. Here, the second approach effectively incorporates the dimension reduction and can lead under certain assumptions on the forward operator to a significant reduction in computational complexity. This can be particularly useful, if one is only interested in the dynamics of the system and not the full reconstruction. 

This paper is organised as follows. In Section \ref{sec:Framework} we discuss our setting for dynamic inverse problems and continue to discuss low-rank decomposition approaches. Specifically, principal component analysis (PCA) and non-negative matrix factorisation (NMF), which is the focus in this study. As a baseline we first present a low-rank reconstruction method followed by either of the decomposition methods. We then continue to present the proposed framework of joint reconstruction and decomposition with the NMF. In particular, we prove that the proposed framework leads to a monotonic decrease of the cost functions.
We then proceed in Section \ref{sec:Application to Dynamic CT} to evaluate the algorithms under considerations with the use case of dynamic X-ray tomography and two simulated phantoms with different characteristics. We conclude the study in Section \ref{sec:Conclusion and Outlook} with some thoughts on the extension of the proposed framework.


\section{Reconstruction and Low-rank Decomposition Methods} \label{sec:Framework}
    \subsection{A Setting for Dynamic Inverse Problems}\label{subsec:Multi--Channel Inverse Problem}
        
        In this work, we consider a general multi-dimensional inverse problem, where the unknown $x(s,\tau)$ is defined on a spatial domain $\Omega_1 \subset \R^{d_1}$ with dependence on a secondary variable $t\in \R_{\geq 0}$ defined in a bounded interval $\mathcal{T}:=[0,T]$. This setting admits some quite general applications where the secondary variable could have other physical interpretations, notably wavelength for hyper-spectral problems; however, to fix our ideas, we henceforth consider $t$ to explicitly represent time, and our application to be that of \emph{dynamic inverse problems}. Consequently, the underlying equation of the resulting inverse problem can be described in the following form
	    \begin{equation}\label{eq:dynamic IP General Continuous}
	        \mathcal{A}(x(s,t);t) = y(\sigma,t) \quad \text{for } t \in \mathcal{T},
	    \end{equation}
	    where $\mathcal{A}:\Omega_1\times \mathcal{T}\to \Omega_2\times \mathcal{T}$ 
	    is a time-dependent linear bounded operator between suitable Hilbert spaces and $\Omega_2\subset\R^{d_2}$ is the domain for the measurement data. We will primarily consider the non-stationary case here, where the forward operator $\mathcal{A}$ is dependent on $t$.
	    In the special case of a stationary operator $\mathcal{A}(\cdot;t) \equiv \mathcal{A}$ for all $t\in \mathcal{T}$, where for each $t$ the operator follows the same sampling process, we can achieve possible computational improvements. The resulting implications will be discussed later in Section \ref{subsec:Complexity Reduction for Stationary Operator}.
	    
\begin{figure}
        	\centering
        	\includegraphics[width=\textwidth]{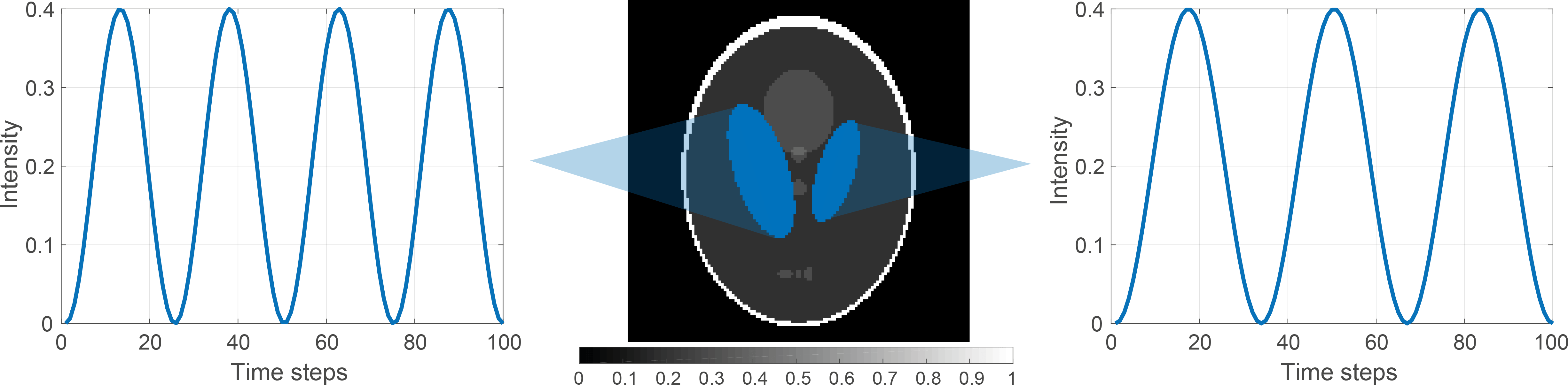}
        	\caption{Illustration of a phantom that can be represented by the decomposition in \eqref{eqn:continuousDeompostion}. The phantom consists of $K=3$ components: the background and two dynamic components with periodically changing intensity (left and right plot). As such, this phantom can be efficiently represented by a low-rank decomposition considered in this study.} \label{fig:dynSheppIllustration}
        \end{figure}

	    Furthermore, the underlying assumption in this work is that the unknown $x:\Omega_1\times\mathcal{T}\to\R_{\geq 0 }$ 
	    can be decomposed into a set of spatial $b^k:\Omega_1\to \R_{\geq 0}$ and channel basis functions $c^k(t):\mathcal{T}\to\R_{\geq 0}$ for $1 \leq k\leq K$. Then the unknown can be represented by the decomposition
	    \begin{equation}\label{eqn:continuousDeompostion}
	        x(s,t) = \sum_{k=1}^K b^k(s) c^k(t).
	    \end{equation}
	    This formulation naturally gives rise to the reconstruction and low-rank decompostion framework to extract the relevant features given by $b^k$ and $c^k$.  An illustration for a possible phantom represented by \eqref{eqn:continuousDeompostion} is shown in Figure \ref{fig:dynSheppIllustration}.
	    
	    We intentionally keep the formulation general here to allow for applications different to dynamic inverse problems, such as multi-spectral imaging. Nevertheless, the derived reconstruction and feature extraction framework in this paper will be used in Section \ref{sec:Application to Dynamic CT} for the specific application to dynamic computed tomography.
	    
	    Furthermore, a suitable discretisation of the continuous formulation \eqref{eq:dynamic IP General Continuous} is needed to introduce the feature extraction methods in the forthcoming sections. 
	    Let us first discretise the secondary variable, such that $t\in\N$ with $1\leq t \leq T$. For the spatial domain, we assume a vectorised representation such that the resulting unknown can be represented as a matrix $X\in\R^{N\times T}$, which leads to the matrix equation 
	    \begin{equation}\label{eq:dynamic IP General Discrete}
	        A_t X_{\bullet, t} = Y_{\bullet, t} \quad \text{for} \quad 1\leq t\leq T,
	    \end{equation}
	    where $A_t\in \R^{M\times N} $ is the discretised forward operator, $X_{\bullet, t}$ the $t$-th column of $X$ and $Y_{\bullet, t}$ the $t$-th column of the data matrix $Y\in \R^{M\times T}.$ Analogously, we will write $M_{n, \bullet}$ for the $n$-th row of an arbitrary matrix $M.$
	    
	    Suitable restrictions to the matrices in Equation \eqref{eq:dynamic IP General Discrete} will be made in the following sections to ensure the applicability of the considered frameworks and, if possible, to properly represent the decomposition \eqref{eqn:continuousDeompostion}. 
	    
    \subsection{Feature Extraction Methods} \label{subsec:Feature Extraction Methods}
        In this section, we introduce two feature extraction methods, namely the \textit{\acl{PCA}}~(\acs{PCA}) and the \textit{\acl{NMF}}\ (\acs{NMF}). These approaches are used to compute the latent components of the reconstruction $X.$ 
        The \ac{NMF} will be used in Section \ref{subsec:Joint Reconstruction and low--rank Decomposition} to introduce a joint reconstruction and low-rank decomposition framework to tackle the problem stated in \eqref{eq:dynamic IP General Discrete}.
        \subsubsection{\acl{PCA}} \label{subsubsec: Principal Component Analysis}
            Large and high dimensional datasets demand modern data analysis approaches to reduce the dimensionality and increase the interpretability of the data while keeping the loss of information as low as possible. Many different techniques have been developed for this purpose, but \ac{PCA} is one of the most widely used and goes back to \cite{pearson1901-PCA}.
            
            For a given matrix $X\in \R^{N\times T}$ with $N$ different observations of an experiment and $T$ features, the PCA is a linear orthogonal transformation given by the weights $C_{\tilde{k}, \bullet} = (C_{\tilde{k}1}, \dots, C_{\tilde{k}T})$ with $C \in \R^{\tilde{K}\times T},$ which transforms each observation $X_{n, \bullet}$ to \textit{principal component scores} given by $B_{n\tilde{k}} \coloneqq \sum_t  X_{nt} C_{\tilde{k}t} $ with $B=[B_{\bullet, 1}, \dots, B_{\bullet, \tilde{K}}]\in \R^{N\times \tilde{K}}$ and $\tilde{K}=\min(N-1, T),$ such that
            \begin{itemize}
                \item the sample variance $\Var(B_{\bullet, \tilde{k}})$ is maximised for all $\tilde{k},$
                \item each row $C_{\tilde{k}, \bullet}$ is constrained to be a unit vector
                \item and the sample covariance $\cov(B_{\bullet, k}, B_{\bullet, \tilde{k}}) = 0 $ for $k\neq \tilde{k}.$
            \end{itemize}
            Together with the usual assumption that the number of observations is higher than the underlying dimension, this leads to $\tilde{K}=T$ and the full transformation $B = XC^\intercal,$ where $C$ is an orthogonal matrix. The $t$-th column vector $(C_{t, \bullet})^\intercal$ defines the $t$-th \textit{principal direction} and is an eigenvector of the covariance matrix $S=X^\intercal X/(N-1).$ The corresponding $t$-th largest eigenvalue of $S$ denotes the variance of the $t$-th principal component.
            
            The above transformation is equivalent to the factorisation of the matrix $X$ given by
            \begin{equation}\label{eq:PCA:Factorisation}
                X = BC,
            \end{equation}
            which allows to decompose each observation into the principal components, such that $X_{n, \bullet} = \sum_{t=1}^T B_{nt} C_{t, \bullet}.$ Therefore, we also have $X = \sum_{t=1}^{T} B_{\bullet, t} C_{t, \bullet}$ similarly to the continuous case in \eqref{eqn:continuousDeompostion}.
            
            Furthermore, it is possible to obtain an approximation of the matrix $X$ by truncating the sum at the first $K<T$ principle components for all $n,$ which yields a rank $K$ matrix $X^{(K)}$ given by
            \begin{equation*}
                X^{(K)} = \sum_{k=1}^K B_{\bullet, k} C_{k, \bullet}.
            \end{equation*}
            Based on the Eckart–Young–Mirsky theorem \cite{GolubVanLoan13-EckardYoung}, $X^{(K)}$ is the best rank $K$ approximation of $X$ in the sense that it minimises the discrepancy $\Vert X - X^{(K)}\Vert$ for both the Frobenius and spectral norm.
            
            One typical approach to compute the \ac{PCA} is based on the \ac{SVD} of the data matrix $X=U\Sigma V^\intercal$ and will be used in this work. Setting $B\coloneqq U\Sigma$ and $C=V^\intercal$ gives already the desired factorisation in \eqref{eq:PCA:Factorisation} based on the \ac{PCA}.
        \subsubsection{\acl{NMF}} \label{subsubsec:Nonnegative Matrix Factorisation}
            \acf{NMF}, originally introduced as positive matrix factorisation by Paatero and Tapper in 1994 \cite{paatero94NMFIntro}, is an established tool to obtain low-rank approximations of nonnegative data matrices. It has been widely used in the machine learning and data mining community for compression, basis learning, clustering and feature extraction for high-dimensional classification problems with applications in music analysis \cite{fevotte09NMFMusicAnalysis}, document clustering \cite{Ding05K-Means-Equivalence} and medical imaging problems such as tumor typing in matrix-assisted laser desorption/ionisation (MALDI) imaging in the field of bioinformatics \cite{leuschner18}.
            
            Different from the \ac{PCA} approach above, the \ac{NMF} enforces nonnegativity constraints on the factor matrices without any orthogonality restrictions. This makes the \ac{NMF} the method of choice for application fields, where the underlying physical model enforces the solution to be nonnegative assuming that each datapoint can be described as a superposition of some unknown characteristic features of the dataset. The \ac{NMF} makes it possible to extract these features while constraining the matrix factors to have nonnegative entries, which simplifies their interpretation.
            These data assumptions are true for many application fields including the ones mentioned above but also especially our considered problem of dynamic computed tomography, where the measurements consist naturally of the nonnegative absorption of photons.
            Mathematically, the basic \ac{NMF} problem can be formulated as follows: For a given nonnegative matrix $X\in \R_{\geq 0}^{N\times T},$ find nonnegative matrices $B\in \R_{\geq 0}^{N\times K}$ and $C\in \R_{\geq 0}^{K\times T}$ with $K\ll \min\{N, T\},$ such that
            \begin{equation*}
                X\approx BC.
            \end{equation*}
            The factorisation allows to approximate the rows $X_{n, \bullet}$ and columns $X_{\bullet, t}$ of the data matrix as a superposition of the $K$ columns $B_{\bullet, k}$ of $B$ and rows $C_{k, \bullet}$ of $C$ respectively, such that $X_{n,\bullet} \approx \sum_{k=1}^{K} B_{nk}C_{k,\bullet} $ and $ X_{\bullet, t} \approx \sum_{k=1}^{K} C_{kt} B_{\bullet, k}.$ Similarly, it holds that
            \begin{equation*}
                X \approx BC = \sum_{k=1}^{K} B_{\bullet, k} C_{k, \bullet},
            \end{equation*}
            where the $K$ terms of the sum are rank-one matrices. Hence, the sets $\{ B_{\bullet, k} \}_k$ and $ \{ C_{k, \bullet} \}_k $ can be interpreted as a low-dimensional basis to approximate the data matrix, i.e.\ the \ac{NMF} performs the task of basis learning with additional nonnegativity constraints.
            
            The usual approach to compute the factorisation is to define a suitable discrepancy term $\mathcal{D}_{\text{NMF}},$ which has to be chosen according to the noise assumption of the underlying problem, and to reformulate the \ac{NMF} as a minimisation problem. Typical discrepancies include the default case of the Frobenius Norm on which we will focus on, the Kullback-Leibler divergence, the Itakura-Saito distance or other generalized divergences \cite{cichocki09bookNMF}.
            
            Furthermore, \ac{NMF} problems are usually ill-posed due to the non-uniqueness of the solution \cite{Klingenberg09-NMFIllposedness} and require the application of suitable regularisation techniques. One common method is to include penalty terms in the minimisation problem to tackle the ill-posedness of the problem but also to enforce desirable properties of the factorisation matrices. Typical examples range from $\ell_1, \ell_2$ and total variation regularisation terms \cite{lecharlierDeMol13NMFTV} to more problem specific terms, which enforce additional orhogonality of the matrices or even allow supervised classification workflows if the \ac{NMF} is used as a prior feature exctraction method \cite{FM18,leuschner18}.
            
            Hence, the general regularised \ac{NMF} problem can be written as
            \begin{equation}\label{eq:NMF Minimisation Problem}
	            \min_{B, C\geq 0} \mathcal{D}_{\text{NMF}}(X, BC) + \sum_{\ell = 1}^{L} \gamma_\ell \mathcal{P}_\ell(B, C) \eqqcolon \min_{B, C\geq 0} \mathcal{F}(B, C),
            \end{equation}
		    where $\mathcal{P}_\ell$ denote the penalty terms, $\gamma_\ell\geq 0$ the corresponding regularisation parameters and $\mathcal{F}$ the cost function of the \ac{NMF}.
		    
\begin{figure}
            	\centering
            	\includegraphics[width=\textwidth]{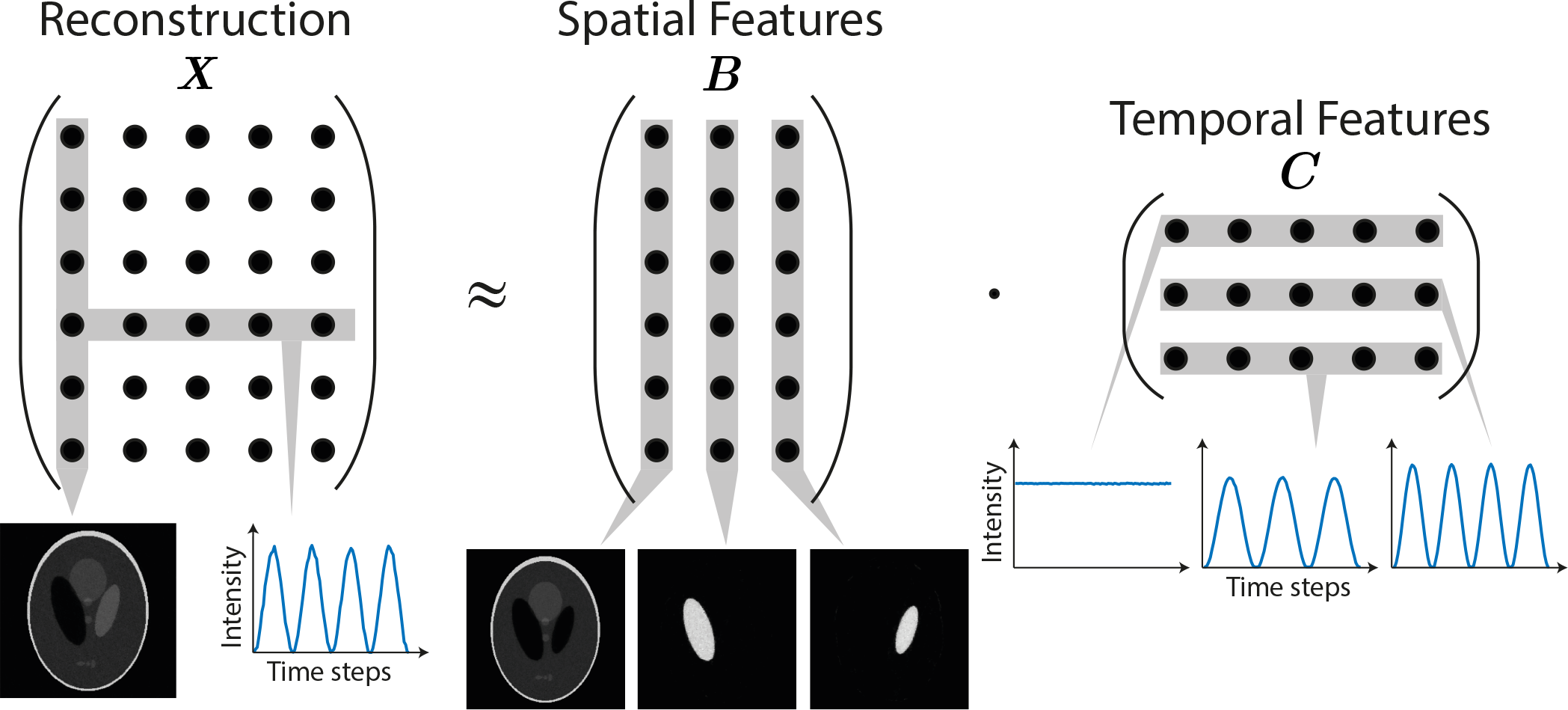}
            	\caption{Structure of the NMF in the context of the dynamic Shepp-Logan phantom as shown in Figure \ref{fig:dynSheppIllustration}. Here, the nonnegative spatial and temporal basis functions can be naturally represented by the matrices $B$ and $C$. }
            	\label{fig:NMFandCT}
            \end{figure}

		    The considered optimisation approach in this work is based on the so-called Majorise-Minimisation principle and gives rise to multiplicative update rules of the matrices in \eqref{eq:NMF Minimisation Problem}, which automatically preserve the nonnegativity of the iterates provided that they are initialised nonnegative. For more details on this optimisation technique, we refer the reader to Appendix \ref{app:sec:Optimisation Techniques for NMF Problems}. 
            The idea of the feature extraction procedure based on the NMF can be well illustrated by considering the example from Figure \ref{fig:dynSheppIllustration} that satisfies the decomposition assumption from \eqref{eqn:continuousDeompostion}. Here, the highlighted spatial regions change their intensities according to the given dynamics. The NMF allows a natural interpretation of the factorisation matrices $B$ and $C$ as the spatial and temporal basis functions for this case, as illustrated in Figure \ref{fig:NMFandCT}. The column $X_{\bullet, t}$ of $X$ denotes the reconstruction of the $t$-th time step of the inverse problem in \eqref{eq:dynamic IP General Discrete}. The NMF allows to decompose the spatial and temporal features of the data: The matrix $B$ contains the spatial features in its columns with the corresponding temporal features in the rows of $C.$ 
		    
    \subsection{Separated Reconstruction and Low-rank Decomposition} \label{subsec:Separated Reconstruction and Low--rank Decomposition}
        Let us first discuss a separated reconstruction and feature extraction approach to solve the inverse problem in \eqref{eq:dynamic IP General Discrete}, that means we compute first a reconstruction and perform then subsequently the feature extraction with one of the previously discussed methods. We consider this method as baseline for our comparison.
        
        The considered reconstruction method for this separated framework involves a basic gradient descent approach together with a regularisation step and a subsequent total variation denoising, which will be henceforth referred to as \texttt{gradTV}. The details on the algorithm are provided in Algorithm \ref{alg:gradTV}. 
        In particular, we aim to compute solutions to the least squares problem and incorporate the low-rank assumptions as additional penalty of the nuclear norm of $X_{\bullet, t}$, that is
         \[
         \min_{X_{\bullet, t}\geq 0} \| Y_{\bullet, t} - A_t X_{\bullet, t} \|^2_2 + \alpha \|X_{\bullet, t}\|_*
         \]
         for all $t$; see e.g.~\cite{Lingala2011,Tremoulheac2014}. This can then be efficiently 
         solved by a proximal gradient descent scheme with a soft-thresholding 
         on the singular values and hence enforcing the low-rank structure. 
         Ideally, one would like to include the total variation regularisation 
         as penalty term, but as this tends to be computationally expensive for 
         the fine temporal sampling, we included this as a subsequent denoiser.
         
        In practice, after a suitable initialisation of the reconstruction matrix, the gradient descent step is computed with an, \emph{a priori} defined, fixed stepsize $\rho_{\text{grad}}.$ For the proximal step, the truncated \ac{SVD} of $X$ is computed and a soft thresholding of the singular values is performed with a fixed threshold $\rho_{\text{thr}}.$ 
        Afterwards, we enforce the nonnegativity with a projection step on the reconstruction $X.$ When the stopping criterion is satisfied, a TV denoising algorithm\footnote{\tiny \url{https://www.mathworks.com/matlabcentral/fileexchange/36278-split-bregman-method-for-total-variation-denoising}} based on \cite{GoldsteinOsher-TVDenoiser,Tremoulheac14PhDThesis} with the corresponding parameter $\rho_{\text{TV}}$ is applied.
        \begin{algorithm}
            \caption{\texttt{gradTV}}
            \label{alg:gradTV}
            \begin{algorithmic}[1] 
                \State \textbf{Initialise:} $X$
                \State \textbf{Input:} $\rho_{\text{grad}}, \rho_{\text{thr}}, \rho_{\text{TV}} >0$
                \Repeat
                \State $X_{\bullet, t} \gets X_{\bullet, t} - \rho_{\text{grad}} (A_t^\intercal A_tX_{\bullet, t} - A_t^\intercal Y_{\bullet, t}) \ \ \ \text{for all}\ t$
                \State $(U,\Sigma, V) \gets \textsc{SVD}(X)$
                \State $\Sigma \gets \textsc{SoftThresh}_{\rho_{\text{thr}}}(\Sigma)$
                \State $X \gets U\Sigma V^\intercal$ 
                \State $X\gets \max(X, 0) $
                \Until{\textsc{StoppingCriterion} satisfied }
                \State $X \gets \textsc{TVDenoiser}_{\rho_{\text{TV}}}(X)$\\
                \Return $X$
            \end{algorithmic}
        \end{algorithm}
        
        After the reconstruction procedure given by Algorithm \ref{alg:gradTV}, we perform the feature extraction of the reconstruction $X$ via both the \ac{PCA} and the \ac{NMF} and call the approach \texttt{gradTV\_PCA} and \texttt{gradTV\_NMF} respectively.
        
        For \texttt{gradTV\_PCA}, we simply compute the \ac{PCA} of $X$ based on its \ac{SVD}. Concerning the method \texttt{gradTV\_NMF}, we consider the standard \ac{NMF} model
        \begin{equation}\label{eq:NMF model gradTV_NMF}
            \min_{B, C\geq 0} \Vert X - BC\Vert_F^2 + \dfrac{\tilde{\mu}_C}{2} \Vert C \Vert_F^2
        \end{equation}
        with the parameter $\tilde{\mu}_C.$ The $\ell_2$ regularisaton penalty term on $C$ is motivated by our application in Section \ref{sec:Application to Dynamic CT}. The corresponding multiplicative algorithms to solve \eqref{eq:NMF model gradTV_NMF} are well-known \cite{demol12,FM18} and a special case of the derived update rules in the next Section.
	\subsection{Joint Reconstruction and Low-rank Decomposition}\label{subsec:Joint Reconstruction and low--rank Decomposition}
        Instead of the previously discussed separated reconstruction, we now aim to include the feature extraction into the reconstruction procedure. This gives rise to consider a joint reconstruction and low-rank decomposition approach based on the \ac{NMF}, rather than one based on a low-rank plus sparsity approach based on \ac{PCA}~\cite{Candes2011,Tao2013,Yuan2013}. The basic idea of the method is to incorporate the reconstruction procedure of the inverse problem in \eqref{eq:dynamic IP General Discrete} into the NMF workflow. To do this, we have to additionally assume that $A_t\in \R_{\geq 0}^{M\times N}, Y\in \R_{\geq 0}^{M\times T}$ and $X\in \R_{\geq 0}^{N\times T}$ to ensure the desired nonnegativity of the factorisation matrices $B$ and $C,$ which corresponds to the assumptions of the decomposition in \eqref{eqn:continuousDeompostion}. The main motivation is that this joint approach allows the reconstruction process to exploit the underlying latent \ac{NMF} features of the dataset,  which can therefore enhance the quality of the reconstructions by enabling regularisation of temporal and spatial features separately.
        
	    This can be achieved by including a discrepancy term $\mathcal{D}_{\text{IP}}(Y_{\bullet, t}, A_t X_{\bullet, t})$ of the inverse problem into the NMF cost function in \eqref{eq:NMF Minimisation Problem}. This leads together with some possible penalty terms for the reconstruction $X$ to the model
	    \begin{equation} \label{eq:NMF model BC-X General}
	        \min_{B, C, X\geq 0} \mathcal{D}_{\text{IP}}(Y_{\bullet, t}, A_t X_{\bullet, t}) + \alpha \mathcal{D}_{\text{NMF}}(X, BC) + \sum_{\ell = 1}^{L} \gamma_\ell \mathcal{P}_\ell(B, C, X),
	    \end{equation}
	    with $\alpha\geq 0$ for the joint reconstruction and low-rank decomposition problem, which we will call \texttt{BC-X}.
	    Furthermore, we can enforce $X\coloneqq BC$ as a hard constraint, such that the reconstruction matrix will have at most rank $K.$ In this case, the discrepancy $\mathcal{D}_{\text{NMF}}$ vanishes and we end up with the model \texttt{BC}:
	    \begin{equation} \label{eq:NMF model BC General}
	    \min_{B, C\geq 0} \mathcal{D}_{\text{IP}}(Y_{\bullet, t}, A_t (BC)_{\bullet, t}) + \sum_{\ell = 1}^{L} \gamma_\ell \mathcal{P}_\ell(B, C).
	    \end{equation}
        \subsubsection{Considered \ac{NMF} Models}\label{subsubsec:Considered NMF Models}
    	    For both models \eqref{eq:NMF model BC-X General} and \eqref{eq:NMF model BC General}, we
    	    use the standard Frobenius norm for both the discrepancy terms $\mathcal{D}_{\text{NMF}}$ and $\mathcal{D}_{\text{IP}}.$
    	    Furthermore, the optimisation method discussed in Section \ref{subsubsec:Algorithms} allows to include a variety of penalty terms into the cost function. This makes it possible to construct suitable regularised \ac{NMF} models and to enforce additional properties to the matrices depending on the specific application. In this work, we will consider standard $\ell_1$ and $\ell_2$ regularisation terms on each matrix and an isotropic total variation penalty on the matrix $B.$ The latter is motivated by our considered application, which denoises the spatial features and thus also the reconstruction matrix. Hence, we will focus on the following NMF models in the remainder of this work:
    	    \begin{equation} \tag*{\texttt{BC-X}} \label{eq:NMF model BC-X} 
                \begin{aligned}
                    &\scalebox{0.9}{\text{$\displaystyle\min_{B, C, X\geq 0} \Bigg\{ \sum_{t=1}^{T} \frac{1}{2} \Vert A_t X_{\bullet, t} - Y_{\bullet, t} \Vert_2^2 + \frac{\alpha}{2} \Vert BC - X \Vert_F^2 + \lambda_{B} \Vert B \Vert_1 + \frac{\mu_B}{2} \Vert B\Vert_F^2 $ }} \\
                    &\hspace{10ex}\scalebox{0.9}{\text{$\displaystyle + \lambda_{C} \Vert C \Vert_1 + \frac{\mu_C}{2} \Vert C\Vert_F^2 + \lambda_X \Vert X \Vert_1 + \frac{\mu_X}{2} \Vert X\Vert_F^2 + \dfrac{\tau}{2} \TV(B)\Bigg\}, $ }}
                \end{aligned}
            \end{equation}
            \begin{equation} \tag*{\texttt{BC}} \label{eq:NMF model BC}
                \begin{aligned}
                    &\min_{B, C \geq 0} \Bigg\{ \sum_{t=1}^{T} \frac{1}{2} \Vert A_t (BC)_{\bullet, t} - Y_{\bullet, t} \Vert_2^2 \!+\! \lambda_{C} \Vert C \Vert_1 + \frac{\mu_{C}}{2} \Vert C\Vert_F^2 + \lambda_{B} \Vert B \Vert_1\\
                    &\hspace*{9ex} + \frac{\mu_{B}}{2} \Vert B\Vert_F^2 + \dfrac{\tau}{2} \TV(B) \Bigg\}.
                \end{aligned}
            \end{equation}
            The regularisation parameters $ \alpha, \lambda_{C}, \mu_{C}, \lambda_{B}, \mu_{B}, \lambda_{X}, \mu_{X}, \tau \geq 0, $ 
            are chosen \emph{a priori}. 
            Furthermore, $ \Vert \cdot \Vert_F $ denotes the Frobenius norm, $ \Vert M \Vert_1\coloneqq \sum_{ij} \vert M_{ij} \vert $ the 1-norm for matrices $ M $ and $ \TV(\cdot) $ is the following smoothed isotropic total variation \cite{defrise11TV,FM18,lecharlierDeMol13NMFTV}.
            \begin{definition}\label{def:TV}
            	The total variation of a matrix $ B\in \mathbb{R}^{N\times K} $ is defined as
            	\begin{equation*}
            	    \TV (B) \coloneqq \sum_{k=1}^K \sum_{n=1}^N \vert \nabla_{nk}B\vert \coloneqq \sum_{k=1}^K \sum_{n=1}^N \sqrt{\varepsilon_{\TV}^2 + \sum_{\ell\in \mathcal{N}_n} (B_{nk}-B_{\ell k})^2},
            	\end{equation*}
            	where $ \varepsilon_{\TV}>0 $ is a small positive constant and $ \mathcal{N}_n $ are index sets referring to spatially neighboring pixels.
            \end{definition}
            A typical example for the neighbourhood of the pixel $ (0,0) $ in two dimensions is $ \mathcal{N}_{(0,0)} = \{ (1,0), (0,1) \} $ to get an estimate of the gradient components in both directions of the axes. The parameter $ \varepsilon_{\TV} $ ensures the differentiability of the TV penalty term.
            
            In the following section, we will present the multiplicative update rules for the NMF models \ref{eq:NMF model BC-X} and \ref{eq:NMF model BC} and derive the algorithms in Appendix \ref{app:sec:Derivation of the Algorithms} based on the Majorise-Minimisation principle.
        \subsubsection{Algorithms} \label{subsubsec:Algorithms}
            In this section, we present in Theorem \ref{theorem:Algorithm BC-X} and \ref{theorem:Algorithm BC} the multiplicative algorithms for the NMF problems in \ref{eq:NMF model BC-X} and \ref{eq:NMF model BC}. As mentioned in Section \ref{subsubsec:Nonnegative Matrix Factorisation}, the multiplicative structure of the iteration scheme ensures automatically the nonnegativity of the matrices $B$ and $C$ as long as they are initialised nonnegative. The derivation of such algorithms in this work are based on the so-called Majorise-Minimisation principle. The main idea of this approach is to replace the considered NMF cost function $\mathcal{F}$ with a suitable auxiliary function $\mathcal{Q}_{\mathcal{F}}$, whose minimisation is much easier to handle and leads to a monotone decrease of $\mathcal{F}.$ Furthermore, specific construction techniques of these \textit{surrogate functions} lead to the desired multiplicative update rules which fulfill the nonnegativity constraint. We provide a short description of the main principles in Appendix \ref{app:sec:Optimisation Techniques for NMF Problems}. A more detailed discussion of different construction methods for various kinds of discrepancy and penalty terms of $\mathcal{F}$ can be found in the survey paper \cite{FM18}.
            
            For better readability, we present only the main results here and a detailed construction of the surrogate functions as well as derivation of the algorithms for both cost functions \ref{eq:NMF model BC-X} and \ref{eq:NMF model BC} can be found in Appendix \ref{app:sec:Derivation of the Algorithms}. 
            Consequently, we will only state the main results in Theorem \ref{theorem:Algorithm BC-X} and \ref{theorem:Algorithm BC} here. Nevertheless, due to the construction of a suitable surrogate function for the TV penalty term (see Appendix \ref{app:sec:Derivation of the Algorithms} and \cite{FM18} for more details), we first introduce the following matrices $ P(B), Z(B)\in \mathbb{R}_{\geq 0}^{N\times K} $ as
            \begin{align}
            	P(B)_{n k} &\coloneqq  \dfrac{1}{\vert \nabla_{nk} B \vert} \sum_{\ell \in \mathcal{N}_n} 1 + \sum_{\ell \in \bar{\mathcal{N}}_n} \dfrac{1}{\vert \nabla_{\ell k} B \vert}, \label{eq:TV:MatrixP}\\
        		Z(B)_{n k} &\coloneqq  \dfrac{1}{P(B)_{n k}} \left ( \dfrac{1}{\vert \nabla_{nk} B \vert} \sum_{\ell \in {\mathcal{N}}_n} \dfrac{B_{n k} + B_{\ell  k}}{2} + \sum_{\ell \in \bar{\mathcal{N}}_n} \dfrac{B_{n k} + B_{\ell  k}}{2 \vert \nabla_{\ell k} B \vert} \right ), \label{eq:TV:MatrixZ}
            \end{align}
            where $ \bar{\mathcal{N}}_n $ is the set of the so-called \textit{adjoint} neighbourhood pixels, which is given by the relation
    		\begin{equation*}
    	    	\ell\in \bar{\mathcal{N}}_n \Leftrightarrow n\in \mathcal{N}_\ell.
    		\end{equation*}
            Furthermore, we write $\vec{1}_{M\times N}$ for an $M\times N$ matrix with ones in every entry.
            
            We then obtain the two algorithms for both models under consideration. First for the \ref{eq:NMF model BC-X} model that jointly obtains the reconstruction $X$ and the decomposition:
                \begin{theorem}[Algorithm for \ref{eq:NMF model BC-X}] \label{theorem:Algorithm BC-X}
                	For $A_t\in \R_{\geq 0}^{M\times N}, Y\in \R_{\geq 0}^{M\times T}$ and initialisations $ X^{[0]}\in \mathbb{R}_{> 0}^{N\times T}, B^{[0]}\in \mathbb{R}_{> 0}^{N\times K}, C^{[0]}\in \mathbb{R}_{> 0}^{K\times T}, $ the alternating update rules
                	\begin{align*}
                	    X_{\bullet, t}^{[d+1]} &= X_{\bullet, t}^{[d]} \circ \dfrac{A_t^\intercal Y_{\bullet, t} + \alpha B^{[d]}C_{\bullet, t}^{[d]}}{A_t^\intercal A_tX_{\bullet, t}^{[d]} + (\mu_{X} + \alpha)X_{\bullet, t}^{[d]} + \lambda_{X} \vec{1}_{N\times 1}} \\
                	    B^{[d+1]} &= B^{[d]} \circ \dfrac{\alpha X^{[d+1]}{C^{[d]}}^\intercal + \tau P(B^{[d]}) \circ Z(B^{[d]})}{\alpha B^{[d]}C^{[d]}{C^{[d]}}^\intercal + \mu_{B} B^{[d]} + \lambda_{B} \vec{1}_{N\times K} + \tau B^{[d]} \circ P(B^{[d]}) } \\
                	    C^{[d+1]} &= C^{[d]} \circ \dfrac{\alpha {B^{[d+1]}}^\intercal X^{[d+1]}}{\alpha {B^{[d+1]}}^\intercal B^{[d+1]} C^{[d]} + \mu_{C} C^{[d]} + \lambda_{C} \vec{1}_{K\times T}}
                	\end{align*}
                	lead to a monotonic decrease of the cost function in \ref{eq:NMF model BC-X}.
                \end{theorem}
                Similarly, for the \ref{eq:NMF model BC} model we obtain the updates rules without constructing the matrix $X$ during the reconstruction process:
                \begin{theorem}[Algorithm for \ref{eq:NMF model BC}] \label{theorem:Algorithm BC}
                	For $A_t\in \R_{\geq 0}^{M\times N}, Y\in \R_{\geq 0}^{M\times T}$ and initialisations $ B^{[0]}\in \mathbb{R}_{> 0}^{N\times K}, C^{[0]}\in \mathbb{R}_{> 0}^{K\times T}, $ the alternating update rules
                	\begin{align*}
                	    \scalebox{0.82}{\text{$B^{[d+1]}$ }} &=\scalebox{0.82}{\text{$ B^{[d]} \circ \dfrac{\sum_{t=1}^T A_t^\intercal Y_{\bullet, t} \cdot ({C^{[d]}}^\intercal)_{t,\bullet} + \tau P(B^{[d]}) \circ Z(B^{[d]})}{\sum_{t=1}^T A_t^\intercal A_t (B^{[d]} C^{[d]})_{\bullet, t} \cdot ({C^{[d]}}^\intercal)_{t,\bullet} + \mu_{B} B^{[d]} + \lambda_{B}\vec{1}_{N\times K} + \tau B^{[d]} \circ P(B^{[d]})}$ }} \\
                	    C^{[d+1]}_{\bullet,t} &= C^{[d]}_{\bullet,t} \circ \dfrac{{B^{[d+1]}}^\intercal A_t^\intercal Y_{\bullet,t}}{{B^{[d+1]}}^\intercal A_t^\intercal A_t (B^{[d+1]}C^{[d]})_{\bullet,t} + \mu_{C} C^{[d]}_{\bullet,t} + \lambda_{C} \vec{1}_{K\times 1} }
                	\end{align*}
                	lead to a monotonic decrease of the cost function in \ref{eq:NMF model BC}.
                \end{theorem}
                We remind that the derivation is described in Appendix \ref{app:sec:Derivation of the Algorithms}, which leads to the update rules in the Theorems above.
                Due to the multiplicative structure of the algorithms, zero entries in the matrices stay zero during the iteration scheme and can cause divisions by zero. This issue is handled via the strict positive initialisation in both Theorems. Furthermore, very small or high numbers can cause numerical instabilities and lead to undesirable results. As a standard procedure, this problem is handled by suitable projection steps after every iteration step \cite{cichocki09bookNMF}.

    \subsection{Complexity Reduction for Stationary Operator} \label{subsec:Complexity Reduction for Stationary Operator}
        Let us now consider the case of a stationary operator, i.e. $\mathcal{A}(\cdot;t)$ in equation \eqref{eq:dynamic IP General Continuous} does not change with $t$. Then we simply write $\mathcal{A}$ or $A$ for the matrix representation in \eqref{eq:dynamic IP General Discrete}. If further the number of channels $T$ is large, the application of the forward operator represented a major computational burden per channel. 
        In particular, we make use here of the assumption $T\gg K$, i.e. the number of channels is much larger than the basis functions for the decomposition. In this case, we can effectively reduce the computational cost by shifting the application of the forward operator to the spatial basis functions contained in $B$. That means, we make essential use of the decomposition $X\approx BC$ in the reconstruction task and as such avoid to construct the approximation to $X$. Consequently, we will only consider the case of \ref{eq:NMF model BC} here. Since $A$ is independent from $t,$ the NMF model \ref{eq:NMF model BC} becomes
        
        \begin{equation} \tag*{\texttt{sBC}} \label{eq:NMF model sBC}
            \begin{aligned}
                &\min_{B, C \geq 0} \Big\{ \frac{1}{2} \Vert ABC - Y \Vert_F^2 + \lambda_{C} \Vert C \Vert_1 + \frac{\mu_{C}}{2} \Vert C\Vert_F^2 + \lambda_{B} \Vert B \Vert_1\\
                &\hspace*{9ex} + \frac{\mu_{B}}{2} \Vert B\Vert_F^2 + \dfrac{\tau}{2} \TV(B) \Big\}.
            \end{aligned}
        \end{equation}
        
        To illustrate this, let us consider the update equation in Theorem \ref{theorem:Algorithm BC} for $B,$ where we can simplify the first term in the denominator as follows:
        \begin{equation*}
            \scalebox{0.86}{\text{$\sum_{t=1}^T A^\intercal A (B^{[d]} C^{[d]})_{\bullet, t} \cdot ({C^{[d]}}^\intercal)_{t,\bullet} = A^\intercal A \sum_{t=1}^T (B^{[d]} C^{[d]})_{\bullet, t} \cdot ({C^{[d]}}^\intercal)_{t,\bullet} = A^\intercal A B^{[d]} C^{[d]} {C^{[d]}}^\intercal.$ }}
        \end{equation*}
        The other terms in the update rules can be simplified similarly, such that we obtain the following reduced update equations:
        
        \begin{corollary}[Algorithm for \ref{eq:NMF model sBC}] \label{cor:Algorithm sBC}
            For $A\in \R_{\geq 0}^{M\times N}, Y\in \R_{\geq 0}^{M\times T}$ and initialisations $ B^{[0]}\in \mathbb{R}_{> 0}^{N\times K},$ $C^{[0]}\in \mathbb{R}_{> 0}^{K\times T}, $ the alternating update rules
        	\begin{align*}
        	    B^{[d+1]} &= B^{[d]} \circ \dfrac{A^\intercal Y {C^{[d]}}^\intercal + \tau P(B^{[d]}) \circ Z(B^{[d]})}{A^\intercal A B^{[d]} C^{[d]} {C^{[d]}}^\intercal + \mu_{B} B^{[d]} + \lambda_{B}\vec{1}_{N\times K} + \tau B^{[d]} \circ P(B^{[d]})} \\
        	    C^{[d+1]} &= C^{[d]} \circ \dfrac{{B^{[d+1]}}^\intercal A^\intercal Y}{{B^{[d+1]}}^\intercal A^\intercal A B^{[d+1]}C^{[d]} + \mu_{C} C^{[d]} + \lambda_{C} \vec{1}_{K\times T} }.
	        \end{align*}
        	lead to a monotonic decrease of the cost function in \ref{eq:NMF model sBC}.
        \end{corollary}
    	
    	Finally, the order of application is essential here to obtain the complexity reduction. 
    	In particular, we implemented the algorithm such that $A$ acts on the basis functions in $B$. That means, we compute first the product $A^\intercal A B$ followed by multiplication with $C$. That means, we can expect a reduction of computational complexity by a factor $T/K$ with the \ref{eq:NMF model sBC} model and hence is especially useful for dimension reduction under fine temporal sampling.

\section{Application to Dynamic CT} \label{sec:Application to Dynamic CT}



In the following we will apply the presented methods to the use case of dynamic computerised tomography (CT).  Here, the quantity of interest is given as the attenuation coefficient $x(s,t)$ at time $t\in [0,T]$ on a bounded domain in two dimensions $s\in\Omega_1\subset\R^2$. Following the formulation in \eqref{eq:dynamic IP General Continuous}, the time-dependent forward operator is given by the Radon transform 
    	\begin{equation}\label{eq:DCT Continuous}
    		y(\theta, \sigma, t) \coloneqq (\mathcal{R}_{\mathcal{I}(t)} x(s, t))(\theta, \sigma) = \int_{s\cdot \theta = \sigma} x(s, t) \diff s 
    	\end{equation}
	Here, the measurement $y(\theta,\sigma,t)$ consist of line integrals over the domain $\Omega_1$ for each time point $t\in\mathcal{T}$, and is referred to as the sinogram. This measurement depends on two parameters, the angle $\theta \in S^1$ on the unit circle and a signed distance to the origin $\sigma\in\R$. Consequently, the measurements depend on a set of angles at each time step $\mathcal{I}(t)$, such that $(\theta,\sigma)\in \mathcal{I}(t)$ at time $t$, we will refer to this as the sampling patterns. In a slight abuse of notation, we will use $| \mathcal{I}(t) |$ for the number of angles, i.e. directions for the line integrals, at each time point.
	
	In the following we consider two scenarios for the choice of angles in $\mathcal{I}(t)$ and by that defining the nature of the forward operator, as discussed in Section \ref{subsec:Multi--Channel Inverse Problem}. In the general case of a nonstationary forward operator, that means the sampling patterns are time-dependent, we assume that the angles change but the amount of angles is constant over time $| \mathcal{I}(t) | \equiv c$. Additionally, we will consider the case for stationary operators, which in our setting means that the set of angles does not change over time, we can write for instance $\mathcal{I}(t) \equiv \mathcal{I}(t=0)$, and hence this leads to a stationary measurement operator of the dynamic process in \eqref{eq:DCT Continuous}. We note that even though the measurement process is stationary, the obtained measurement $y(\theta,\sigma,t)$ itself is still time dependent. 
	
	For the computations, we discretise \eqref{eq:DCT Continuous} to obtain a matrix vector representation as in \eqref{eq:dynamic IP General Discrete}. 	In the following we will write $R_t$ for the discrete Radon transform with respect to the sampling pattern $\mathcal{I}(t)$ at time point $t$, which gives rise to the discrete reconstruction problem for dynamic CT
    \begin{equation}\label{eq:DCT Discrete}
	        R_t X_{\bullet, t} = Y_{\bullet, t} \quad \text{for} \quad 1\leq t\leq T.
	\end{equation}
	We note, that due to the definition of the Radon transform by line integrals, the matrix $R_t\in \R_{\geq 0}^{M\times N}$ has only nonnegative entries and hence satisfies the assumption for Theorem \ref{theorem:Algorithm BC-X} and \ref{theorem:Algorithm BC}. Furthermore, $N$ denotes here the number of pixels in the original image and $M$ is given by the product $M\coloneqq \vert\mathcal{I}(t)\vert n_S,$ where $n_S$ is the number of detection points. 



    \subsection{Results and Discussion} \label{subsec:Results and Discussion}
        For a qualitative evaluation of the proposed \ac{NMF} approaches, we consider in the following sections two simulated datasets. Due to the known ground truth in both cases, we are able to measure the performance of each method via computing the mean of the \acf{PSNR} and the mean of the  \acf{SSIM} index \cite{BVW12-SSIM} over all time steps for every experiment.
        
        For each dataset, the parameters of all methods are chosen empirically to provide good reconstructions. For the \ac{NMF} models of the joint reconstruction and low-rank decomposition approach, we restrict ourselves to the total variation penalty term on $B$ to provide some denoising effect on the spatial features and the $\ell^2$ penalty on $C$ for the time features, since we expect and enforce smooth changes in time. We consider the standard case for the TV term with the default pixel neighbourhood and choose the smoothing parameter $\varepsilon_{\text{TV}} = 10^{-5}$ relatively small.
        
        Furthermore, for both datasets we measure different angles at each time step based on a tiny golden angle sampling \cite{wundrak2014small} using consecutive projections with increasing angle of $\varphi = 32.039\dots$, such that projection angles are not repeated. Nevertheless, we remind that we keep the total number of observed angles constant for each time step.
        
        For all considered approaches we use the unfiltered backprojection, given by the adjoint of the Radon transform,
        applied to the noisy data matrix $Y$ as the initialisation for the reconstruction matrix $X.$ In case of the NMF approaches, the matrices $B$ and $C$ are initialised via \ac{SVD} of $X$ based on \cite{BG08-NMFInit}. After the initialisation and at every iteration of the NMF algorithm, a suitable projection step for small values is performed to prevent numerical instabilities and zero entries during the multiplicative algorithm \cite{cichocki09bookNMF}.
        
        The algorithms were implemented with MATLAB\textsuperscript{\textregistered} R2019b and run on an Intel\textsuperscript{\textregistered} Core\texttrademark \ i7-7700K quad core CPU @4.20 GHz with 32 GB of RAM.
        
        To this end we present a summary and short explanation of all considered algorithms in this experimental section in Table \ref{tab:listOfAlgorithms}.
        \begin{table}[H]
            \ra{1.1}
            \resizebox{\columnwidth}{!}{%
            \pgfplotstabletypeset[
                col sep=&, row sep=\\,
                string type,
                every head row/.style={
                    before row={
                        \toprule
                    },
                    after row=\midrule,
                },
                every last row/.style={
                    after row=\bottomrule},
                every even row/.style={
                    before row={\rowcolor[gray]{0.9}}
                },
                columns/para/.style     ={column name=\textbf{Algorithm},column type=l},
                columns/descr/.style    ={column name=\textbf{Description},column type=l},
                ]{
                    para & descr \\
                    \ref{eq:NMF model BC} & \makecell[{{p{12cm}}}]{Joint reconstruction and feature extraction with the NMF model \ref{eq:NMF model BC} without constructing $X$, see algorithm in Theorem \ref{theorem:Algorithm BC}}\\
        			\ref{eq:NMF model BC-X} & \makecell[{{p{12cm}}}]{Joint reconstruction and feature extraction with NMF model \ref{eq:NMF model BC-X} and explicit construction of $X$, see algorithm in Theorem \ref{theorem:Algorithm BC-X}} \\
        			\ref{eq:NMF model sBC} & \makecell[{{p{12cm}}}]{Joint reconstruction and feature extraction method with NMF model \ref{eq:NMF model sBC} for \\ stationary operator, see algorithm in Corollary \ref{cor:Algorithm sBC}} \\
        			\texttt{gradTV} & \makecell[l]{Low-rank based reconstruction method for $X$, see Algorithm \ref{alg:gradTV}} \\
        			\texttt{gradTV\_PCA} & \makecell[{{p{12cm}}}]{Separated reconstruction and feature extraction with Algorithm \ref{alg:gradTV} and \\ subsequent PCA computation} \\
        			\texttt{gradTV\_NMF} & \makecell[{{p{12cm}}}]{Separated reconstruction and feature extraction with Algorithm \ref{alg:gradTV} and \\ subsequent NMF computation based on the model in \eqref{eq:NMF model gradTV_NMF}} \\
                }
            }
    		\caption{Summary and short explanation of considered algorithms in the experimental section.}
    		\label{tab:listOfAlgorithms}
        \end{table}
    
        
        \subsubsection{Shepp-Logan Phantom} \label{subsubsec:Shepp--Logan Phantom}
            This synthetic dataset consists of a dynamic two-dimensional Shepp-Logan phantom with $T=100$ and spatial size $ 128 \times 128$, see Figure \ref{fig:dynSheppIllustration} for the ground-truth.  During the whole time, two of the inner ellipsoids change their intensities sinusoidally with different frequencies while the rest of the phantom remains constant.
            
            In the following, we perform a variety of experiments for $\vert \mathcal{I}_t\vert \in \{ 2,\dots, 12\}$ with 1\% and 3\% Gaussian noise respectively. For all cases, we choose $K=5$ for the number of the NMF features. In such a way, the NMF is also capable to approximate minor characteristics such as noise or other artefacts of the reconstruction matrix besides the three main features.
            
            The parameters of all methods were determined empirically and are displayed in Table \ref{tab:dynShepp:parameterChoice} in Appendix \ref{app:sec:parameter Choice} for both noise levels. The stopping criterion for all methods is met, if 1200 iteration steps are reached or if the relative change of all matrices $B, C$ and $X$ goes below $5\cdot 10^{-5}.$
            \begin{figure}
                \centering
                \subfloat[\ref{eq:NMF model BC}]{\includegraphics[width=0.45\textwidth]{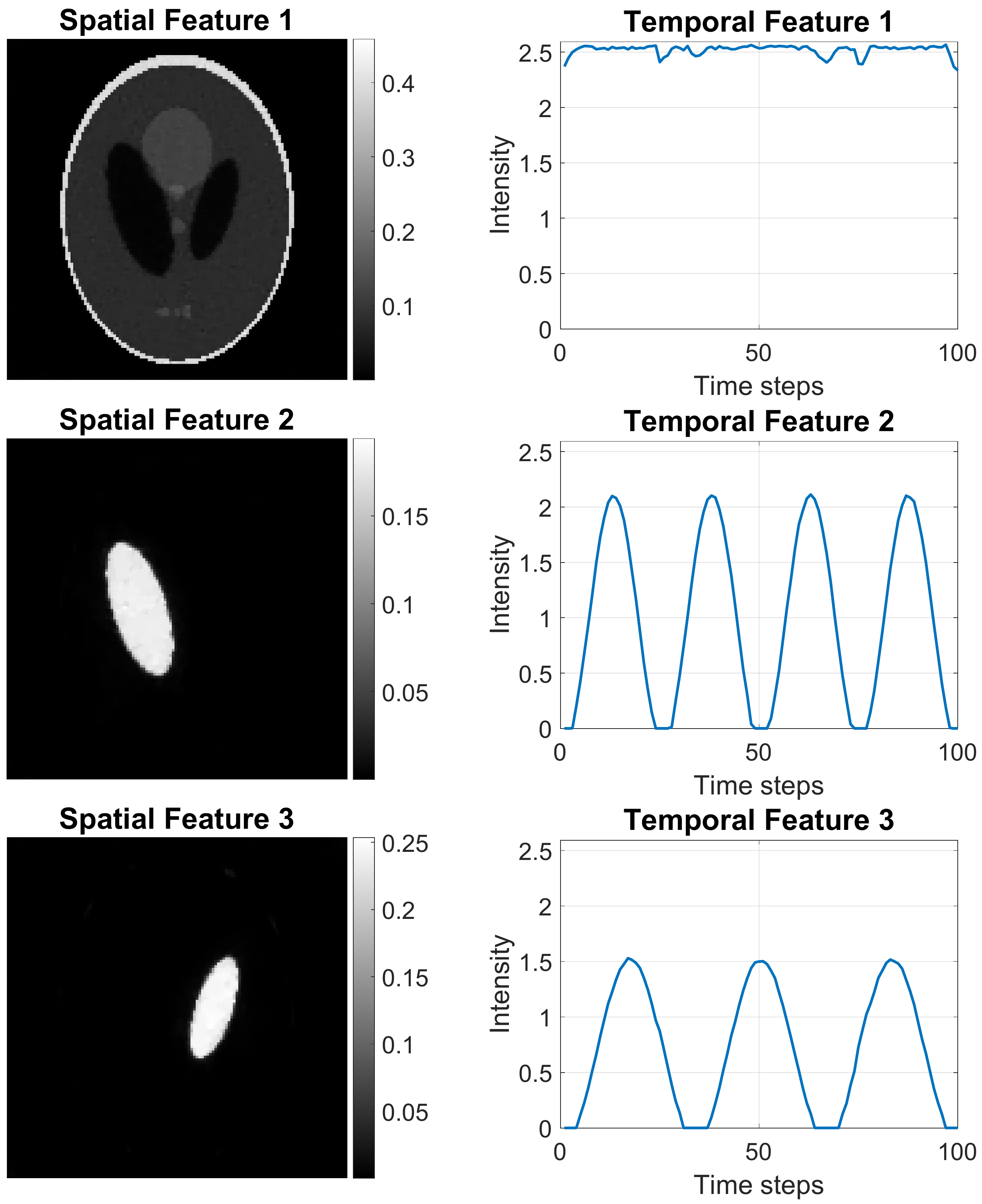}%
                \label{fig:dynShepp_BC_0.01_nTheta6}}
                \hfill
                \subfloat[\ref{eq:NMF model BC-X}]{\includegraphics[width=0.45\textwidth]{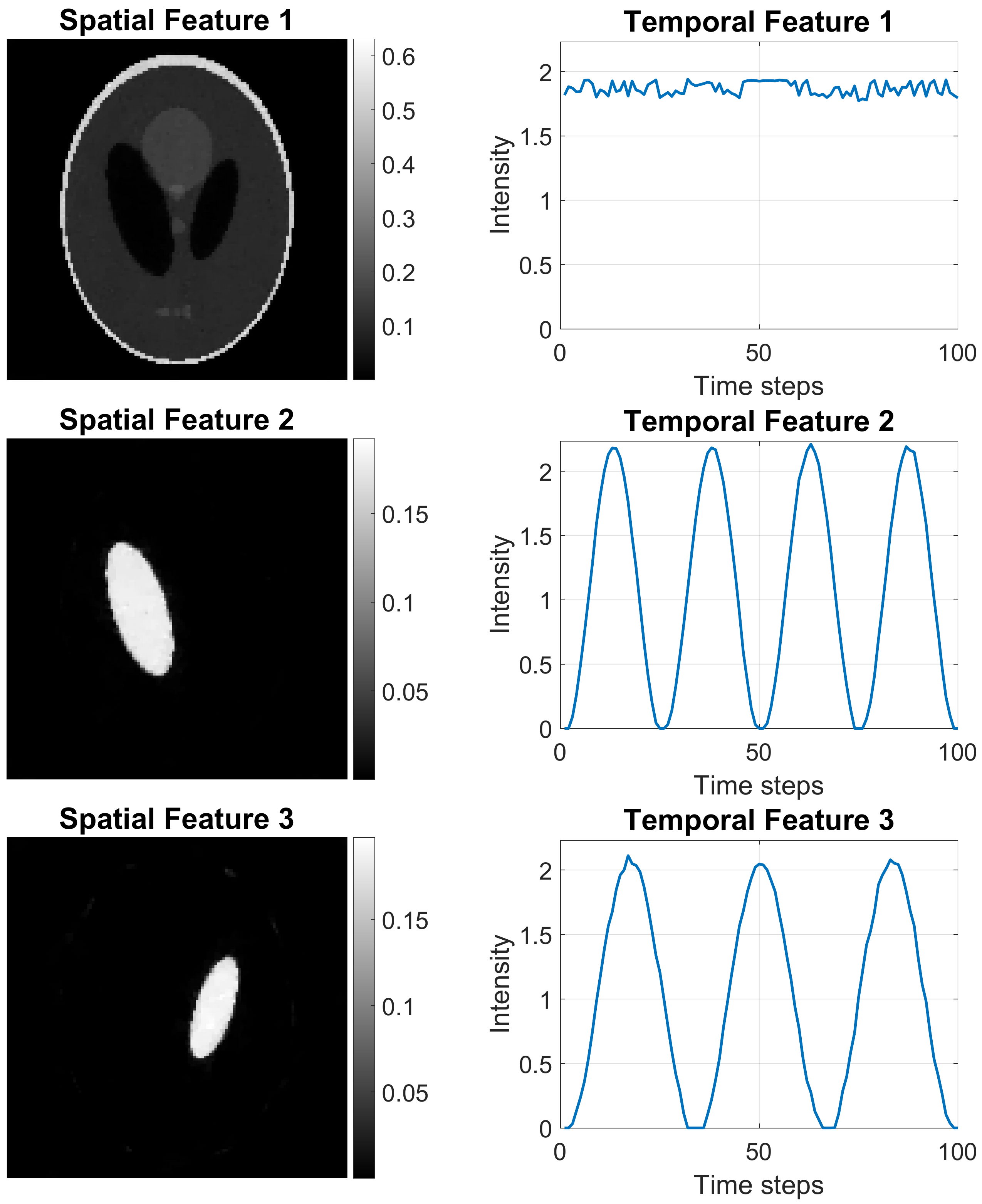}%
                \label{fig:dynShepp_BC-X_0.01_nTheta6}}
                \caption{Results for the dynamic Shepp-Logan phantom with $\vert \mathcal{I}_t \vert =6$ angles per time step and 1\% Gaussian noise. Shown are the leading extracted features for the \ref{eq:NMF model BC} model (left) and  for \ref{eq:NMF model BC-X} (right).}
                \label{fig:dynShepp:joint:noise1:nTheta6}
            \end{figure}
            \begin{figure}
                \centering
                \subfloat[\texttt{gradTV\_PCA}]{\includegraphics[width=0.45\textwidth]{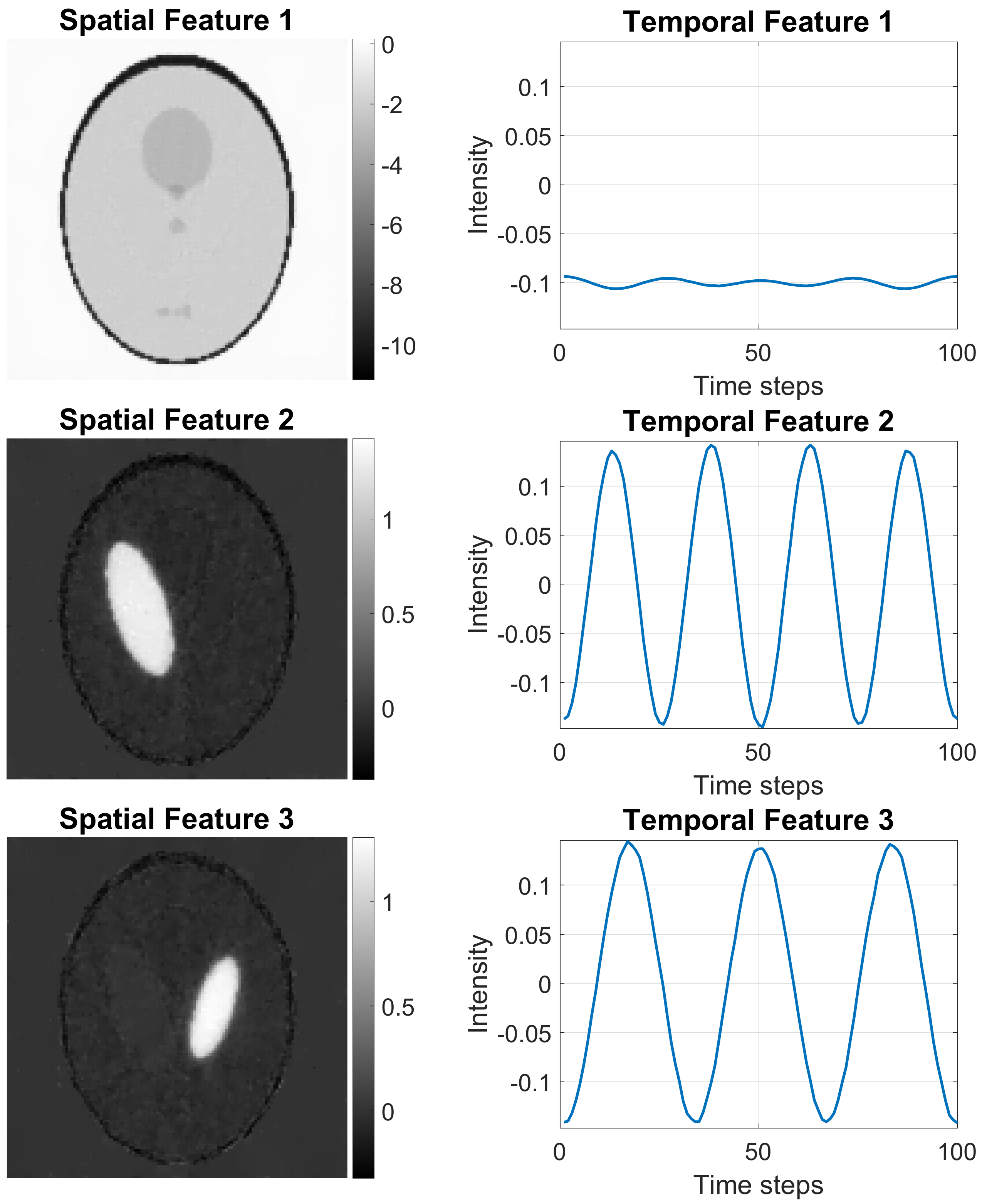}%
                \label{fig:dynShepp_gradTV_PCA_0.01_nTheta6}}
                \hfill
                \subfloat[\texttt{gradTV\_NMF}]{\includegraphics[width=0.45\textwidth]{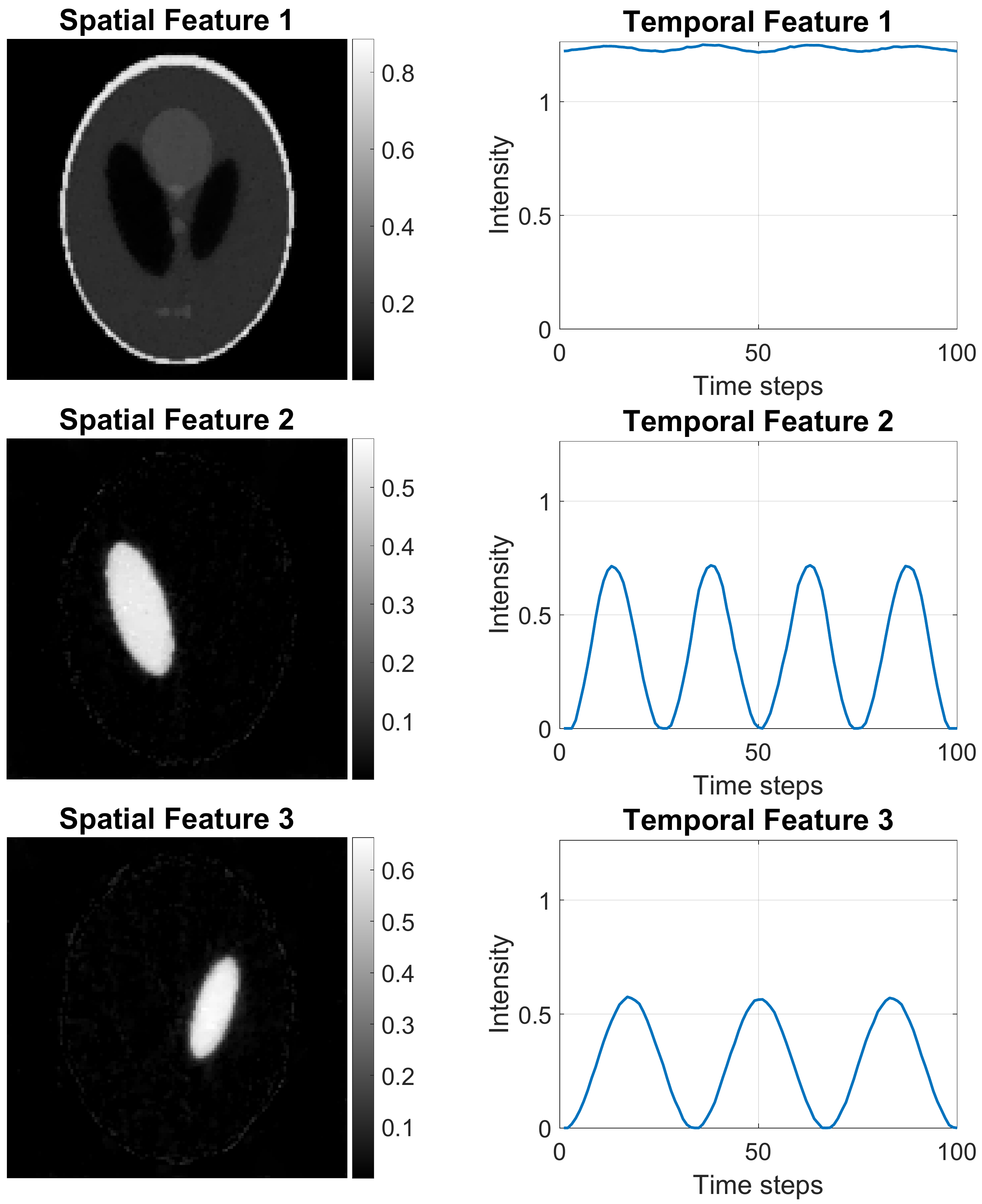}%
                \label{fig:dynShepp_gradTV_NMF_0.01_nTheta6}}
                \caption{Results for the dynamic Shepp-Logan phantom with $\vert \mathcal{I}_t \vert =6$ angles per time step and 1\% Gaussian noise. Shown are the leading extracted features for the \texttt{gradTV\_PCA} model (left) and  for \texttt{gradTV\_NMF} (right).}
                \label{fig:dynShepp:sep:noise1:nTheta6}
            \end{figure}
   
            
            We show first some results for the case with  $\vert \mathcal{I}_t \vert =6$ and 1\% Gaussian noise in Figure \ref{fig:dynShepp:joint:noise1:nTheta6} for the joint NMF methods and Figure \ref{fig:dynShepp:sep:noise1:nTheta6} for the separate reconstruction and extraction. The order of shown features is based on the singular values of $B$ for \texttt{gradTV\_PCA} and on the $\ell_2$-norm of the spatial features for NMF approaches.
            
            In this case, all considered approaches are able to successfully identify the constant and dynamic parts of the dataset and extract meaningful spatial and temporal features. The extracted spatial features of \ref{eq:NMF model BC}, \ref{eq:NMF model BC-X} and \texttt{gradTV\_NMF} show very clearly the dynamic and non-dynamic parts of the Shepp-Logan phantom. However, the spatial features of \texttt{gradTV\_NMF} are slightly more blurred and affected by minor artefacts especially in both dynamic features. This underlines the positive effect of the separate TV regularisation on the spatial feature matrix $B$ in the joint methods. In contrast, \texttt{gradTV\_PCA} is able to identify the main components of the dataset correctly, but there is a clear corruption of the dynamic features with other parts from the phantom. Furthermore, all spatial features contain negative parts due to the non-existent nonnegativity constraint of the \texttt{gradTV\_PCA} approach which makes their interpretation more challenging. Hence, the additional nonnegativity constraint of the NMF methods improve significantly the quality and interpretability of the extracted components in comparison with the \ac{PCA} based extraction method.
            
            The temporal features of all methods are clearly extracted and are consistent with the underlying ground truth of the dataset. However, we note that \ref{eq:NMF model BC} and \ref{eq:NMF model BC-X} have a slight difficulty to resolve the lower intensity part close to 0, which is probably caused by the multiplicative structure of the algorithms.

                \begin{figure}
                \centering
                \subfloat[\ref{eq:NMF model BC}]{\includegraphics[width=0.45\textwidth]{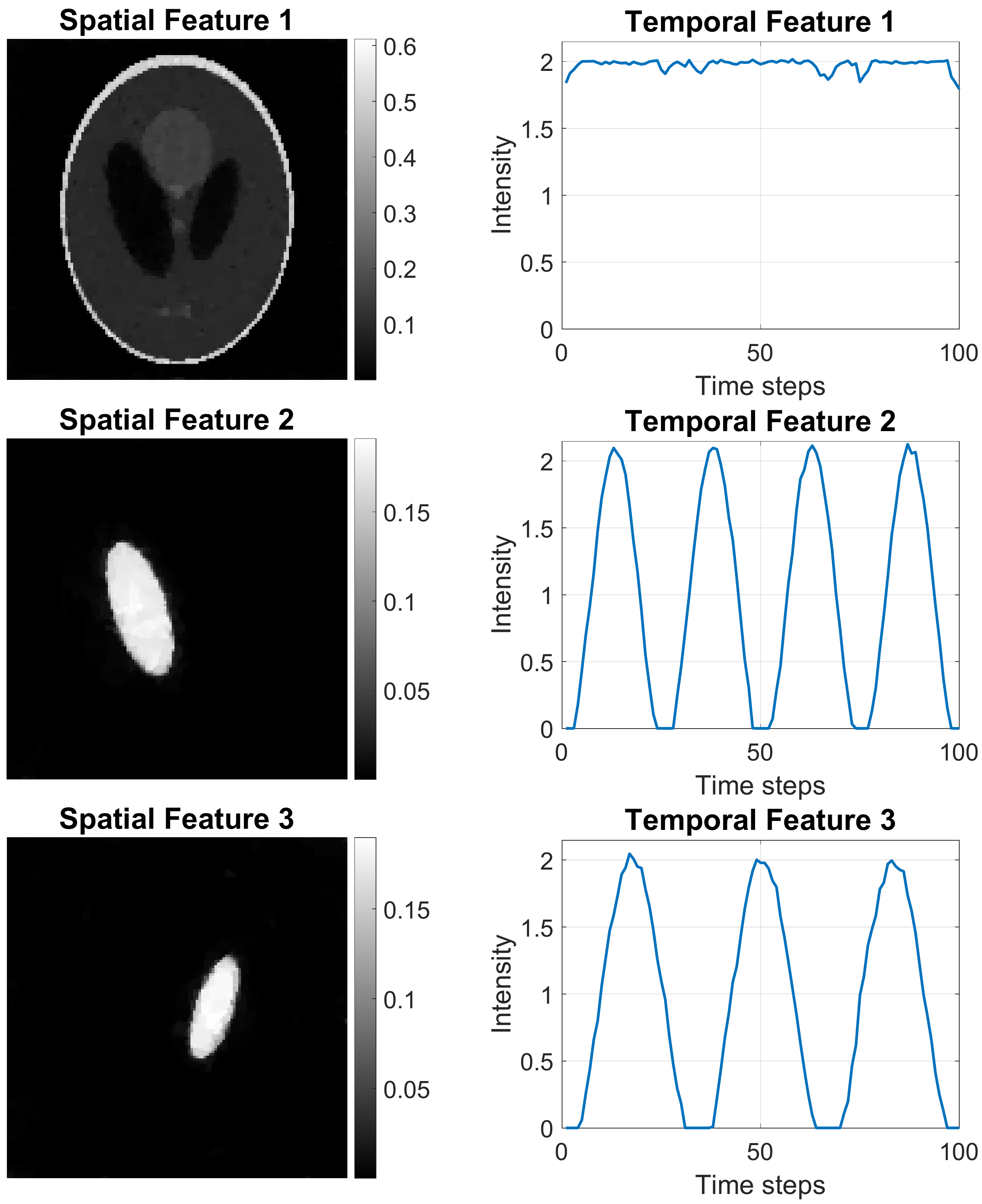}%
                \label{fig:dynShepp_BC_0.03_nTheta6}}
                \hfill
                \subfloat[\texttt{gradTV\_PCA}]{\includegraphics[width=0.45\textwidth]{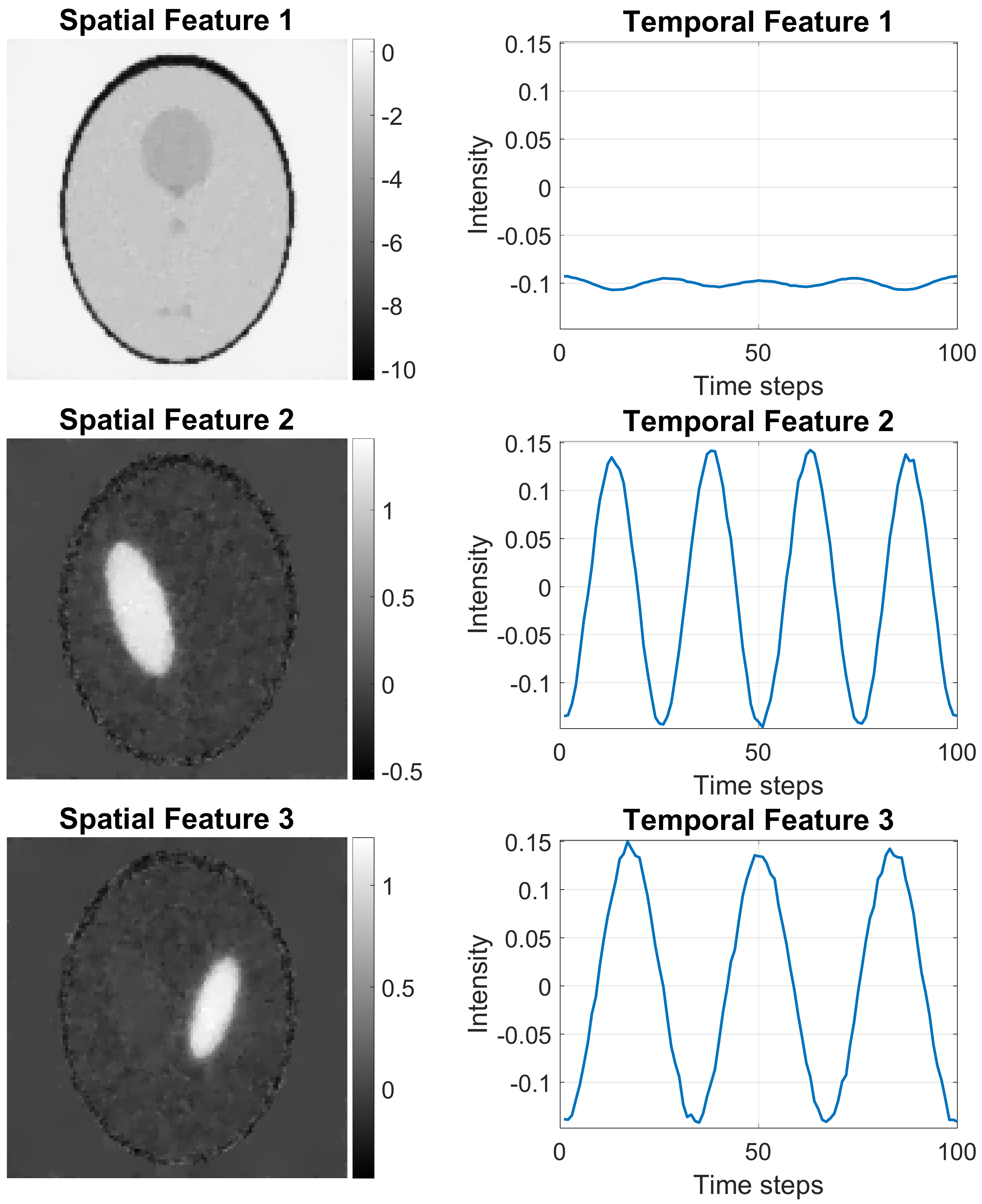}%
                \label{fig:dynShepp_gradTV_PCA_0.03_nTheta6}}
                \caption{Results for the dynamic Shepp-Logan phantom with $\vert \mathcal{I}_t \vert =6$ angles per time step and 3\% Gaussian noise. Shown are the leading extracted features for the \ref{eq:NMF model BC} model (left) and  for \texttt{gradTV\_PCA} (right).}
                \label{fig:dynShepp:noise3:nTheta6}
            \end{figure}
            \begin{figure}
                \centering
                \subfloat[\ref{eq:NMF model BC}]{\includegraphics[width=0.4\textwidth]{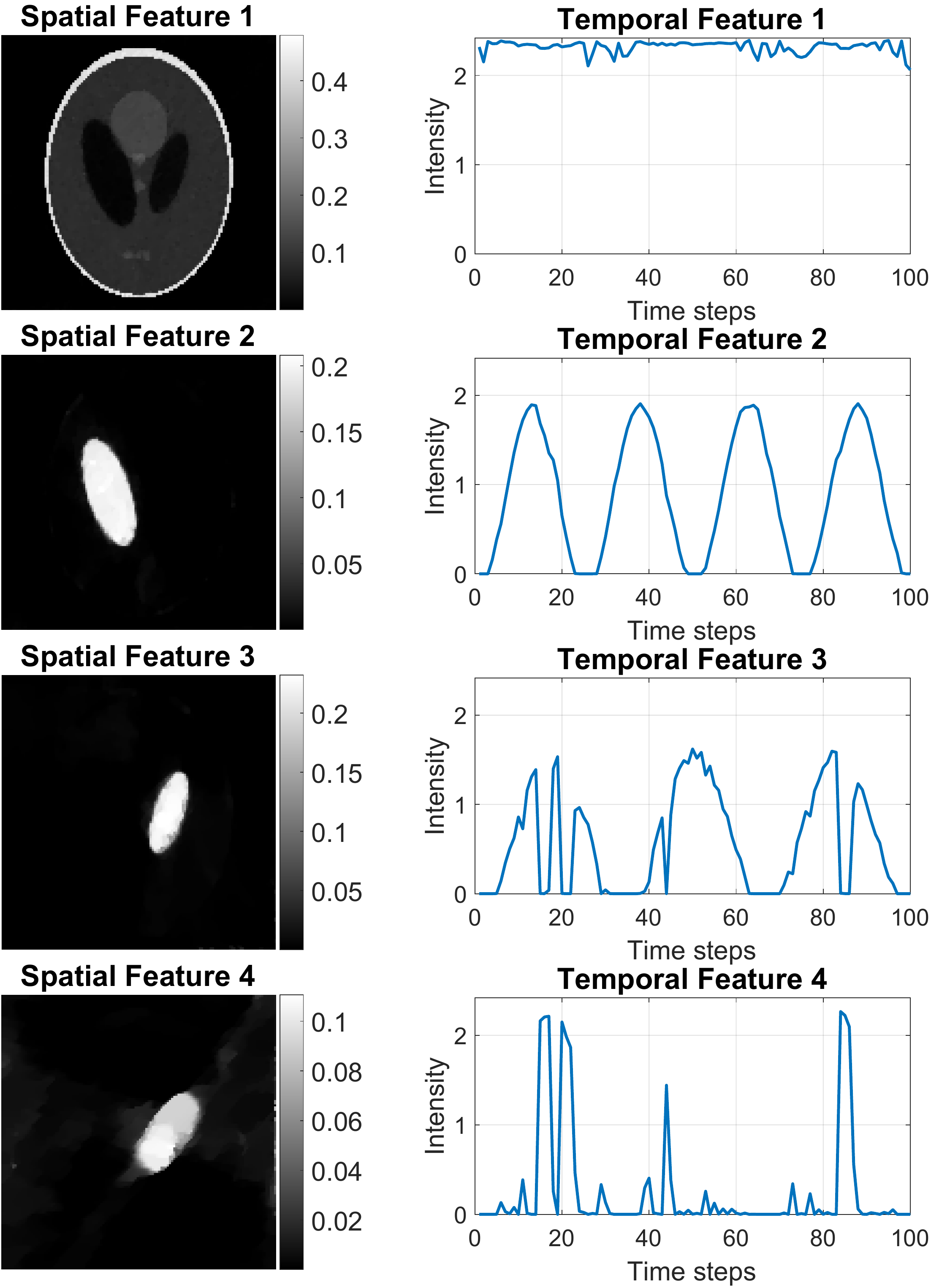}%
                \label{fig:dynShepp_BC_0.01_nTheta3}}
                \hfill
                \subfloat[\ref{eq:NMF model BC-X}]{\includegraphics[width=0.46\textwidth]{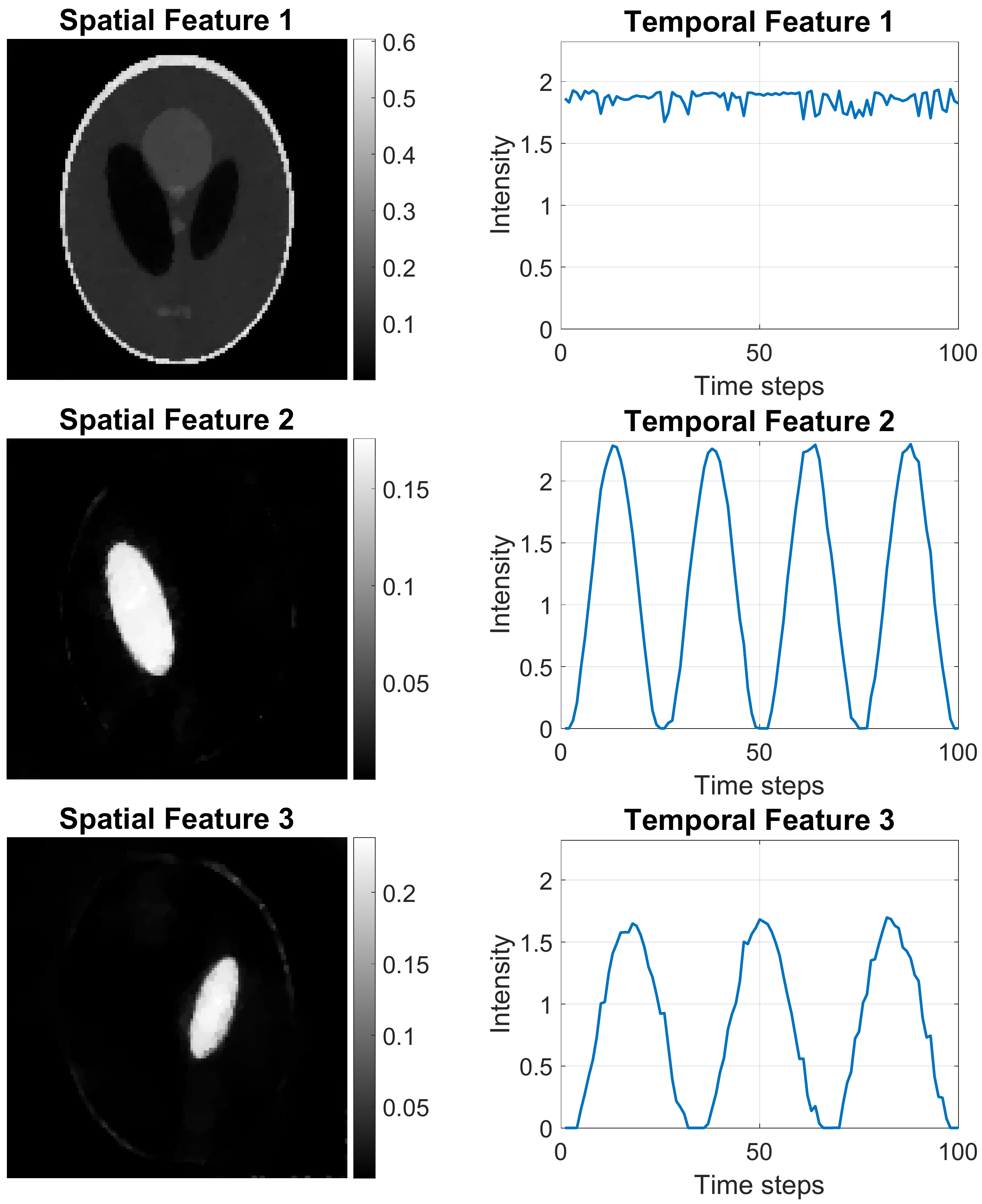}%
                \label{fig:dynShepp_BC-X_0.01_nTheta3}}
                \caption{Results for the dynamic Shepp-Logan phantom with $\vert \mathcal{I}_t \vert =3$ angles per time step and 1\% Gaussian noise. Shown are the leading extracted features for the \ref{eq:NMF model BC} model (left) and  for \ref{eq:NMF model BC-X} (right).}
                \label{fig:dynShepp:joint:noise1:nTheta3}
            \end{figure}
            
            \noindent Similar observations can be made for the case $\vert \mathcal{I}_t \vert =6$ and 3\% Gaussian noise. We present the reconstructed features in Figure \ref{fig:dynShepp:noise3:nTheta6} for \ref{eq:NMF model BC} and \texttt{gradTV\_PCA} only. The higher amount of noise can be observed especially in the spatial features of \texttt{gradTV\_PCA}, whereas it only has a slight effect in the \ref{eq:NMF model BC} model.
            
            Finally, we present the reconstructed features with \ref{eq:NMF model BC} and \ref{eq:NMF model BC-X} in Figure \ref{fig:dynShepp:joint:noise1:nTheta3} for $\vert \mathcal{I}_t \vert =3$, i.e.\ only three three angles per time step with a noise level of 1\%.
            The major difference to the previous cases can be seen in the results of the \ref{eq:NMF model BC} model. Here, the method splits up the dynamics of the right ellipse into two different temporal features, such that the true dynamics are not retained. However, the \ref{eq:NMF model BC-X} approach perform remarkably well with respect to the feature extraction despite the rather low number of projection angles. This might indicate, that enforcing the reconstruction $X$ to have small data error helps in the \ref{eq:NMF model BC-X} model to stabilise the reconstruction in highly sparse data settings.
            
            Let us shortly discuss other considered values of $\vert \mathcal{I}_t \vert,$ that are not shown here. First of all, the performance of \texttt{gradTV\_PCA} and \texttt{gradTV\_NMF} with respect to the feature extraction behaves very similar for both noise cases. Besides the above mentioned drawbacks, both approaches give remarkably consistent results especially for low number of angles and do not tend as much to split up features like in \ref{eq:NMF model BC} and \ref{eq:NMF model BC-X}. The latter occurs in different degrees for several numbers of angles. For 1\% noise, it occurs for $\vert \mathcal{I}_t \vert \in \{ 3, 7, 8, 10\} $ in \ref{eq:NMF model BC} and for $\vert \mathcal{I}_t \vert = 10 $ in \ref{eq:NMF model BC-X}. In the case of a noise level of 3\%, the split up effect only occurs for $\vert \mathcal{I}_t \vert = 10 $ in \ref{eq:NMF model BC}. However, for $\vert \mathcal{I}_t \vert = 10, $ it is possible to partially recover the correct temporal feature by simply adding up both features. Nevertheless, both approaches provide better reconstruction quality of $X$ than \texttt{gradTV} as we will discuss in the following.
            
\begin{figure}
                \centering
                \subfloat[Mean PSNR]{\includegraphics[width=0.48\textwidth]{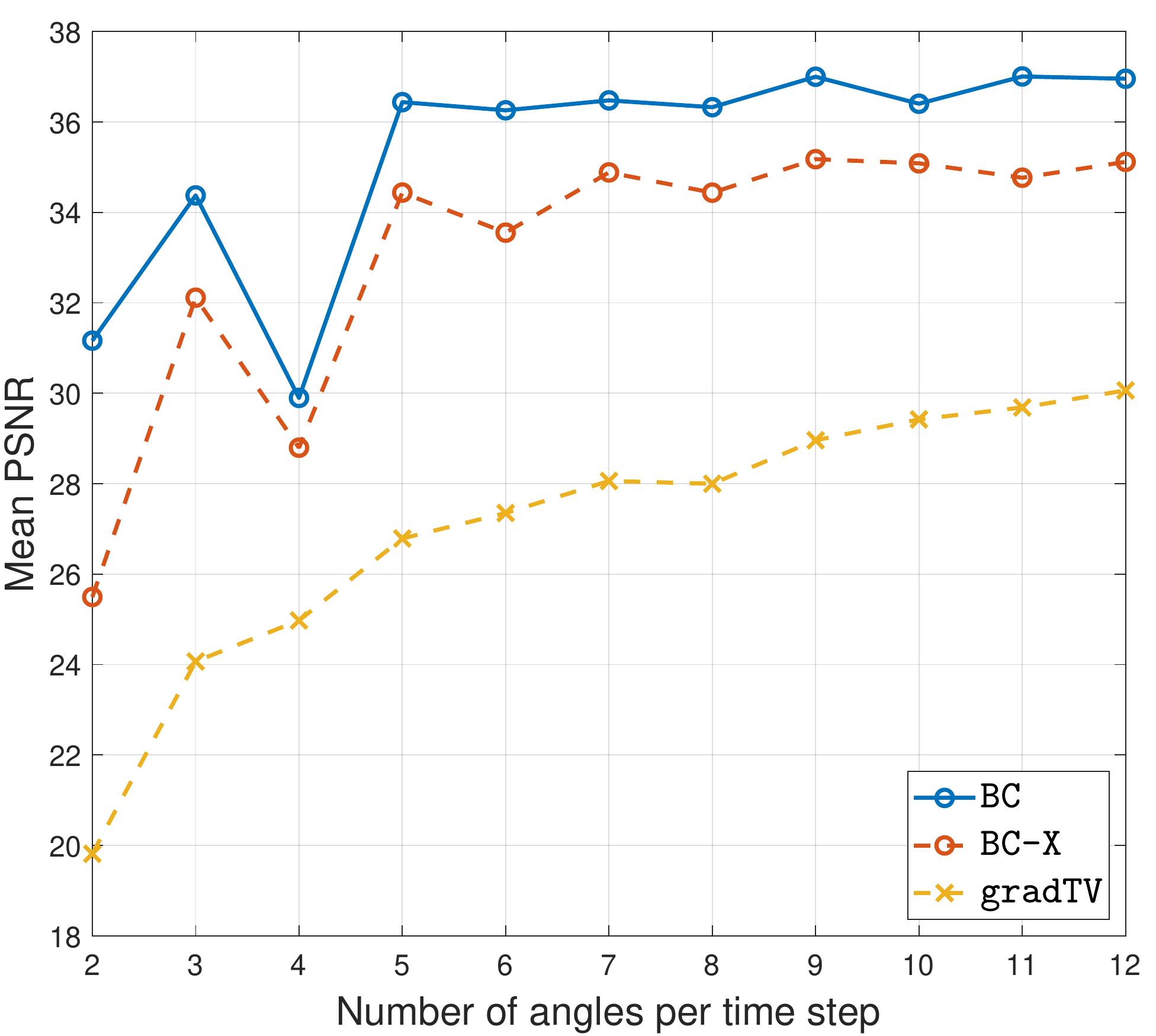}%
                \label{fig:dynShepp1PSNRPlot}}
                \hfill
                \subfloat[Mean SSIM]{\includegraphics[width=0.48\textwidth]{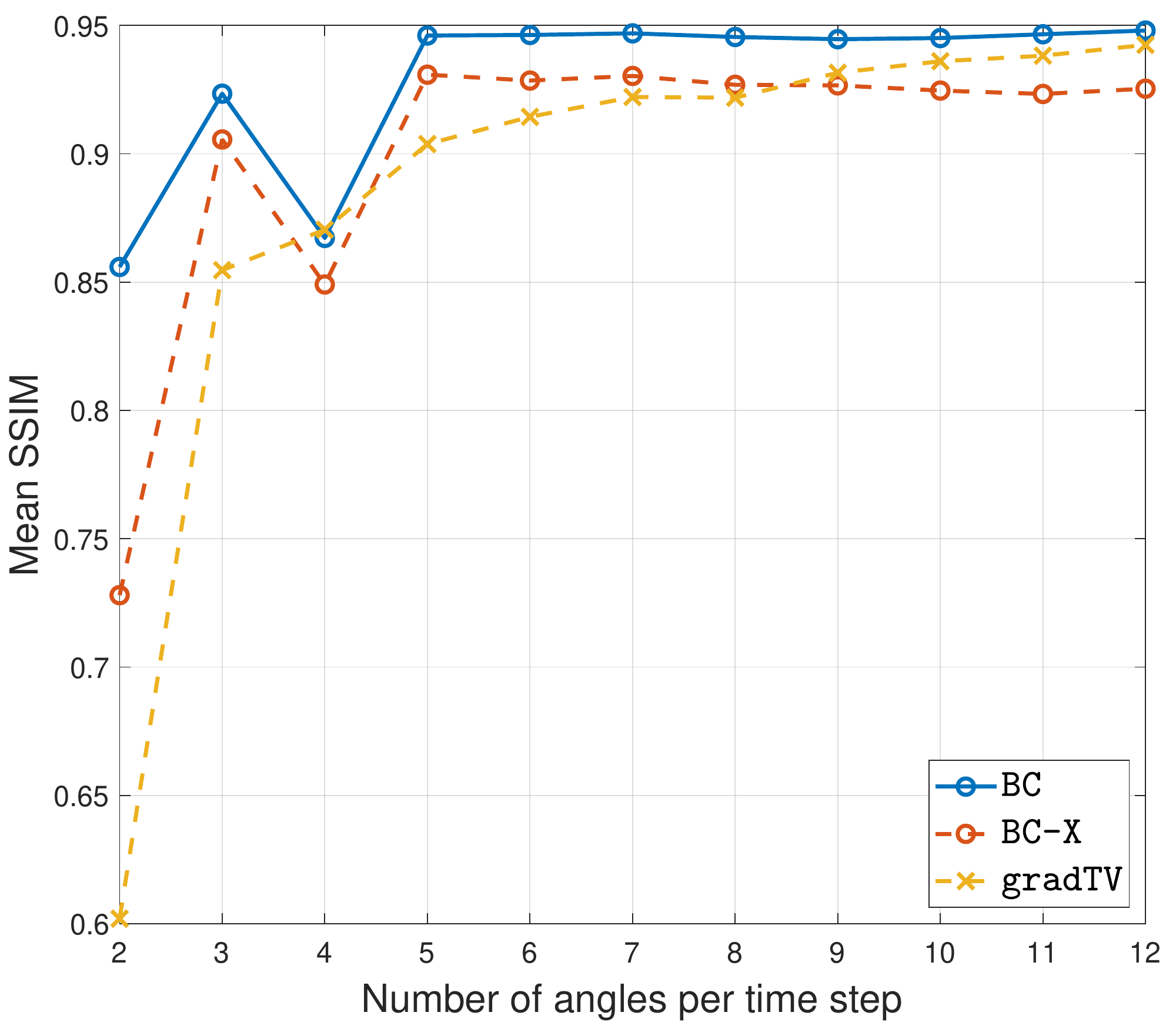}%
                \label{fig:dynShepp1SSIMPlot}}
                \caption{Mean PSNR and SSIM values of the reconstructions of the dynamic Shepp-Logan phantom with 1\% Gaussian noise for different numbers of projection angles.}
                \label{fig:dynShepp:noise1:qualityMeasures}
            \end{figure}

\begin{figure}
                \centering
                \subfloat[Mean PSNR]{\includegraphics[width=0.48\textwidth]{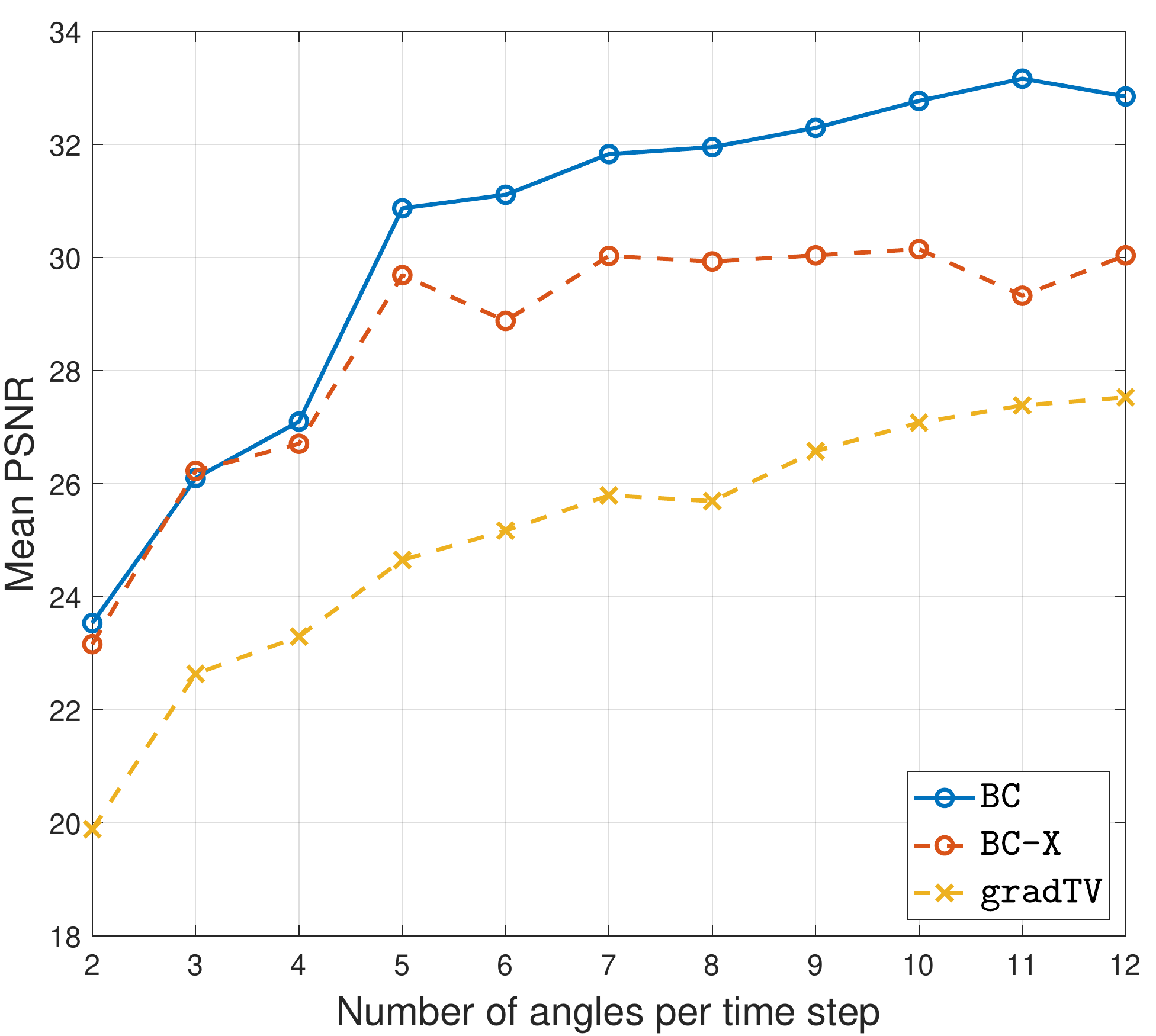}%
                \label{fig:dynShepp3PSNRPlot}}
                \hfill
                \subfloat[Mean SSIM]{\includegraphics[width=0.48\textwidth]{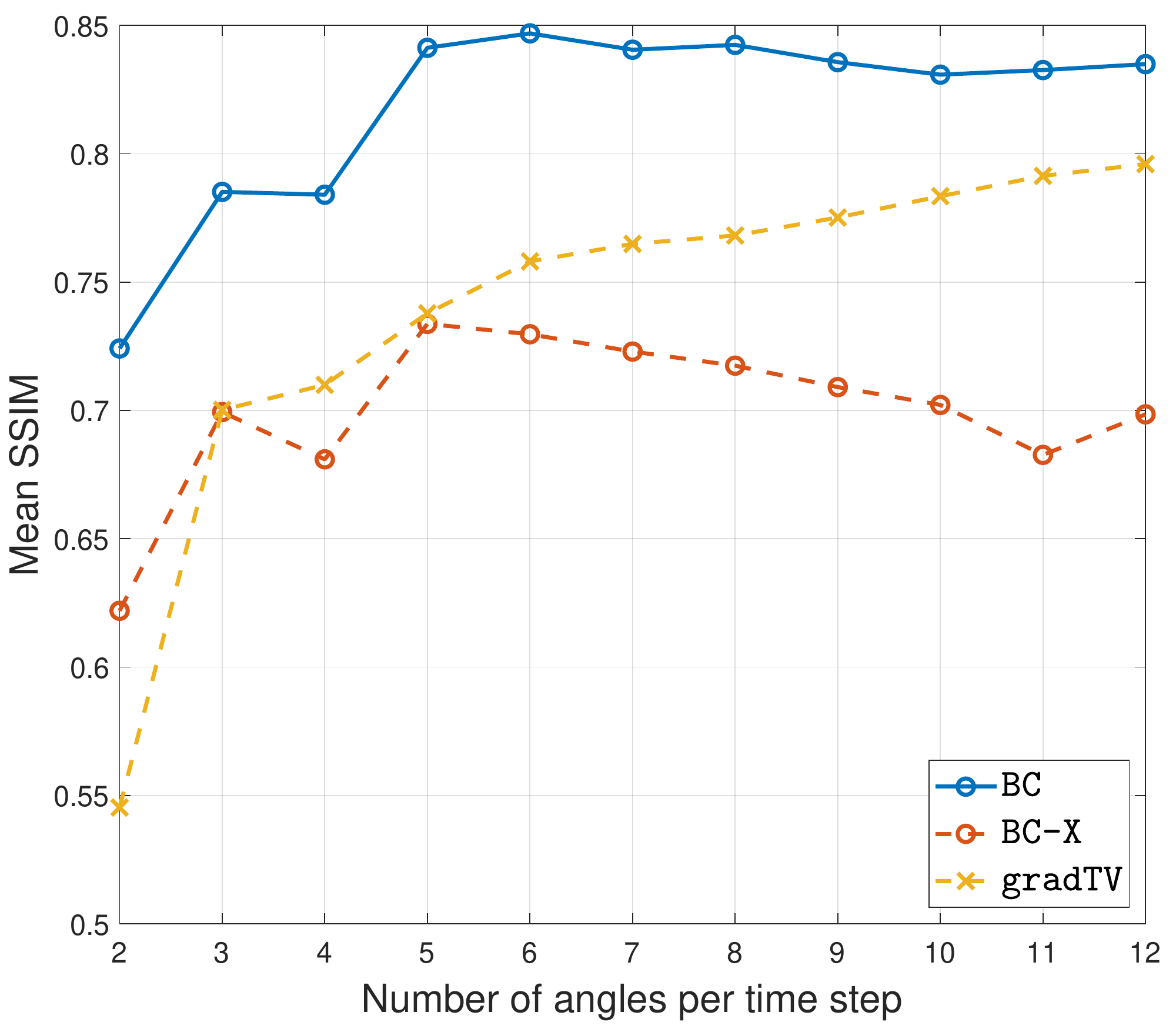}%
                \label{fig:dynShepp3SSIMPlot}}
                \caption{Mean PSNR and SSIM values of the reconstructions of the dynamic Shepp-Logan phantom with 3\% Gaussian noise for different numbers of projection angles.}
                \label{fig:dynShepp:noise3:qualityMeasures}
            \end{figure}

            \paragraph{Quantitative Evaluation}
            Let us now discuss the quantitative reconstruction quality for all methods. In Figure \ref{fig:dynShepp:noise1:qualityMeasures} and \ref{fig:dynShepp:noise3:qualityMeasures}, we show the mean \ac{PSNR} and \ac{SSIM} of the reconstructions for 1\% and 3\% noise over all time steps for all considered numbers of projection angles. Note that for the NMF model \ref{eq:NMF model BC-X}, we compute the quality measures for $X.$ The same goes for \texttt{gradTV}, where we only compute the quality measures of $X$ after the reconstruction procedure independently of the subsequent feature extraction method. In the case of \ref{eq:NMF model BC}, the reconstruction is computed as $X=BC$.
            
            As expected, the reconstruction quality tends to get better if more angles per time step are considered. More importantly, we see that it is possible to obtain reasonable reconstructions with just a few projections per time step especially in the case of the joint reconstruction and feature extraction method via the NMF approach. In particular, we reach a stable reconstruction quality already with 5 or more angles for both joint methods and 1\% noise.

            \begin{figure}
            	\centering
            	\includegraphics[width=0.5\textwidth]{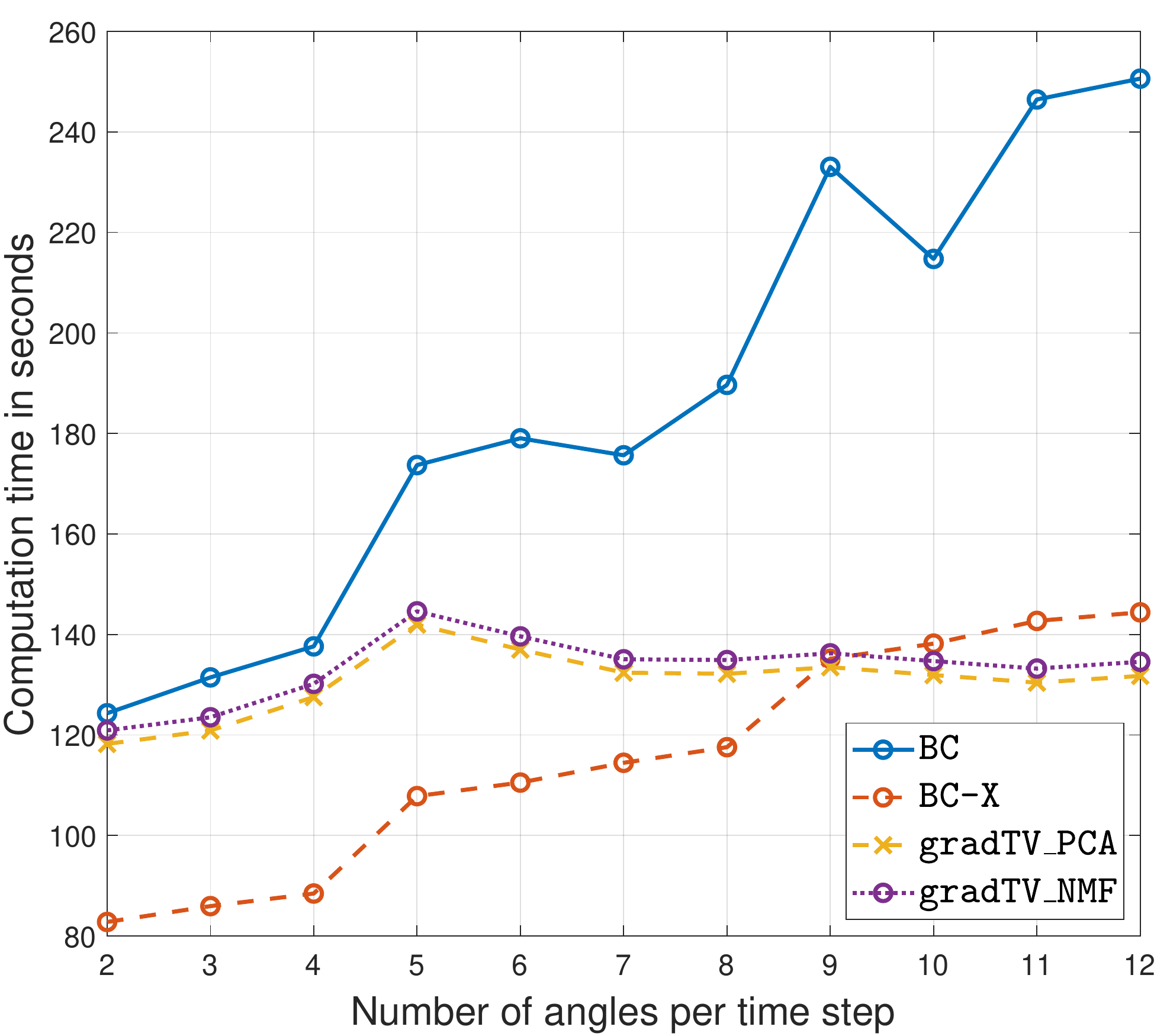}
            	\caption{Needed time in seconds for the reconstruction and feature extraction of the dynamic Shepp-Logan phantom with 1\% Gaussian noise.}
            	\label{fig:dynShepp:computationTime}
            \end{figure}
            The \ref{eq:NMF model BC} model 
            clearly performs best with respect to the reconstruction quality. For almost every number of angles, the mean \ac{PSNR} and \ac{SSIM} values outperform the ones of the \ref{eq:NMF model BC-X} and \texttt{gradTV} method for both noise levels. 
            In the case of 3\% noise (see Figure \ref{fig:dynShepp:noise3:qualityMeasures}) we can see that \texttt{gradTV} performs slightly better than \ref{eq:NMF model BC-X} in most of the cases in terms of their \ac{SSIM} values. Still, the mean \ac{PSNR} values of \texttt{gradTV} are significantly lower than the ones in \ref{eq:NMF model BC-X} for all numbers of angles. A selection of reconstructions for the experiments in Figure \ref{fig:dynShepp:noise1:qualityMeasures} and \ref{fig:dynShepp:noise3:qualityMeasures} are provided as videos in the Supplementary files.
            
            Note that for \ref{eq:NMF model BC-X}, it is also possible to compute the reconstruction based on the decomposition $B\cdot C$ instead of the joint reconstruction $X$ in the algorithm. 
            Interestingly, our experiments showed that the reconstruction quality of $B\cdot C$ is in almost all cases better than the one of the matrix $X$ itself and also mostly outperforms the \texttt{gradTV} approach. We believe, that this is due to the stronger regularising effect on the components $B$ and $C$, which especially influences SSIM.
            
            The computation times for the reconstruction and feature extraction with 1\% noise for all algorithms until the stopping criterion is fulfilled are shown in Figure \ref{fig:dynShepp:computationTime}.
            As expected, the computation time tends to increase with the number of projection angles and, considering all methods, ranges approximately from 1 to 5 minutes. For $ \vert \mathcal{I}_t \vert\leq 8, $ the \ref{eq:NMF model BC-X} method is the fastest while it is outperformed by \texttt{gradTV\_PCA} for $ \vert \mathcal{I}_t \vert \geq 9.$ \texttt{gradTV\_NMF} and \ref{eq:NMF model BC} needs more time in all experiments compared to \texttt{gradTV\_PCA}. The significant temporal difference between \ref{eq:NMF model BC-X} and \ref{eq:NMF model BC} is due to its higher computational complexity: Owing to the model formulation of \ref{eq:NMF model BC} with the discrepancy term $\Vert R_t (BC)_{\bullet, t} - Y_{\bullet, t} \Vert_2^2,$ the update rules in Theorem \ref{theorem:Algorithm BC} for both matrices $B$ and $C$ contain the discretised Radon transform $R_t.$ This is in contrast to the \ref{eq:NMF model BC-X} algorithm, where $R_t$ only appears in the update rule of $X.$
            
            Based on the presented results for the dynamic Shepp-Logan phantom, we can conclude 
            that the joint approaches \ref{eq:NMF model BC} and \ref{eq:NMF model BC-X} outperform both other methods with respect to the reconstruction quality and for most cases of the extracted features. Nevertheless, the models \texttt{gradTV\_PCA} and \texttt{gradTV\_NMF} give remarkably consistent and stable results of the extracted components throughout all numbers of angles. Furthermore, the nonnegativity constraint of the NMF improves significantly the interpretability and quality of the extracted spatial features.
            \paragraph{Stationary Operator}
            As we have seen, the computational complexity of the \ref{eq:NMF model BC} model with the non-stationary operator is clearly higher than for all other cases. Thus, let us now consider the possibility to speed up the reconstructions with a stationary operator, which leads us to the complexity reduced formulation presented in Corollary \ref{cor:Algorithm sBC} as the \ref{eq:NMF model sBC} model.
            Here, we present similarly to the case above experiments with the dynamic Shepp-Logan phantom for $\vert \mathcal{I}_t \vert \in \{ 2,\dots,30 \} $ and 1\% Gaussian noise, as we primarily aim to illustrate the reduction of the computational cost. Furthermore, the same hyperparameters and stopping criteria are used as before.
            
\begin{figure}
                    \centering
                    \subfloat[$\vert \mathcal{I}_t \vert = 6 $]{\includegraphics[width=0.45\textwidth]{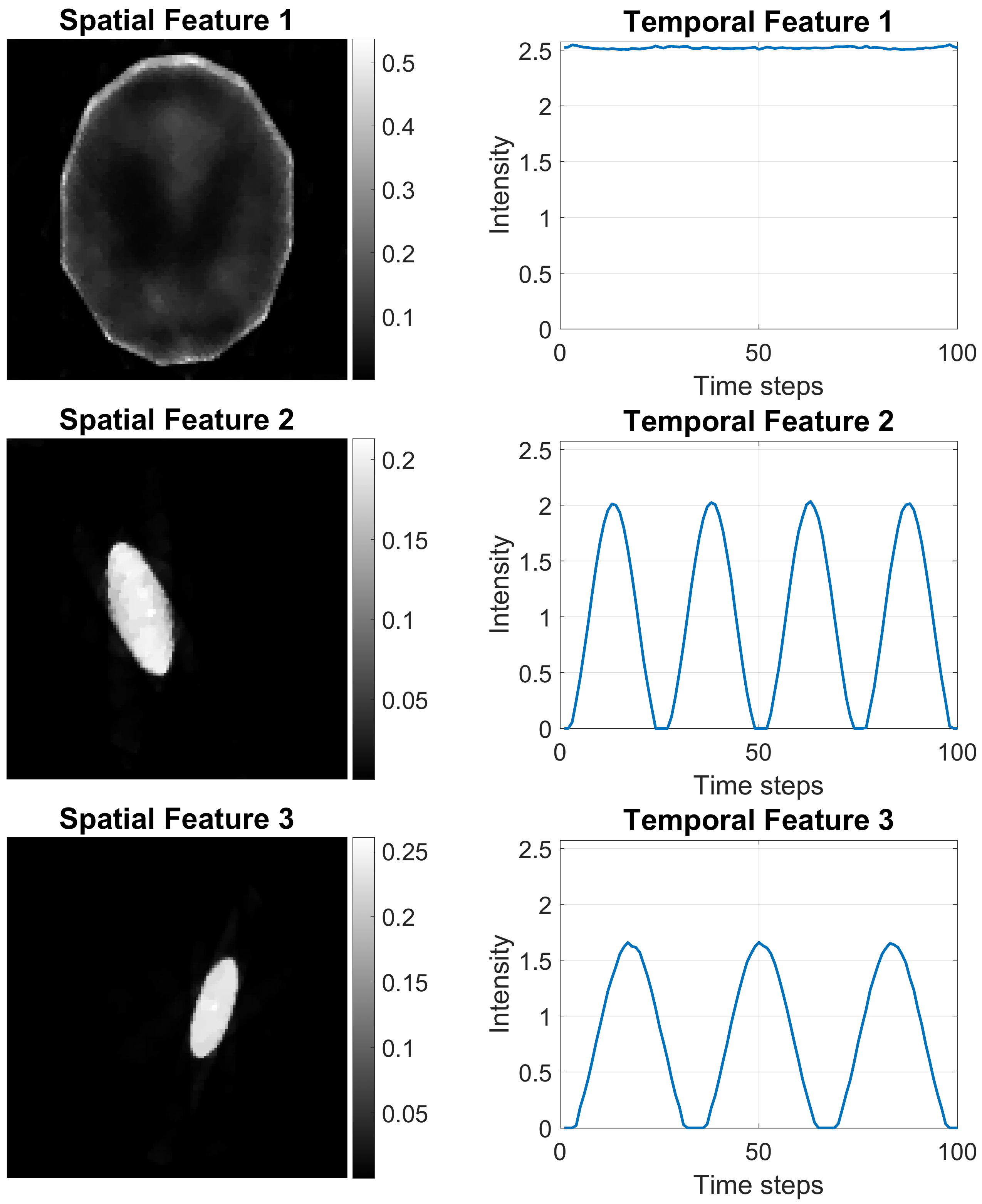}%
                    \label{fig:dynShepp_sBC_0.01_nTheta6}}
                    \hfill
                    \subfloat[$\vert \mathcal{I}_t \vert = 30 $]{\includegraphics[width=0.45\textwidth]{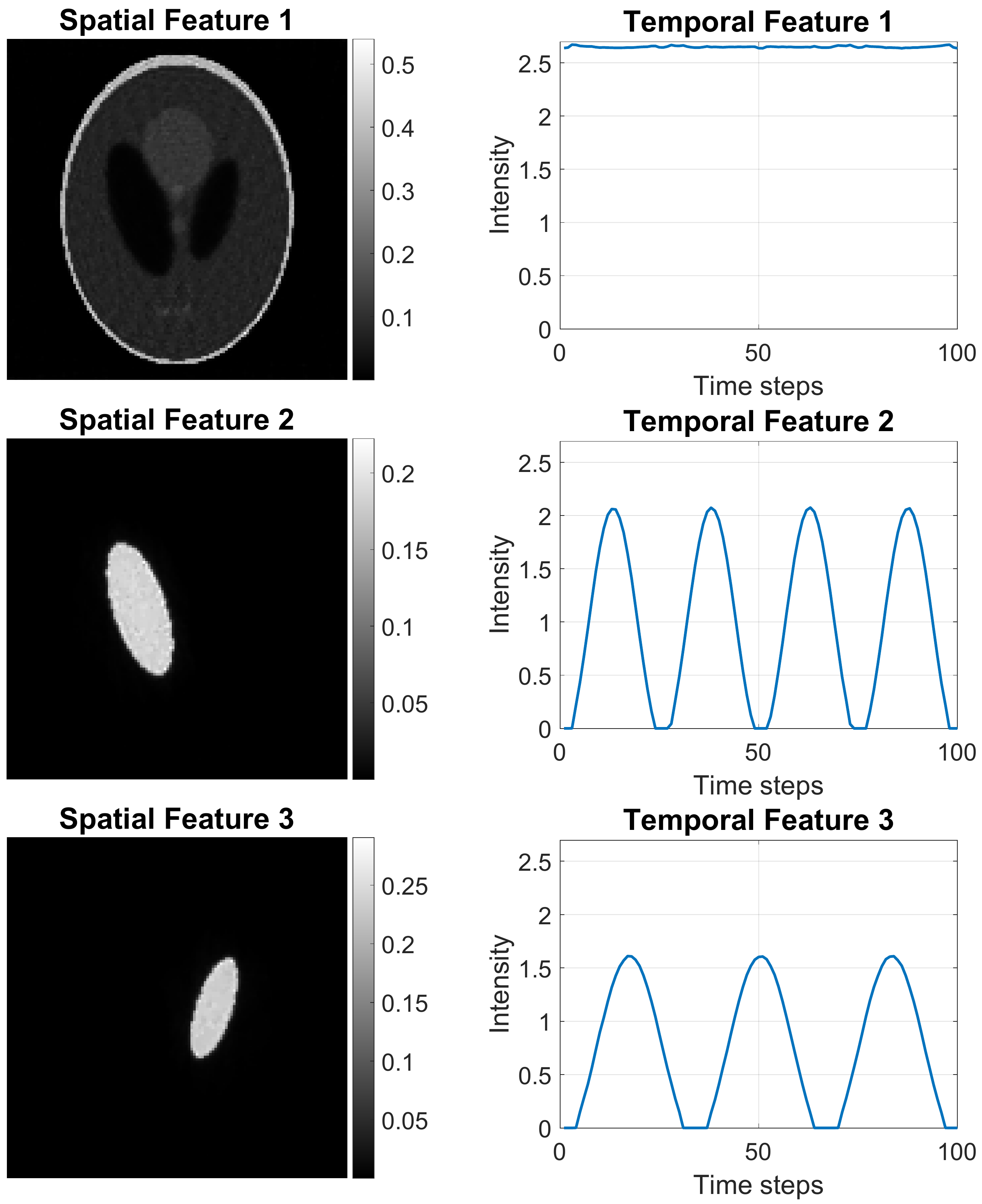}%
                    \label{fig:dynShepp_sBC_0.01_nTheta30}}
                    \caption{Results for the dynamic Shepp-Logan phantom with a stationary operator and 1\% Gaussian noise. Shown are the leading extracted features with the \ref{eq:NMF model sBC} model for $\vert \mathcal{I}_t \vert = 6$ angles per time step (left) and $\vert \mathcal{I}_t \vert = 30$ (right).}
                    \label{fig:dynShepp:sBC:noise1}
                \end{figure}

            The reconstructed features for the cases $\vert \mathcal{I}_t \vert = 6 $ and $\vert \mathcal{I}_t \vert = 30$ are shown in 
            Figure \ref{fig:dynShepp:sBC:noise1}. In particular, comparing the results in Figure \ref{fig:dynShepp_sBC_0.01_nTheta6} to the corresponding results of \ref{eq:NMF model BC} in Figure \ref{fig:dynShepp_BC_0.01_nTheta6}, one can immediately see a significant difference between the extracted spatial features. This is clearly due to the fact that the same projection angles are used at every time step and the individual projection directions are clearly visible for the stationary model \ref{eq:NMF model sBC}.
            Consequently, the details in the Shepp-Logan phantom are not well recovered, such that the extracted constant feature is significantly inferior to the one of \ref{eq:NMF model BC}. As one would expect, more projection angles per time step are needed to reconstruct finer details. This effect can be clearly seen for 30 angles in Figure \ref{fig:dynShepp_sBC_0.01_nTheta30}.
            
            However, all temporal basis functions with \ref{eq:NMF model sBC} for $\vert \mathcal{I}_t \vert = 6 $ are remarkably well reconstructed despite the low number of projection angles. This is also true for the other considered values of $\vert \mathcal{I}_t \vert.$ Moreover, we observe that \ref{eq:NMF model sBC} is able to extract the correct three main features for every $\vert \mathcal{I}_t \vert \in \{ 2,\dots,30 \}.$ Even for $\vert \mathcal{I}_t \vert = 2,$ the quality of the dynamic temporal features are similar to the ones in Figure \ref{fig:dynShepp:sBC:noise1}.
            
            This behaviour is different from the dynamic case discussed above. The reason for this is probably based on the different projection directions at every time step in the dynamic case, which results in directional dependencies of the occurring reconstruction artefacts in contrast to the stationary case. This can make it difficult for the NMF to distinguish the main features in the non-stationary case and thus leads to a more stable feature extraction in the here presented stationary case. 
                
\begin{figure}
                    \centering
                    \subfloat[Computation time]{\includegraphics[width=0.32\textwidth]{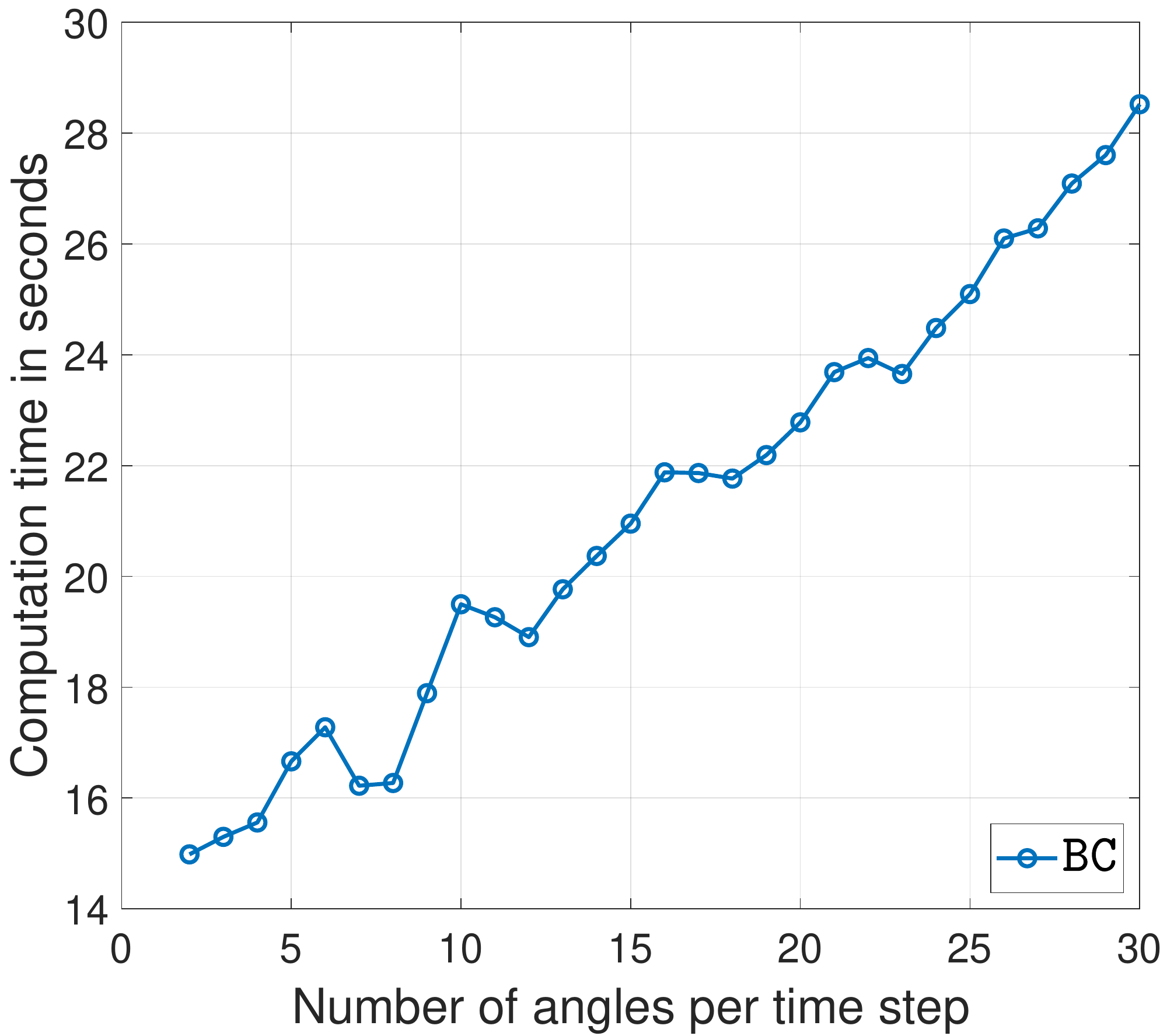}%
                    \label{fig:dynShepp:Stationary:computationTime}}
                    \hfill
                    \subfloat[Mean PSNR]{\includegraphics[width=0.32\textwidth]{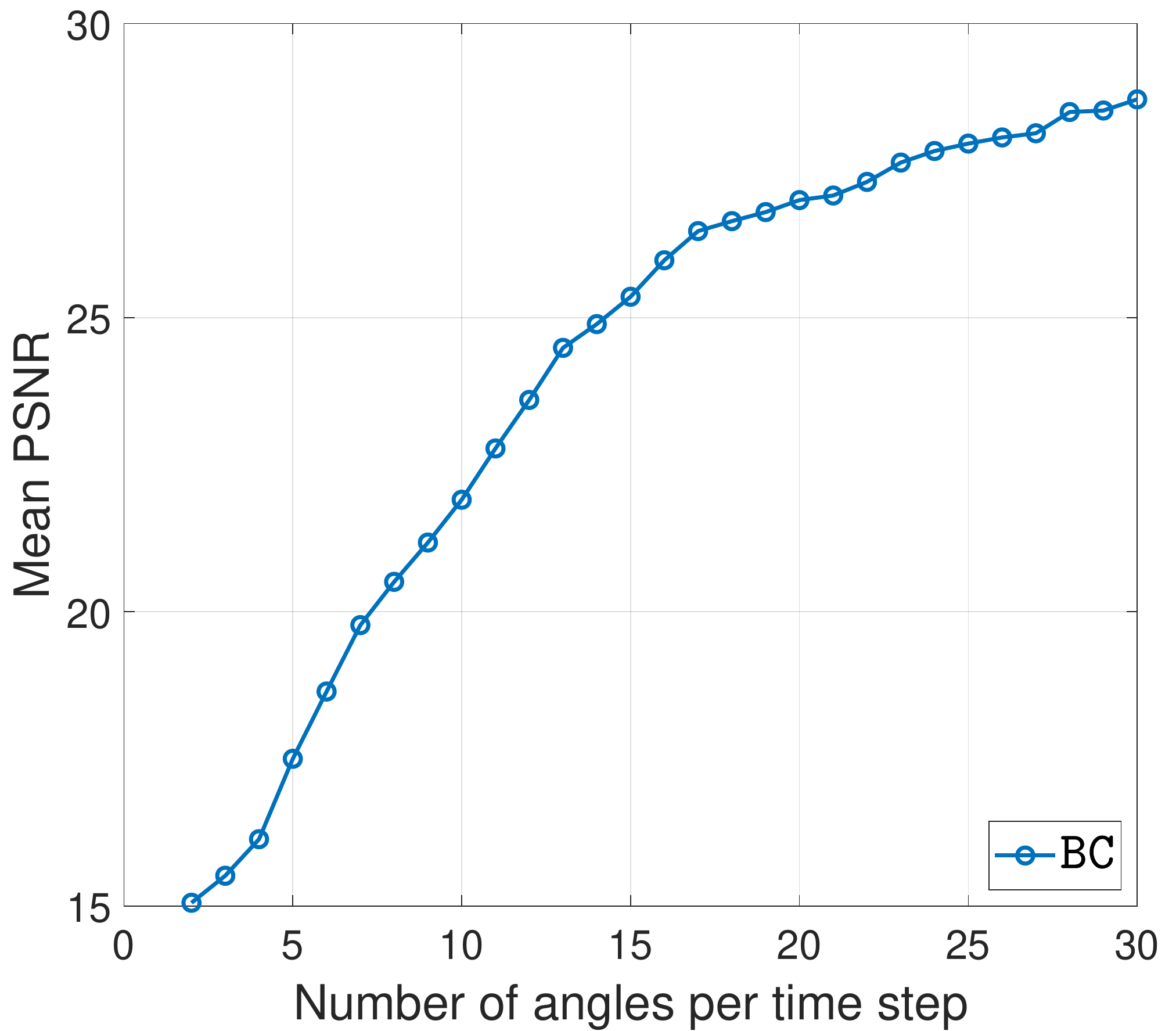}%
                    \label{fig:dynShepp1StationaryPSNRPlot}}
                    \hfill
                    \subfloat[Mean SSIM]{\includegraphics[width=0.32\textwidth]{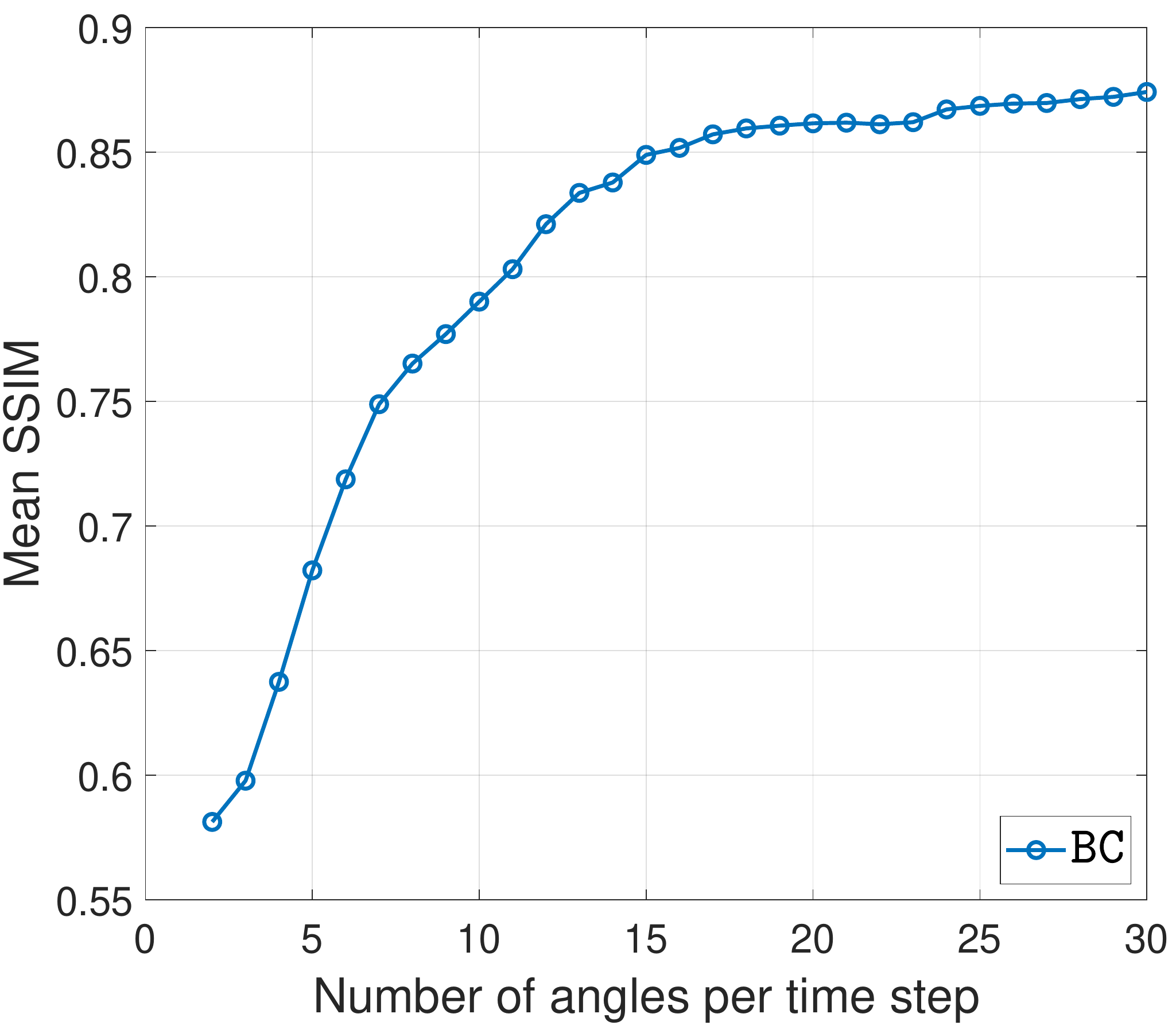}%
                    \label{fig:dynShepp1StationarySSIMPlot}}
                    \caption{Needed time in seconds, mean \ac{PSNR} and mean \ac{SSIM} values of the reconstructions of the dynamic Shepp-Logan phantom with 1\% Gaussian noise for the stationary case \ref{eq:NMF model sBC} and different numbers of projection angles.}
                    \label{fig:dynShepp:noise1:Stationary}
                \end{figure}

                The quantitative measures are shown in
                Figure \ref{fig:dynShepp:noise1:Stationary} for all experiments. Comparing the computation time of \ref{eq:NMF model BC} with the one of \ref{eq:NMF model sBC}, we obtain a clear speed-up by a factor of 10--20 with the stationary model. However, as expected, comparing Figure \ref{fig:dynShepp1StationaryPSNRPlot} and \ref{fig:dynShepp1StationarySSIMPlot} with the quality measures of \ref{eq:NMF model BC} in Figure \ref{fig:dynShepp:noise1:qualityMeasures}, one can observe that significantly more projection angles per time step are needed in the stationary case to provide a sufficient reconstruction quality. In conclusion, we can say that the \ref{eq:NMF model sBC} model is especially recommended if one is primarily interested in the dynamics of the system under consideration, as we could extract the temporal basis functions stably for all considered angles with $\vert \mathcal{I}_t| \geq2$.
            

        \subsubsection{Vessel Phantom} \label{subsubsec:Vessel Phantom}
\begin{figure}
            	\centering
            	\includegraphics[width=1.0\textwidth]{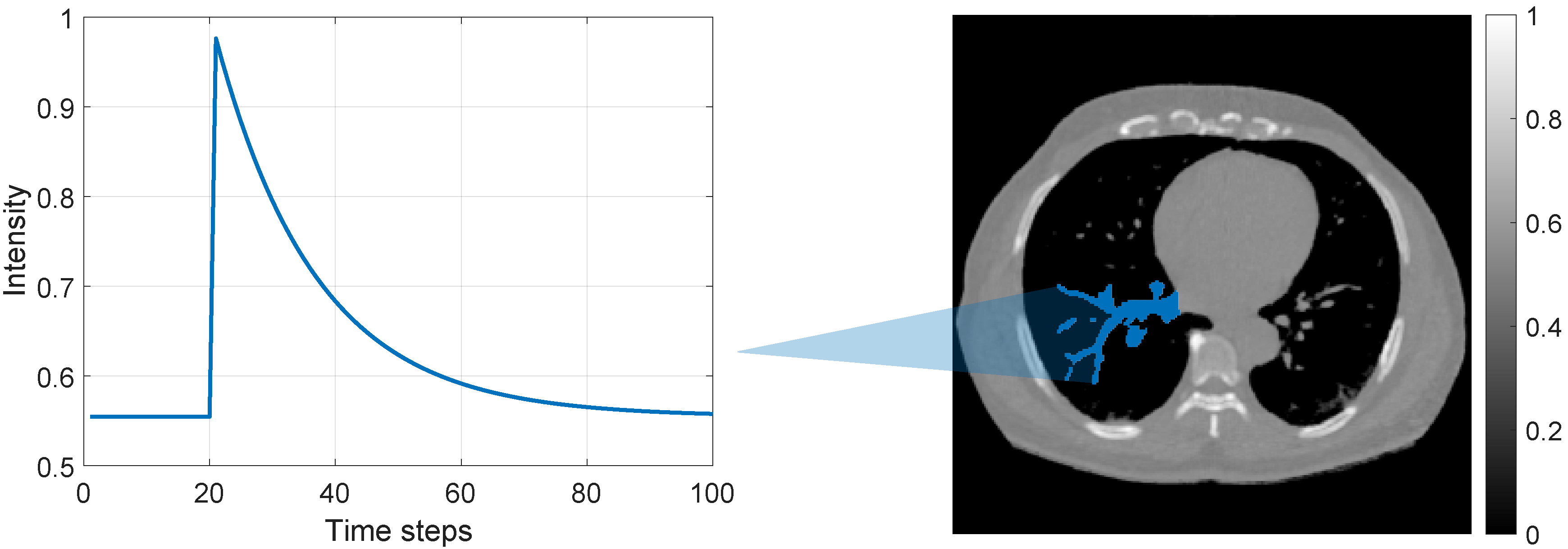}
            	\caption{Illustration of the vessel phantom dataset consisting of $ T=100 $ phantoms of dimension $264\times264,$ where the intensity of the blue highlighted area changes according to blue curve on the left.}
            	\label{fig:phantomVesselIllustration}
            \end{figure}

            The second test case is based on a CT scan of a human lung\footnote{The phantom is based on the CT scans in the \emph{ELCAP Public Lung Image database}: \url{http://www.via.cornell.edu/lungdb.html}}, see Figure  \ref{fig:phantomVesselIllustration}. Here, the decomposition is given by the constant background and a segmented vessel that exhibits a sudden increase in attenuation followed by an exponential decay. This could for instance represent the injection of a tracer to the blood stream.
            
            In contrast to the previous dataset, we perform only selected experiments for specific choices of noise levels and numbers of projection angles. More precisely, we present results for 1\% Gaussian noise together with $\vert \mathcal{I}_t \vert \in \{7, 12\}$ and 3\% Gaussian noise with $\vert \mathcal{I}_t \vert =12.$ In all cases, we choose $K=4$ NMF features. Furthermore, the stopping criterion from the experiments with the dynamic Shepp-Logan phantom is changed for this dataset in such a way, that the maximum number of iterations is raised to 1400 to ensure sufficient convergence. The regularisation parameters of all methods are chosen empirically and are displayed in Table \ref{tab:phantomVessel:parameterChoice} in Appendix \ref{app:sec:parameter Choice}.
            \begin{figure}
                \centering
                \subfloat[\ref{eq:NMF model BC}]{\includegraphics[width=0.45\textwidth]{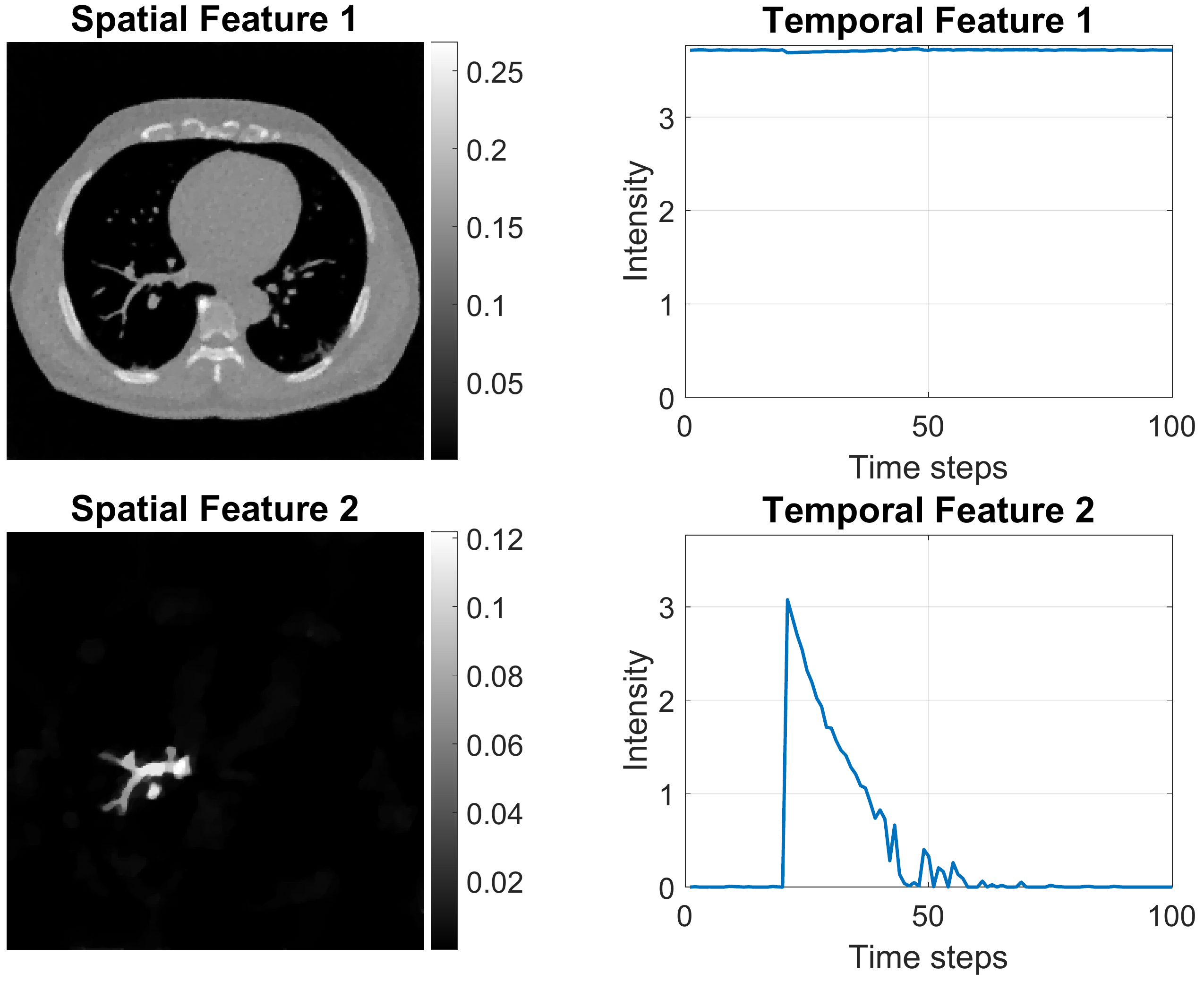}%
                \label{fig:phantomVessel_BC_0.01_nTheta12}}
                \hfill
                \subfloat[\ref{eq:NMF model BC-X}]{\includegraphics[width=0.45\textwidth]{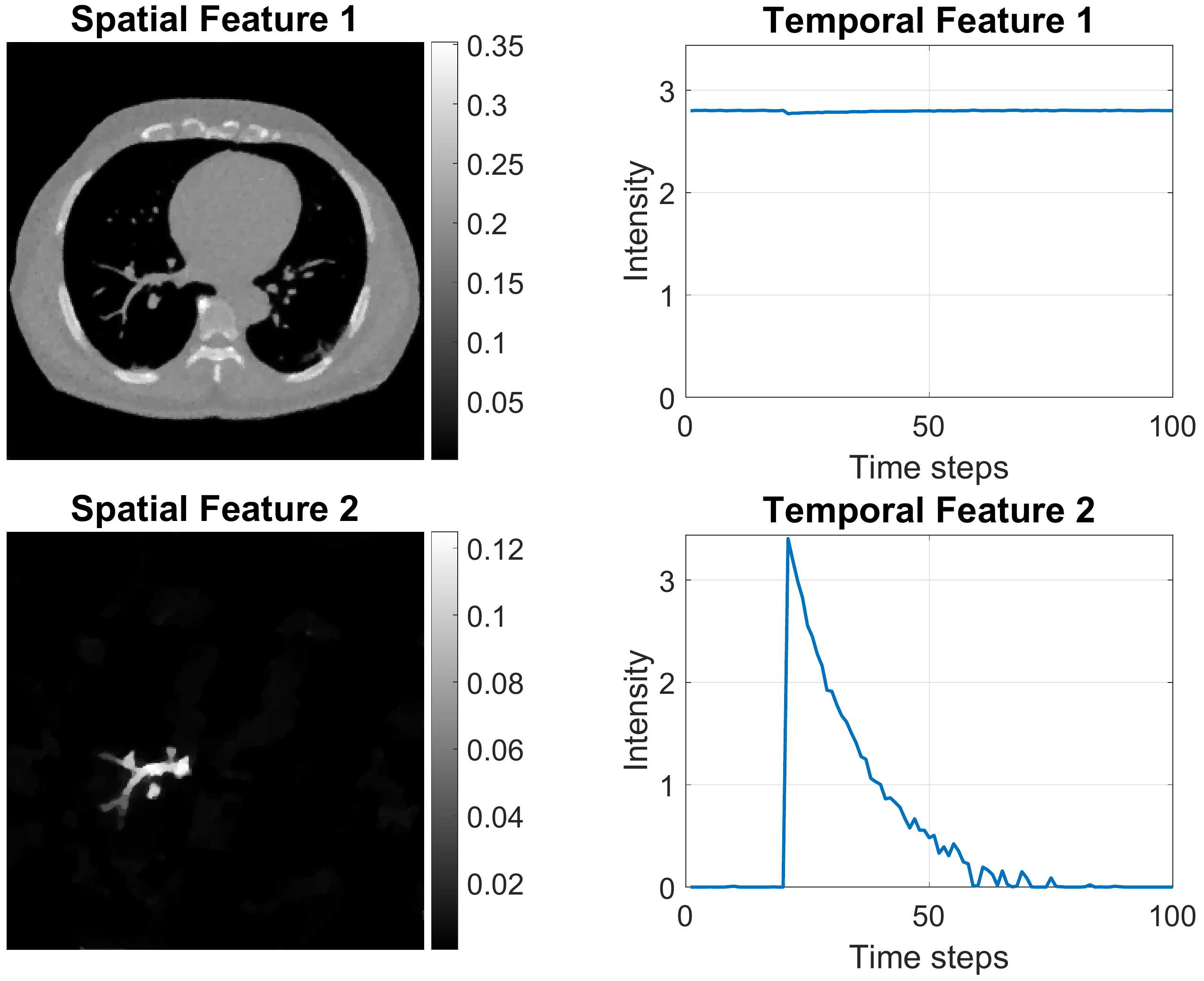}%
                \label{fig:phantomVessel_BC-X_0.01_nTheta12}}
                \caption{Results for the vessel phantom with $\vert \mathcal{I}_t \vert =12$ angles per time step and 1\% Gaussian noise. Shown are the leading extracted features for the \ref{eq:NMF model BC} model (left) and  for \ref{eq:NMF model BC-X} (right).}
                \label{fig:phantomVessel:joint:noise1:nTheta12}
            \end{figure}
            \begin{figure}
                \centering
                \subfloat[\texttt{gradTV\_PCA}]{\includegraphics[width=0.45\textwidth]{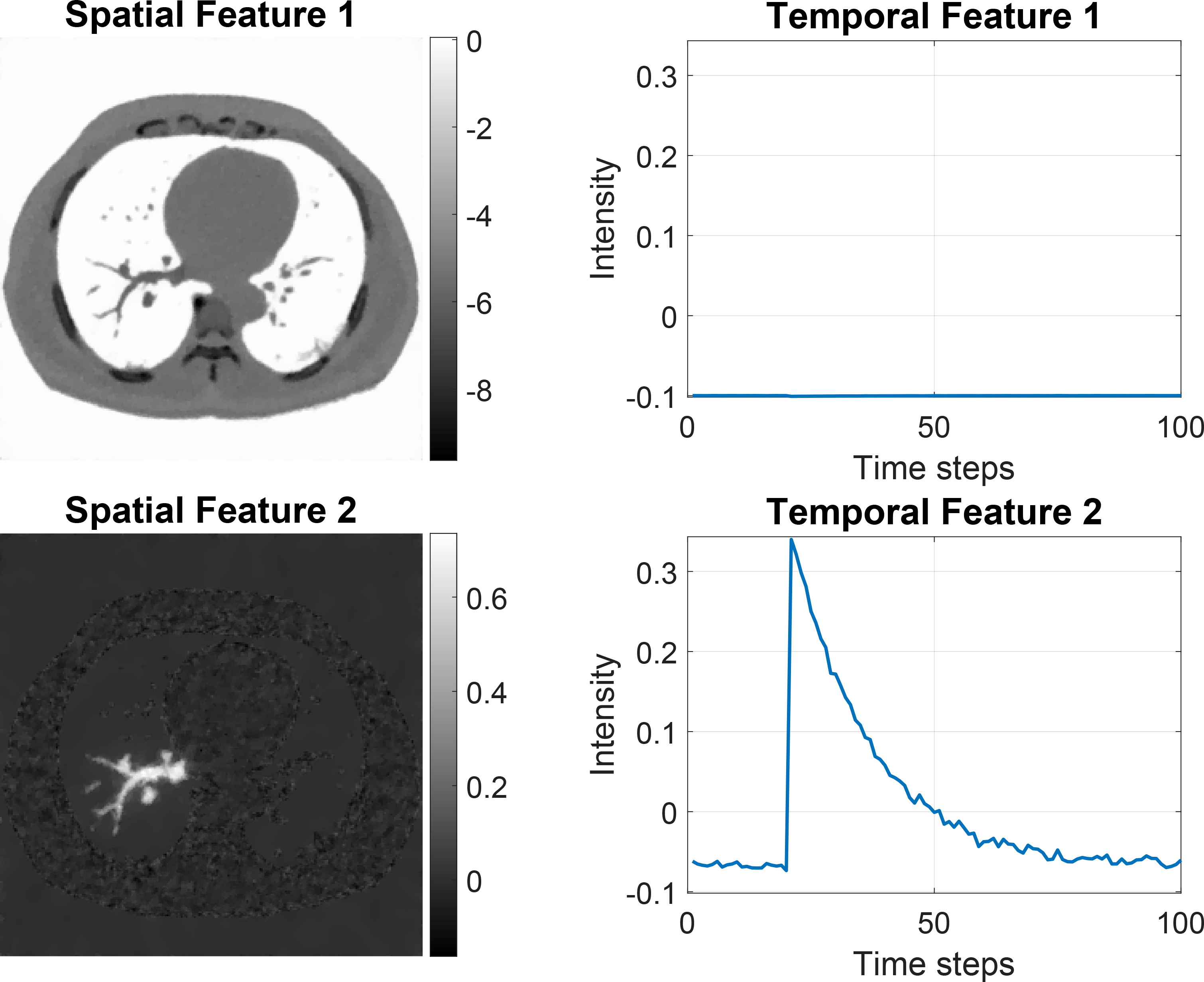}%
                \label{fig:phantomVessel_gradTV_PCA_0.01_nTheta12}}
                \hfill
                \subfloat[\texttt{gradTV\_NMF}]{\includegraphics[width=0.45\textwidth]{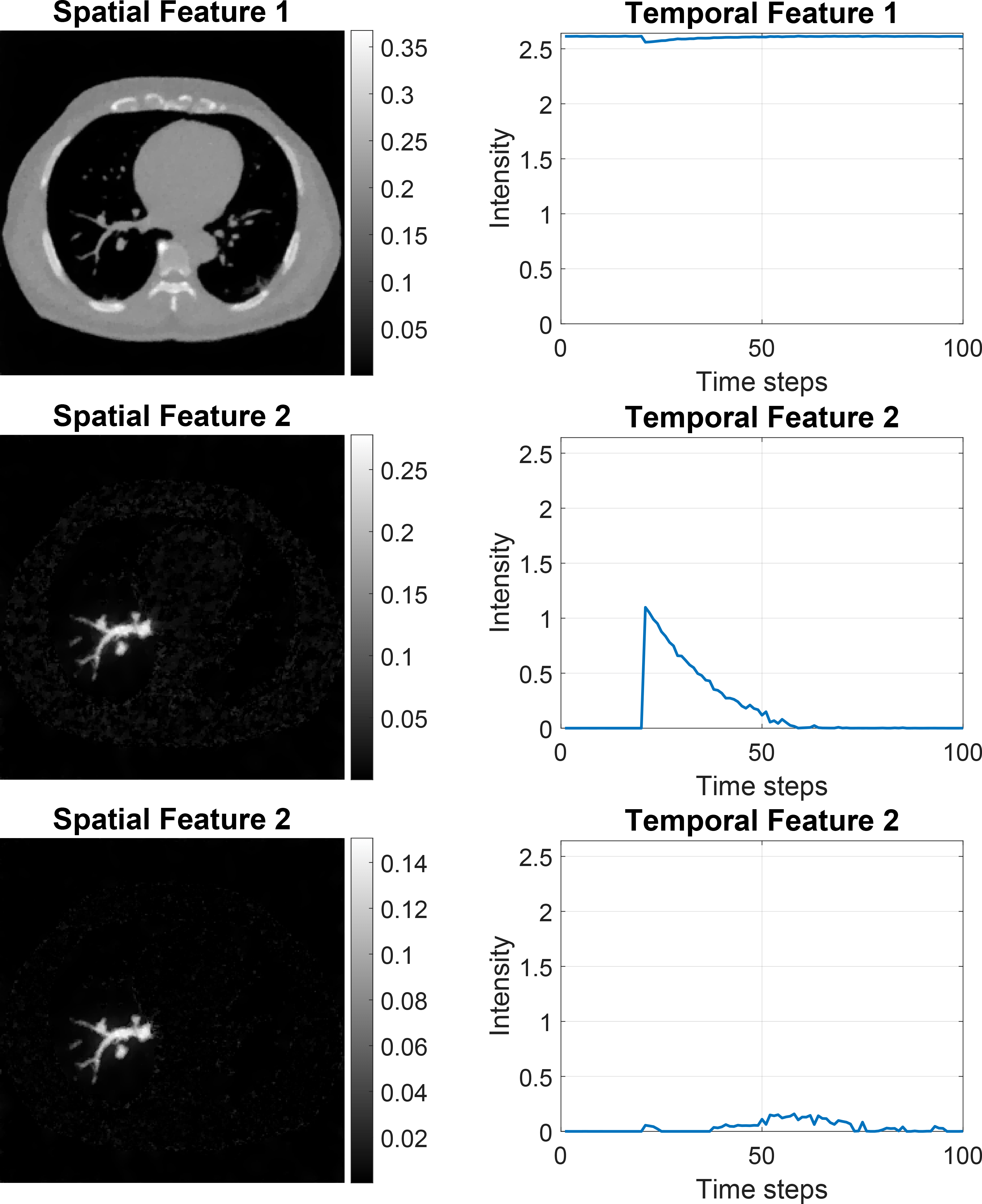}%
                \label{fig:phantomVessel_gradTV_NMF_0.01_nTheta12}}
                \caption{Results for the vessel phantom with $\vert \mathcal{I}_t \vert =12$ angles per time step and 1\% Gaussian noise. Shown are the leading extracted features for the \texttt{gradTV\_PCA} model (left) and  for \texttt{gradTV\_NMF} (right).}
                \label{fig:phantomVessel:sep:noise1:nTheta12}
            \end{figure}
            Figure \ref{fig:phantomVessel:joint:noise1:nTheta12} and \ref{fig:phantomVessel:sep:noise1:nTheta12} show the feature extraction results for the noise level of 1\% and $\vert \mathcal{I}_t \vert =12,$ where all approaches are able to extract both the main constant and dynamic component of the underlying ground truth. The order of the features here is based on a manual sorting.
            
            Similar to the results of the Shepp-Logan phantom in Section \ref{subsubsec:Shepp--Logan Phantom}, the joint methods \ref{eq:NMF model BC} and \ref{eq:NMF model BC-X} have difficulties to recover the lower intensities in the temporal features, whereas \texttt{gradTV\_PCA} produce slight artefacts in the dynamic spatial feature due to the missing nonnegativity constraint. In addition, \texttt{gradTV\_NMF} is able to recover more details in the vessel compared to the joint approaches. This is due to the relatively high choice of the total variation regularisation parameter $\tau$ in \ref{eq:NMF model BC} and \ref{eq:NMF model BC-X} to ensure a sufficient denoising effect on the matrix $B.$ The low peak in the second temporal feature of \texttt{gradTV\_NMF} is likely caused by the choice of the $\ell_2$ regularisation parameter $\tilde{\mu}_C.$
            \begin{figure}
                \centering
                \subfloat[\ref{eq:NMF model BC}]{\includegraphics[width=0.45\textwidth]{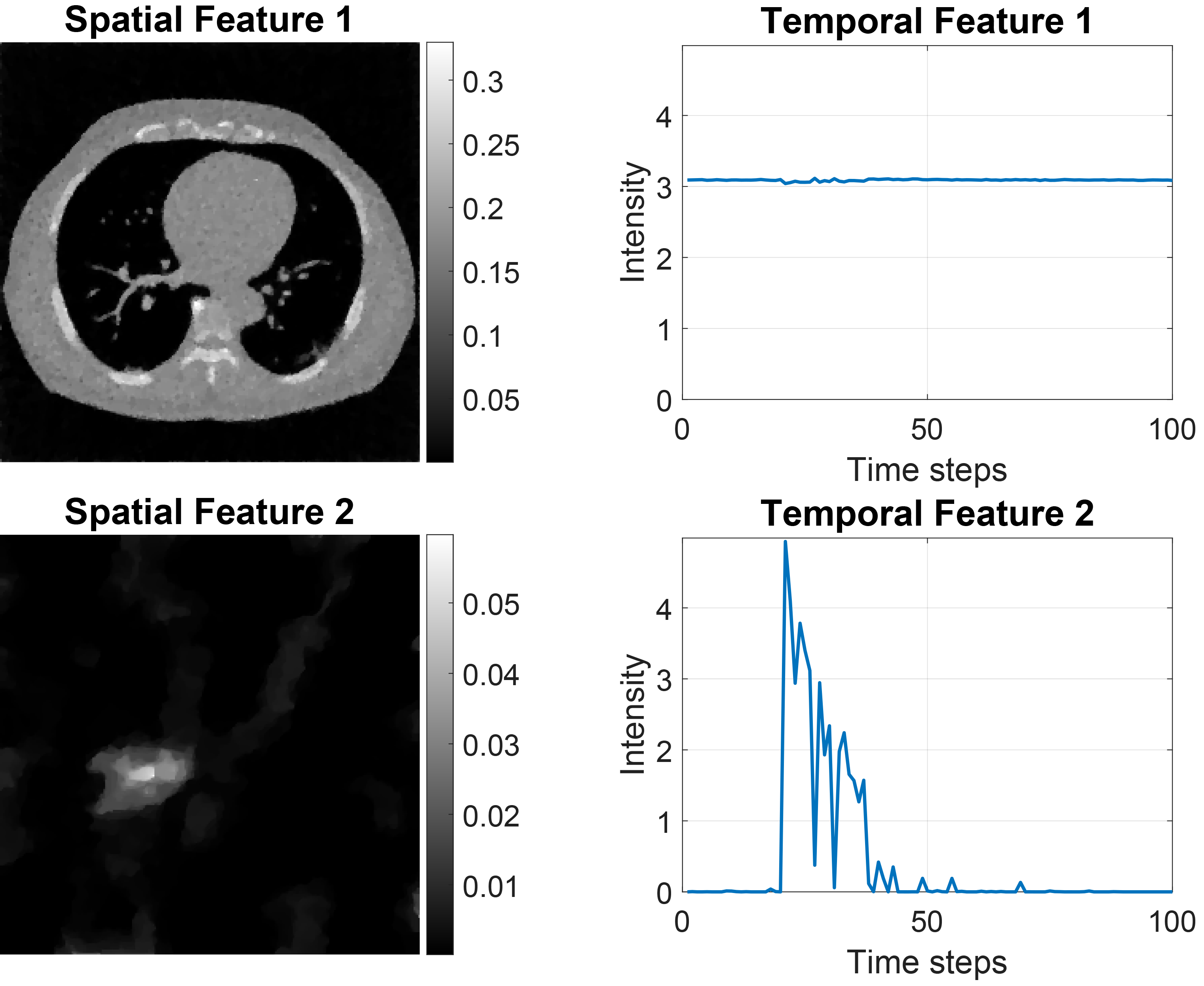}%
                \label{fig:phantomVessel_BC_0.03_nTheta12}}
                \hfill
                \subfloat[\texttt{gradTV\_PCA}]{\includegraphics[width=0.45\textwidth]{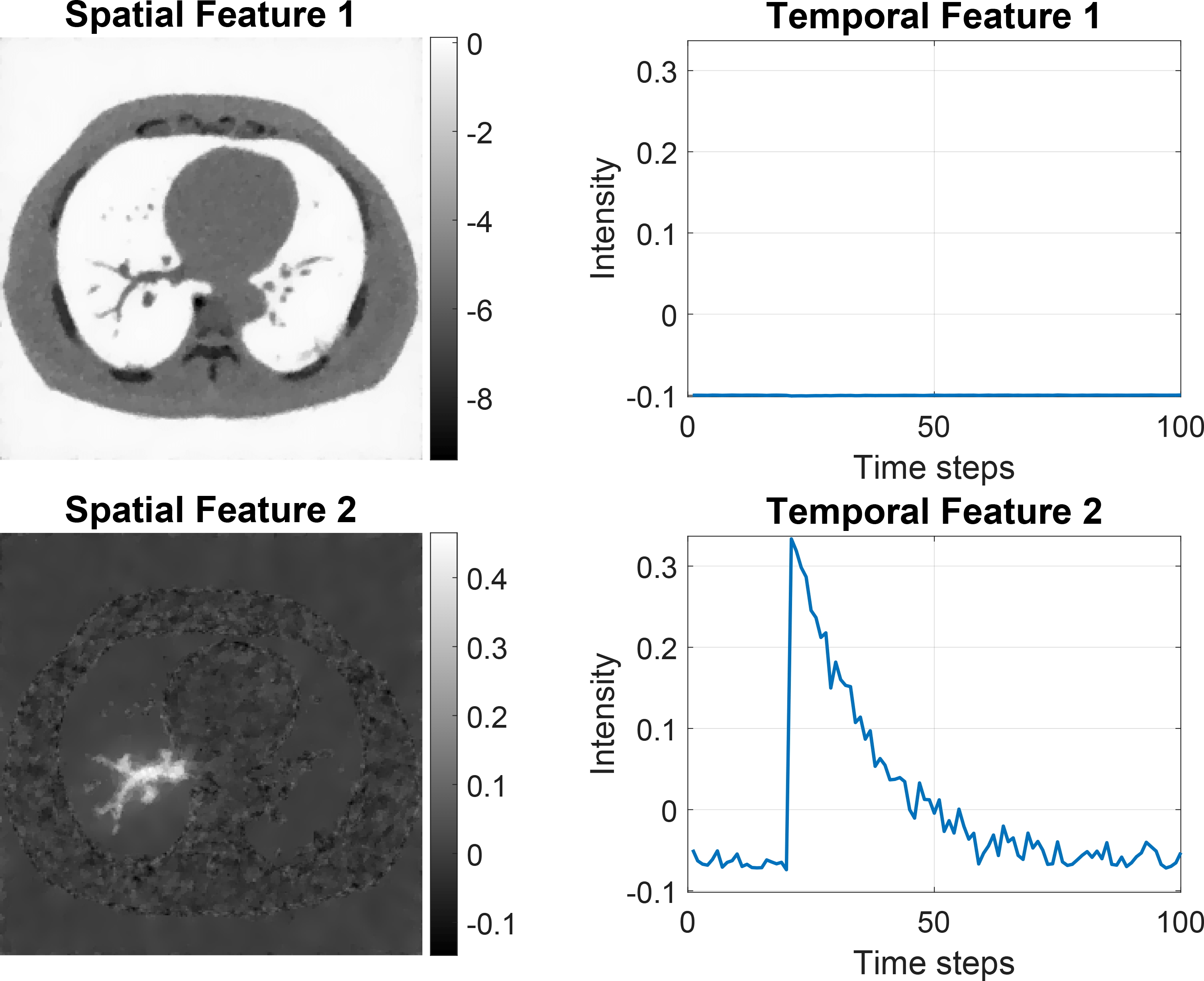}%
                \label{fig:phantomVessel_gradTV_PCA_0.03_nTheta12}}
                \caption{Results for the vessel phantom with $\vert \mathcal{I}_t \vert =12$ angles per time step and 3\% Gaussian noise. Shown are the leading extracted features for the \ref{eq:NMF model BC} model (left) and  for \texttt{gradTV\_PCA} (right).}
                \label{fig:phantomVessel:noise3:nTheta12}
            \end{figure}
            Further experiments show that the quality of the extracted components of \ref{eq:NMF model BC-X} decreases steadily for lower angles until the main features cannot be identified anymore for $\vert \mathcal{I}_t \vert \leq 8.$ \ref{eq:NMF model BC} produces inferior results and cannot extract reasonable components anymore for $\vert \mathcal{I}_t \vert \leq 10.$
            
            In comparison, both separated approaches \texttt{gradTV\_PCA} and \texttt{gradTV\_NMF} are still able to extract decent features for $\vert \mathcal{I}_t \vert = 7.$ For $\vert \mathcal{I}_t \vert \leq 6,$ the performance of both methods decreases significantly.
            
            Similar results for \texttt{gradTV\_PCA} could be obtained for 3\% noise and $\vert \mathcal{I}_t \vert = 12,$ which are shown in Figure \ref{fig:phantomVessel_gradTV_PCA_0.03_nTheta12}. Its constant feature is inferior to the one of \ref{eq:NMF model BC} in Figure \ref{fig:phantomVessel_BC_0.03_nTheta12} due to the additional nonnegativity constraint of the NMF model. However, the details of the vessel in the dynamic spatial feature of \ref{eq:NMF model BC} are lost due to the choice of large regularisation parameter $\tau$ and the temporal features are affected by several disturbances. Further tests with the noise level of 3\% showed that both joint methods are not able to recover the underlying features for $\vert \mathcal{I}_t \vert \leq 10,$ while the separated approaches gives still acceptable results for $\vert \mathcal{I}_t \vert =  6.$ 

            
            \begin{table}[t!]
                \ra{1.25}
                \resizebox{\columnwidth}{!}{%
                \begin{filecontents*}{tables/qualityMeasures_phantomVesselWithoutBCX_BC.csv}
noise,numAng,BCPSNR,BCSSIM,phantom1,BCXPSNRX,BCXSSIMX,phantom2,gradTVPCAPSNR,gradTVPCASSIM
1\%,7,(34.130),(0.9016),,32.969,0.8382,,31.463,0.8414
1\%,12,35.050,0.9068,,33.919,0.8496,,34.309,0.8839
3\%,12,30.148,0.7484,,28.119,0.5708,,29.375,0.6698
\end{filecontents*}
\pgfplotstabletypeset[
                    col sep=comma,
                    every head row/.style={
                        before row={
                            \toprule
                            && \multicolumn{2}{c}{\ref{eq:NMF model BC}} & \phantom{i} & \multicolumn{2}{c}{\ref{eq:NMF model BC-X}} & \phantom{i} & \multicolumn{2}{c}{\texttt{gradTV}}\\
                            \cmidrule{3-4} \cmidrule{6-7} \cmidrule{9-10}
                        },
                        after row=\midrule,
                    },
                    every last row/.style={
                        after row=\bottomrule},
                    every even row/.style={
                        before row={\rowcolor[gray]{0.9}}
                    },
                    columns/noise/.style       ={column name=Noise, string type},
                    columns/numAng/.style       ={column name=$\vert \mathcal{I}_t \vert$},
                    columns/BCPSNR/.style       ={column name=PSNR, string type},
                    columns/BCSSIM/.style       ={column name=SSIM, string type},
                    columns/phantom1/.style     ={column name={}},
                    columns/BCXPSNRX/.style     ={column name=PSNR, string type},
                    columns/BCXSSIMX/.style     ={column name=SSIM, string type},
                    columns/phantom2/.style     ={column name={}},
                    columns/gradTVPCAPSNR/.style      ={column name=PSNR, string type},
                    columns/gradTVPCASSIM/.style      ={column name=SSIM, string type},
                    ]
                    {tables/qualityMeasures_phantomVesselWithoutBCX_BC.csv}
                }
                \caption{Mean \ac{PSNR} and \ac{SSIM} values of the reconstruction results of the vessel phantom for different noise levels and numbers of projection angles. Values in brackets indicate that the dynamic part of the dataset in the corresponding experiment could not be reconstructed sufficiently well.}
    			\label{tab:phantomVessel:recQual}
            \end{table}
            
            The reconstruction quality of the experiments are shown in Table \ref{tab:phantomVessel:recQual}. Similar to the Shepp-Logan phantom, the joint approach \ref{eq:NMF model BC} produces the best results compared to all other methods in terms of the mean PSNR and SSIM values. Further experiments confirm this observation for $4 \leq \vert \mathcal{I}_t \vert \leq 11.$
            
            However, these observations have to be treated with caution. \ref{eq:NMF model BC} is not able to recover the dynamics for $\vert \mathcal{I}_t \vert \leq 10$ and 1\% noise. In the case of \ref{eq:NMF model BC-X}, the dynamics can be reconstructed to some degree within the angle range $9 \leq \vert \mathcal{I}_t \vert \leq 11,$ but are not recognizable anymore for $\vert \mathcal{I}_t \vert \leq 8.$ In the case of 3\% Gaussian noise, \texttt{gradTV} is still able to give acceptable reconstruction results for $\vert \mathcal{I}_t \vert = 10.$ For less angles, the reconstructed dynamics of \texttt{gradTV} get constantly worse until they are not apparent anymore for $\vert \mathcal{I}_t \vert \leq 6 .$
            
            The computation times of the experiments in Table \ref{tab:phantomVessel:recQual} range approximately from 7 to 15 minutes. The corresponding reconstructions can be found as video files in the Supplementary information.
\section{Conclusion}\label{sec:Conclusion and Outlook}
    In this work we consider dynamic inverse problems with the assumption that 
    the target of interest has a low-rank structure and can be efficiently 
    represented by spatial and temporal basis functions. This assumption leads 
    naturally to a reconstruction and low-rank decomposition framework. In 
    particular, we concentrate here on the \acl{NMF} as decomposition because it exhibits  three main advantages:
    \begin{enumerate}
    \item[i.)] It naturally incorporates the physical assumption of nonnegativity
    \item[ii.)] Basis functions are not restricted to being strictly orthogonal and therefore correspond more naturally to actual components
    \item[iii.)] It allows the flexibility to incorporate separate regularisation on each of the factorisation matrices 
    \end{enumerate}
    In particular, the last point is of importance, as it allows to consider different regularisers for spatial and temporal basis functions, and as such can be tailored to different applications.
    
    We then proposed two approaches to obtain a joint reconstruction and low-rank decomposition based on the \ac{NMF}, termed  \ref{eq:NMF model BC-X} and  \ref{eq:NMF model BC}. Both methods performed better than a baseline method, that computes a reconstruction with low-rank constraint followed by a subsequent decomposition.    
    In particular, the second \ref{eq:NMF model BC} model has shown to have a stronger regularising effect on the reconstructed features as well as the reconstruction, which can be simply obtained as $X=BC.$ 
    We believe this is due to the fact, that only the decomposition is recovered during the reconstruction without the need to build the reconstruction $X$ explicitly and hence the resulting features at the end exhibit a higher regularity. More importantly, if one considers a stationary operator in the complexity reduced \ref{eq:NMF model sBC} model we can obtain a considerable computational speed-up. Even though, due to constant projection angles the spatial basis functions are not as well recovered as in the non-stationary case, but the temporal features can be nicely extracted even for as low as 2 angles. This might be especially of interest in applications, where one is primarily interested in the underlying dynamics of the imaged target. 
    
    The primary limitation of the presented approach is the assumption on the decomposition of the target into spatial and temporal basis functions, as this does not allow for spatial movements in the target. 
    However, it opens up the possibility of combination with other methods, that do in fact allow for movements but assume a brightness consistency in the target, such as the optical flow constraint in CT \cite{burger2017variational}. Furthermore, the presented low-rank decomposition may be combined with a morphological motion model \cite{Gris:2019aa} to allow for a flexible and general model for dynamic inverse problems.



\section*{Acknowledgments} \label{sec:Acknowledgements}
This project was supported by the Deutsche Forschungsgemeinschaft (DFG, German Research Foundation) within the framework of RTG ``$\pi^3$: Parameter Identification -- Analysis, Algorithms, Applications'' -- Projektnummer 281474342/\\GRK2224/1.
This work was partially supported by the Academy of Finland Project 312123 (Finnish Centre of Excellence in Inverse Modelling and Imaging, 2018--2025), EPSRC grant EDCLIRS (EP/N022750/1) as well as CMIC-EPSRC platform grant (EP/M020533/1).

\bibliographystyle{siam}
\bibliography{references/references,references/Inverse_problems_references_2018,references/dyntomo_refs,references/lowrank_refs}


\newpage
\appendix

\section{Optimisation Techniques for NMF Problems} \label{app:sec:Optimisation Techniques for NMF Problems}
    The majority of optimisation techniques for NMF problems are based on alternating minimisation schemes. This is due to the fact that the corresponding cost function in \eqref{eq:NMF Minimisation Problem} is usually convex in $B$ for fixed $C$ and $C$ for fixed $B$ and non-convex in $(B, C)$ together, which yields algorithms of the form 
    \begin{align*}
        B^{[d+1]} &\coloneqq \arg\min_{B\geq 0}  \mathcal{F}(B, C^{[d]}), \\
        C^{[d+1]} &\coloneqq \arg\min_{C\geq 0}  \mathcal{F}(B^{[d+1]}, C).
    \end{align*}
    Typical minimisation approaches are based on alternating least squares methods, multiplicative algorithms as well as projected gradient descent and quasi-newton methods \cite{cichocki09bookNMF}. In this work, we focus on the derivation of multiplicative update rules based on the so-called \textit{Majorise-Minimisation} (MM) principle \cite{Lange13}. This approach allows the derivation of multiplicative update rules for non-standard NMF cost functions and gives therefore the flexibility to adjust the discrepancy and penalty terms according to the NMF model motivated by the corresponding application \cite{FM18}. What is more, the update rules consist only of multiplications and summations of matrices, which allow very simple implementations of the algorithms and ensure automatically the nonnegativity of the iterates $ B $ and $ C $ without the need of any inversion process, provided they are initialised nonnegative.
    \subsection{Multiplicative Algorithms} \label{app:subsec:Multiplicative Algorithms}
        The works of Lee and Seung \cite{LS99,LS00} brought much attention to NMF methods in general and, in particular, the multiplicative algorithms, which they derived based on the MM principle for the standard case with the Frobenius norm and the Kullback-Leibler divergence as discrepancy terms.
        
        The main idea of the MM approach is to replace the original cost function $ \mathcal{F} $ by a majorizing so-called \textit{surrogate function} $ \mathcal{Q}_\mathcal{F}, $ which is easier to minimise and leads to the desired multiplicative algorithms due to its tailored construction.
        
\begin{figure}
        	\centering
        	\includegraphics[width=1\textwidth]{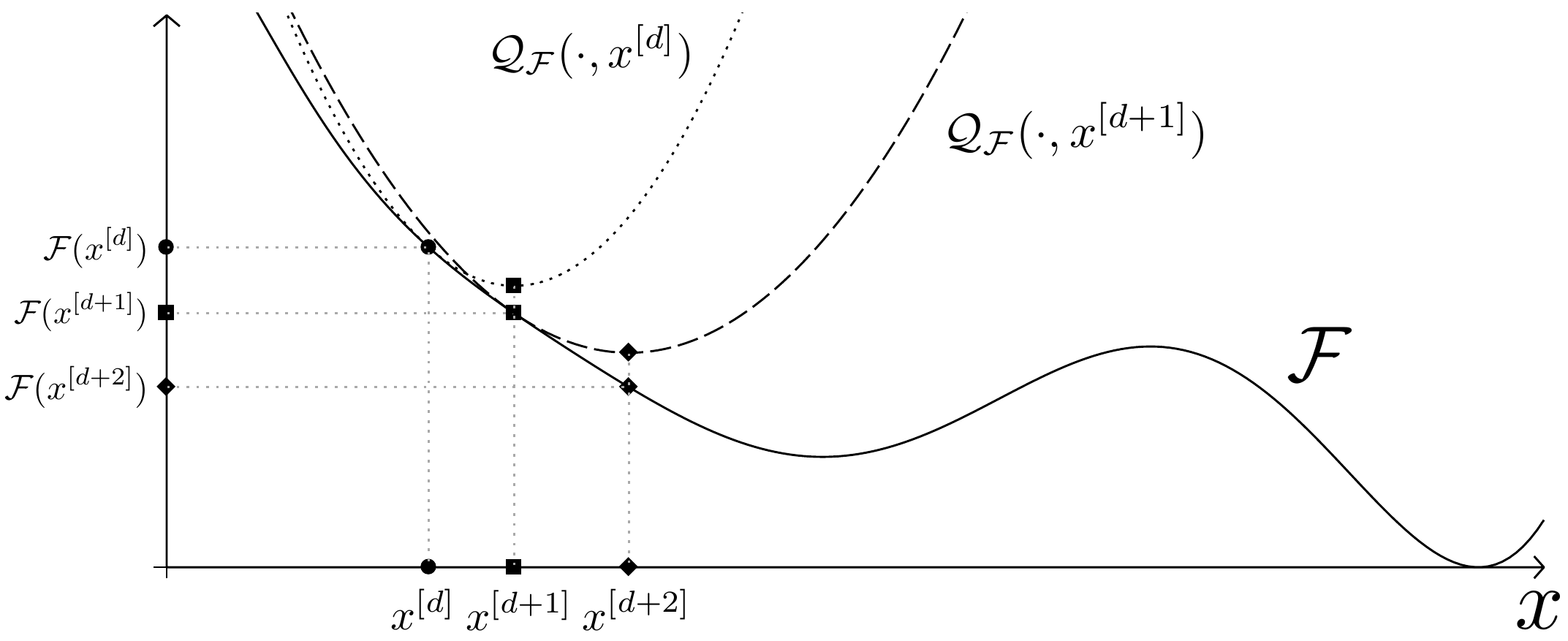}
        	\caption{Illustration of two iteration steps of the MM principle for a cost function $\mathcal{F}$ with bounded curvature and a surrogate function $\mathcal{Q}_\mathcal{F},$ which is strictly convex in the first argument.}
        	\label{fig:MMPrinciple}
        \end{figure}

        \begin{definition}[Surrogate Function]\label{def:surrogate function}
        	Let $ \Omega\subset \mathbb{R}^n $ be an open subset and $ \mathcal{F}:\Omega \rightarrow \mathbb{R} $ a function. Then $ \mathcal{Q}_\mathcal{F}:\Omega \times \Omega \rightarrow \mathbb{R} $ is called a \textbf{surrogate function} or \textbf{surrogate} of $ \mathcal{F}, $ if it fulfills the following properties:
        	\begin{enumerate}[i)]
        		\item $ \mathcal{Q}_\mathcal{F}(x, \tilde{x}) \geq \mathcal{F}(x)$ for all $ x, \tilde{x} \in \Omega $ \label{itm:def:surrogate property 1}
        		\item $ \mathcal{Q}_\mathcal{F}(x, x) = \mathcal{F}(x) $ for all $ x\in \Omega $ \label{itm:def:surrogate property 2}
        	\end{enumerate}
        \end{definition}
        The minimisation step of the MM approach is then defined by the update rule
        \begin{equation}\label{eq:Surrogate Update Rule}
            x^{[d+1]} \coloneqq \arg \min_{x\in \Omega} \mathcal{Q}_\mathcal{F}(x, x^{[d]}),
        \end{equation}
        assuming that the $\arg \min_{x\in \Omega} \mathcal{Q}_\mathcal{F}(x, \tilde{x})$ exists for all $\tilde{x} \in \Omega. $ Due to the defining properties of a surrogate function in Definition \ref{def:surrogate function}, the monotonic decrease of the cost function $ \mathcal{F} $ is easily shown:
        \begin{equation} \label{eq:LQBP Monotone Decrease}
            \mathcal{F}(x^{[d+1]}) \leq \mathcal{Q}_\mathcal{F}(x^{[d+1]}, x^{[d]}) \leq \mathcal{Q}_\mathcal{F}(x^{[d]}, x^{[d]}) = \mathcal{F}(x^{[d]}).
        \end{equation}
        This principle is also illustrated in Figure \ref{fig:MMPrinciple}.
        Typical construction techniques lead to surrogate functions, which are strictly convex in the first component to ensure the unique existence of the corresponding minimiser. Furthermore, the surrogates must be constructed in such a way, that the minimisation in Equation \eqref{eq:Surrogate Update Rule} yields multiplicative updates to ensure the nonnegativity of the matrix iterates. Finally, another useful property is the separability of $\mathcal{Q}_\mathcal{F}$ with respect to the first variable. This ensures, that $\mathcal{Q}_\mathcal{F}(x,\tilde{x})$ can be written as a sum, where each component just depends on one entry of $ x $ and allows the derivation of the multiplicative algorithm via the zero gradient condition $\nabla_x \mathcal{Q}_\mathcal{F} = 0.$
        
        One typical construction method is the so-called \textit{Quadratic Upper Bound Principle} (QUBP) \cite{BL88,Lange13}, which forms one of the main approaches to construct suitable surrogate functions for NMF problems. Nice overviews of other construction principles, which will not be used in this work, can be found in \cite{Lange13,LHY00}. The QUBP is described in the following Lemma.
    	\begin{lemma} \label{lem:Quadratic Upper Bound Principle}
    		Let $ \Omega\subset \mathbb{R}^n $ be an open and convex subset, $ \mathcal{F}:\Omega \rightarrow \mathbb{R} $ twice continuously differentiable with bounded curvature, i.e.\ there exists a matrix $ \Lambda\in \mathbb{R}^{n\times n}, $ such that $ \Lambda - \nabla^2\mathcal{F}(x) $ is positive semi-definite for all $ x\in \Omega. $ We then have
    		\begin{align*}
    		\mathcal{F}(x) &\leq \mathcal{F}(\tilde{x}) + \nabla \mathcal{F}(\tilde{x})^\intercal (x-\tilde{x}) + \dfrac{1}{2} (x-\tilde{x})^\intercal \Lambda (x-\tilde{x}) \quad \forall x, \tilde{x} \in \Omega \\
    		&\eqqcolon \mathcal{Q}_\mathcal{F}(x, \tilde{x}), \nonumber
    		\end{align*}
    		where $ \mathcal{Q}_\mathcal{F} $ is a surrogate function of $ \mathcal{F}. $
    	\end{lemma}
    	This is a classical result based on the second-order Taylor polynomial and will not be proven here.
    	
    	If the matrix $ \Lambda $ is additionally symmetric and positive definite, it can be shown \cite{FM18} that the update rule for $ x $ according to \eqref{eq:Surrogate Update Rule} via the zero gradient condition $ \nabla_{x}\mathcal{Q}_\mathcal{F}(x, \tilde{x})=0 $ gives the unique minimiser
    	\begin{equation} \label{eq:QUBP Update Rule}
    	    x_{\tilde{x}}^* = \tilde{x} - \Lambda^{-1} \nabla \mathcal{F}(\tilde{x}).
    	\end{equation}
    	In this work, we will only apply the QUBP for quadratic cost functions $F,$ whose Hessian is automatically a constant matrix. For these functions, typical choices of $ \Lambda $ are diagonal matrices of the form
    	\begin{equation} \label{eq:Lambda Matrix}
    	    \Lambda (\tilde{x})_{i i} \coloneqq \dfrac{(\nabla^2 f (\tilde{x})\ \tilde{x})_i + \kappa_i}{\tilde{x}_i},
    	\end{equation}
    	which are dependent on the second argument of the corresponding surrogate $\mathcal{Q}_\mathcal{F}(x, \tilde{x}).$ The parameters $ \kappa_i\geq 0, $ are constants and have to be chosen depending on the considered penalty terms of the NMF cost function.
    	
    	The diagonal structure of $ \Lambda(\tilde{x}) $ ensures its simple invertibility, the separability of the corresponding surrogate and the desired multiplicative algorithms based on \eqref{eq:Surrogate Update Rule}. Hence, the update rule in \eqref{eq:QUBP Update Rule} can be viewed as a gradient descent approach with a suitable stepsize defined by the diagnoal matrix $ \Lambda. $

\section{Derivation of the Algorithms} \label{app:sec:Derivation of the Algorithms}
    In this section, we derive the multiplicative update rules for the NMF minimisation problems in \ref{eq:NMF model BC-X} and \ref{eq:NMF model BC}.
    \subsection{Model \ref{eq:NMF model BC-X}} \label{app:subsec:Model BC-X}
		\subsubsection{Algorithm for X}\label{app:subsubsec:Algorithm for X (BC-X)}
			We start first of all with the NMF model \ref{eq:NMF model BC-X} and the minimisation with respect to $X.$ The cost function of the NMF problem in \ref{eq:NMF model BC-X} for the minimisation with respect to $ X $ reduces to
			\begin{equation} \label{eq:Cost Function F wrt X}
			    \mathcal{F}(X) \coloneqq \underbrace{\sum_{t=1}^{T} \frac{1}{2} \Vert A_t X_{\bullet, t} - Y_{\bullet, t} \Vert_2^2 + \frac{\mu_{X}}{2} \Vert X\Vert_F^2 + \lambda_{X} \Vert X \Vert_1}_{\eqqcolon \mathcal{F}_1(X)} + \underbrace{\frac{\alpha}{2} \Vert X - BC \Vert_F^2}_{\eqqcolon \mathcal{F}_2(X)}
			\end{equation}
			by neglecting the constant terms. To apply the QUBP and to avoid fourth-order tensors during the computation of the Hessians, we use the separability of $ \mathcal{F}_1 $ with respect to the columns of $ X, $ i.e.\ it can be written as sum, where each term depends only on the respective column $X_{\bullet, t}.$ Hence, we write
			\begin{equation*}
			    \mathcal{F}_1(X) = \sum_{t=1}^{T} \left [ \frac{1}{2} \Vert A_t X_{\bullet, t} - Y_{\bullet, t} \Vert_2^2 + \frac{\mu_{X}}{2} \Vert X_{\bullet, t}\Vert_2^2 + \lambda_{X} \Vert X_{\bullet, t} \Vert_1 \right ] \eqqcolon \sum_{t=1}^{T} f_t(X_{\bullet, t}).
			\end{equation*}
			We can assume that $ X $ contains only strictly positive entries due to the strict positive initialisations of the multiplicative algorithms. Hence, the functions $ f_t $ are twice continuously differentiable despite the occuring $ \ell^1 $ regularisation term. The computations of the gradient and the Hessian of $ f_t $ are straightforward and we obtain
			\begin{align*}
			    \nabla f_t(X_{\bullet, t}) &= A_t^\intercal A_t X_{\bullet, t} - A_t^\intercal Y_{\bullet, t} + \mu_{X}X_{\bullet, t} + \lambda_{X} \vec{1}_{N\times 1},\\
			    \nabla^2 f_t(X_{\bullet, t})&= A_t^\intercal A_t + \mu_{X} I_{N\times N},
			\end{align*}
			where $ I_{N\times N} $ is the $ N\times N $ identity matrix. Choosing $ \kappa_n=\lambda_{X} $ for all $ n $ in \eqref{eq:Lambda Matrix}, we define the surrogate $ \mathcal{Q}_{f_t} $ according to Lemma \ref{lem:Quadratic Upper Bound Principle}. It is then easy to see, that
			\begin{equation*}
		    	\mathcal{Q}_{\mathcal{F}_1}(X, \tilde{X}) \coloneqq \sum_{t=1}^T \mathcal{Q}_{f_t} (X_{\bullet, t}, \tilde{X}_{\bullet, t})
			\end{equation*}
			defines a separable and convex surrogate function for $ \mathcal{F}_1. $ For $ \mathcal{F}_2, $ we set simply $ \mathcal{Q}_{\mathcal{F}_2}(X,\tilde{X}) \coloneqq \nicefrac{\alpha}{2} \Vert X-BC\Vert_F^2, $ such that we end up with
			\begin{equation*}
		    	\mathcal{Q}_\mathcal{F}(X, A) \coloneqq \mathcal{Q}_{\mathcal{F}_1}(X,A) + \mathcal{Q}_{\mathcal{F}_2}(X,A)
			\end{equation*}
			as a suitable surrogate for $ F. $ Based on the update rule in \eqref{eq:Surrogate Update Rule}, we consider the zero gradient condition $ \nabla_{X}\mathcal{Q}_\mathcal{F}(X, \tilde{X})=0 $ and compute
			\begin{align*}
		    	\dfrac{\partial \mathcal{Q}_F}{\partial X_{nt}}(X,\tilde{X}) &=\dfrac{\partial f_t}{\partial X_{nt}} (\tilde{X}_{\bullet,t}) + \left ( \Lambda(\tilde{X}_{\bullet,t}) (X_{\bullet, t} - \tilde{X}_{\bullet, t}) \right )_n + \dfrac{\alpha}{2} \dfrac{\partial}{\partial X_{nt}} \Vert X-BC \Vert_F^2 \\
		    	&=\left ( A_t^\intercal A_t \tilde{X}_{\bullet, t} \right)_n - \left( A_t^\intercal Y_{\bullet, t} \right)_n + \mu_{X} \tilde{X}_{nt} + \lambda_{X} \\
	    		&+\scaleeq{0.93}{\dfrac{\left ( (A_t^\intercal A_t + \mu_{X} I_{N\times N})\tilde{X}_{\bullet,t} \right )_n + \lambda_{X}}{\tilde{X}_{nt}}(X_{nt} - \tilde{X}_{nt} ) + \alpha (X_{nt} - (BC)_{nt})} \\
		    	&=\scaleeq{0.93}{- \left( A_t^\intercal Y_{\bullet, t} \right)_n + X_{nt} \dfrac{\left ( A_t^\intercal A_t \tilde{X}_{\bullet, t} \right)_n + \mu_{X} \tilde{X}_{nt} + \lambda_{X}}{\tilde{X}_{nt}} + \alpha (X_{nt} - (BC)_{nt})} \\
		    	&=0.
			\end{align*}
			Rearranging the equation leads to
			\begin{equation*}
		    	X_{nt} = \dfrac{\left( A_t^\intercal Y_{\bullet, t} \right)_n + \alpha (BC)_{nt}}{ \dfrac{\left ( A_t^\intercal A_t \tilde{X}_{\bullet, t} \right)_n + \mu_{X}\tilde{X}_{nt} + \lambda_{X}}{\tilde{X}_{nt}} + \alpha}.
			\end{equation*}
			We therefore have
			\begin{equation*}
		    	X_{\bullet, t} = \tilde{X}_{\bullet, t} \circ \dfrac{A_t^\intercal Y_{\bullet, t} + \alpha BC_{\bullet, t}}{A_t^\intercal A_t \tilde{X}_{\bullet, t} + (\mu_{X} + \alpha)\tilde{X}_{\bullet, t} + \lambda_{X} \vec{1}_{N\times 1}},
			\end{equation*}
			which yields the multiplicative update rule
			\begin{equation*}
		    	X_{\bullet, t} \leftarrow X_{\bullet, t} \circ \dfrac{A_t^\intercal Y_{\bullet, t} + \alpha BC_{\bullet, t}}{A_t^\intercal A_t X_{\bullet, t} + (\mu_{X} + \alpha)X_{\bullet, t} + \lambda_{X} \vec{1}_{N\times 1}}
			\end{equation*}
			based on \eqref{eq:Surrogate Update Rule}. Note that the correct choice of the matrix $ \Lambda $ together with the $ \kappa_i $ is crucial to ensure the multiplicative structure of the algorithm.
		\subsubsection{Algorithm for B}\label{app:subsubsec:Algorithm for B (BC-X)}
			The minimisation with respect to $ B $ reduces the cost function in \ref{eq:NMF model BC-X} to
			\begin{equation} \label{eq:Cost Function F wrt B}
		    	\mathcal{F}(B) \coloneqq \underbrace{\frac{\alpha}{2} \Vert BC - X \Vert_F^2 + \frac{\mu_{B}}{2} \Vert B\Vert_F^2 + \lambda_{B} \Vert B \Vert_1}_{\eqqcolon \mathcal{F}_1(B)} + \underbrace{\dfrac{\tau}{2} \TV(B)}_{\eqqcolon \mathcal{F}_2(B)} 
			\end{equation}
			and involves the TV regularisation on $ B $ of the NMF model. Analogously to the previous section, we use the separability of $ \mathcal{F}_1 $ and write
			\begin{align*}
		    	\scaleeq{0.95}{\mathcal{F}_1(B) = \sum_{n=1}^{N} \Big[ \dfrac{\alpha}{2}  \Vert X_{n,\bullet} - B_{n,\bullet} C \Vert_F^2 + \dfrac{\mu_{B}}{2} \Vert B_{n,\bullet} \Vert_2^2 + \lambda_{B} \Vert B_{n,\bullet}\Vert_1  \Big] \eqqcolon \sum_{n=1}^{N} f_n(B_{n,\bullet}).}
			\end{align*}
			By computing the gradients
			\begin{align*}
		    	\nabla f_n(B_{n,\bullet}) &= \alpha (B_{n,\bullet} C - X_{n,\bullet}) C^\intercal + \mu_{B} B_{n,\bullet} + \lambda_{B} \vec{1}_{1\times K} \\
		    	\nabla^2 f_n(B_{n,\bullet}) &= \alpha CC^\intercal + \mu_{B} I_{K\times K}
			\end{align*}
			and choosing $ \kappa_k = \lambda_{B} $ in \eqref{eq:Lambda Matrix}, we define analogously the surrogates $ \mathcal{Q}_{f_n}, $ which leads to the convex surrogate
			\begin{equation*}
		    	\mathcal{Q}_{\mathcal{F}_1} (B, \tilde{B}) \coloneqq \sum_{n=1}^{N} \mathcal{Q}_{f_n}(B_{n,\bullet}, \tilde{B}_{n,\bullet})
			\end{equation*}
			for $ \mathcal{F}_1. $ The derivation of a suitable surrogate for the TV regularisation term $\mathcal{F}_2$ is based on an approach different from the QUBP and shall not be discussed in detail. We just state the result and refer the reader for details to \cite{OBF09,defrise11TV,FM18}. A convex and separable surrogate function for $ \mathcal{F}_2 $ is given by
			\begin{equation} \label{eq:surrogate:TV}
		    	\mathcal{Q}_{\mathcal{F}_2}(B, \tilde{B}) = \dfrac{\tau}{2} \sum_{k=1}^K \sum_{n=1}^N \left [ P(\tilde{B})_{nk} (B_{nk} - Z(\tilde{B})_{nk})^2  \right ] + \mathcal{G}(\tilde{B}),
			\end{equation}
			with the matrices $ P(\tilde{B}), Z(\tilde{B})\in \mathbb{R}_{\geq 0}^{N\times K} $ defined in \eqref{eq:TV:MatrixP} and \eqref{eq:TV:MatrixZ} and a function $ \mathcal{G} $ depending only on the matrix $ \tilde{B}.$ Hence, we finally end up with $ \mathcal{Q}_\mathcal{F}(B, \tilde{B}) \coloneqq \mathcal{Q}_{\mathcal{F}_1}(B, \tilde{B}) + \mathcal{Q}_{\mathcal{F}_2}(B, \tilde{B}) $ as a suitable surrogate for $\mathcal{F}$.
			
			Similar to the computations in the previous paragraph, the zero gradient condition yields then
			\begin{equation*}
		    	\scaleeq{0.83}{\dfrac{\partial \mathcal{Q}_\mathcal{F}}{\partial B_{nk}}(B, \tilde{B}) = - \alpha (XC^\intercal)_{nk} + B_{nk} \dfrac{\alpha (\tilde{B}CC^\intercal)_{nk} \!+\! \mu_{\vec{B}}\tilde{B}_{nk} \!+\! \lambda_{\vec{B}}}{\tilde{B}_{nk}} + \tau P(\tilde{B})_{nk} (B_{nk} \!-\! Z(\tilde{B})_{nk}) = 0}
			\end{equation*}
			and therefore
			\begin{equation*}
		    	B_{nk} = \tilde{B}_{nk} \cdot \dfrac{\alpha (XC^\intercal)_{nk} + \tau P(\tilde{B})_{nk}Z(\tilde{B})_{nk} }{\alpha (\tilde{B}CC^\intercal)_{nk} + \mu_{B} \tilde{B}_{nk} + \lambda_{B} + \tau P(\tilde{B})_{nk}\tilde{B}_{nk} }.
			\end{equation*}
			Hence, we have the update rule
			\begin{equation*}
		    	B \leftarrow B \circ \dfrac{\alpha XC^\intercal + \tau P(B) \circ Z(B) }{\alpha BCC^\intercal + \mu_{B} B + \lambda_{B} \vec{1}_{N\times K} + \tau P(B) \circ B }.
			\end{equation*}
		\subsubsection{Algorithm for C}\label{app:subsubsec:Algorithm for C (BC-X)}
			The optimisation with respect to the matrix $ C $ can be tackled analogously with the QUBP and will not be described in detail. In this case, the cost function can be reduced to well-known regularised NMF problems \cite{demol12}, which leads to the update rule
			\begin{align*}
	    		C \leftarrow C \circ \dfrac{\alpha B^\intercal X}{\alpha B^\intercal BC + \mu_{C} C + \lambda_{C} \vec{1}_{K\times T}}.
			\end{align*}
    \subsection{Model \ref{eq:NMF model BC}} \label{app:subsec:Model BC}
        In this section, we discuss the computation of the optimisation algorithms for the NMF model \ref{eq:NMF model BC}.
        \subsubsection{Algorithm for B}\label{app:subsubsec:Algorithm for B (BC)}
            In this case, the cost function reduces to
            \begin{align*}
                \mathcal{F}(B) \coloneqq \underbrace{\sum_{t=1}^{T} \dfrac{1}{2}  \Vert A_t (BC)_{\bullet,t} - Y_{\bullet, t} \Vert_2^2 + \dfrac{\mu_{B}}{2} \Vert B \Vert_F^2 + \lambda_{B} \Vert B\Vert_1}_{\eqqcolon \mathcal{F}_1(B)} + \underbrace{\dfrac{\tau}{2} \TV(B)}_{\eqqcolon \mathcal{F}_2(B)}.
            \end{align*}
            Analogously to the previous cases, we analyze the functions $ \mathcal{F}_1 $ and $ \mathcal{F}_2 $ separately. The difference is here, that $ \mathcal{F}_1 $ is not separable with respect to the rows of $ B $ due to the discrepancy term and therefore, it is necessary to compute the gradient and the Hessian of the whole function $ \mathcal{F}_1. $ Hence, the gradient $ \nabla \mathcal{F}_1(B)$ is an $ N\times K $ matrix and the Hessian $ \nabla^2 \mathcal{F}_1(B) $ a fourth-order tensor, which are given by their entries
            \begin{align*}
                \nabla \mathcal{F}_1(B)_{n k} &= \scaleeq{0.94}{\sum_{t=1}^T C_{k t} \left( A_t^\intercal A_t (B C)_{\bullet, t} \right)_{n} - \sum_{t=1}^T C_{k t} \left(A_t^\intercal Y_{\bullet, t}\right)_{n} + \mu_{B}B_{n k} + \lambda_{B}},\\
                \nabla^2 \mathcal{F}_1(B)_{(n, k), (\tilde{n}, \tilde{k})} &= \sum_{t=1}^T C_{\tilde{k} t} C_{k t} (A_t^\intercal A_t)_{n \tilde{n}} + \mu_{B} \delta_{(n, k), (\tilde{n}, \tilde{k})},
            \end{align*}
            where $ \delta_{(n, k), (\tilde{n}, \tilde{k})}=1 $ if and only if $ (n, k) = (\tilde{n}, \tilde{k}). $ The natural expansion of the quadratic upper bound principle given in Lemma \ref{lem:Quadratic Upper Bound Principle} is the ansatz function
            \begin{align*}
                \mathcal{Q}_{\mathcal{F}_1}(B, \tilde{B}) &\coloneqq \mathcal{F}_1(\tilde{B}) + \langle B-\tilde{B}, \nabla \mathcal{F}_1(\tilde{B})\rangle_F \\
                &+ \dfrac{1}{2} \sum_{(n, k)} \sum_{(\tilde{n}, \tilde{k})} (B-\tilde{B})_{n k} \Lambda(\tilde{B})_{(n, k), (\tilde{n}, \tilde{k})} (B-\tilde{B})_{\tilde{n} \tilde{k}}
            \end{align*}
            with the fourth order tensor
            \begin{equation*}
                \Lambda(\tilde{B})_{(n, k), (\tilde{n}, \tilde{k})} \coloneqq
                \begin{cases}
                    \dfrac{\sum_{(i, j)} \nabla^2 \mathcal{F}_1(\tilde{B})_{(n, k), (i, j)} \tilde{B}_{ij} + \lambda_{B}}{\tilde{B}_{n k}} &\text{for} \quad (n, k) = (\tilde{n}, \tilde{k}), \\
                    0 &\text{for} \quad (n, k) \neq (\tilde{n}, \tilde{k}),
                \end{cases}
            \end{equation*}
            where $ \langle \cdot, \cdot \rangle_F $ denotes the Frobenius inner product.
            
            Taking the same surrogate $ \mathcal{Q}_{\mathcal{F}_2} $ for the TV penalty term as in \eqref{eq:surrogate:TV}, we end up with the surrogate function
            \begin{equation*}
                \mathcal{Q}_\mathcal{F}(B, \tilde{B}) \coloneqq \mathcal{Q}_{\mathcal{F}_1}(B, \tilde{B}) + \mathcal{Q}_{\mathcal{F}_2}(B, \tilde{B})
            \end{equation*}
            for $ \mathcal{F}. $ Its partial derivative with respect to $ B_{nk} $ is given by
            \begin{align*}
                \dfrac{\partial \mathcal{Q}_\mathcal{F}}{\partial B_{nk}}(B) &= \scaleeq{0.94}{- \sum_{t=1}^T C_{k t} \left(A_t^\intercal Y_{\bullet, t}\right)_{n} + B_{nk} \frac{\sum_{t=1}^T C_{k t} \left( A_t^\intercal A_t (\tilde{B} C)_{\bullet, t} \right)_{n} + \mu_B \tilde{B}_{nk} + \lambda_B}{\tilde{B}_{nk}}} \\
                &+ \tau P(\tilde{B})_{nk} (B_{nk} - Z(\tilde{B})_{nk}).
            \end{align*}
            The zero-gradient condition gives then the equation
            \begin{align*}					
                B_{nk} = \tilde{B}_{nk} \Bigg( \dfrac{\sum_{t=1}^T C_{k t} \left(A_t^\intercal Y_{\bullet, t}\right)_{n} + \tau P(\tilde{B})_{nk} Z(\tilde{B})_{nk}}{\sum_{t=1}^T C_{k t} \left( A_t^\intercal A_t (\tilde{B} C)_{\bullet, t} \right)_{n} + \mu_{B} \tilde{B}_{nk} + \lambda_{B} + \tilde{B}_{nk} \tau P(\tilde{B})_{nk}} \Bigg),
            \end{align*}
            which can be extended to the whole matrix $ B. $ Therefore, based on \eqref{eq:Surrogate Update Rule}, we have the update rule
            \begin{equation*}
                B \leftarrow B \circ \Bigg( \dfrac{\sum_{t=1}^T A_t^\intercal Y_{\bullet, t} (C^\intercal)_{t,\bullet} + \tau P(B) \circ Z(B)}{\sum_{t=1}^T A_t^\intercal A_t (B C)_{\bullet, t} \cdot (C^\intercal)_{t,\bullet} + \mu_{B} B + \lambda_{B}\vec{1}_{N\times K} + \tau B \circ P(B)} \Bigg).
            \end{equation*}
        \subsubsection{Algorithm for C} \label{app:subsubsec:Algorithm for C (BC)}
            In this case, the cost function is separable with respect to the columns of $C,$ such that
            \begin{align*}
                \mathcal{F}(C) &\coloneqq \sum_{t=1}^{T} \dfrac{1}{2} \Vert A_t BC_{\bullet,t} - Y_{\bullet, t} \Vert_2^2 + \dfrac{\mu_{C}}{2} \Vert C_{\bullet,t} \Vert_2^2 + \lambda_{C} \Vert C_{\bullet,t} \Vert_1 \eqqcolon \sum_{t=1}^T f_t(C_{\bullet,t}).
            \end{align*}
            Hence, we can split the minimisation into the columns of $ C $ to use the standard QUBP without considering higher order tensors. We compute
            \begin{align*}
                \nabla f_t (C_{\bullet,t}) &= B^\intercal A_t^\intercal A_t (B C)_{\bullet,t} - B^\intercal A_t^\intercal Y_{\bullet,t} + \mu_{C} C_{\bullet,t} + \lambda_{C} \vec{1}_{K\times 1},\\
                \nabla^2 f_t (C_{\bullet,t}) &= B^\intercal A_t^\intercal A_t B + \mu_{C} I_{K\times K}.
            \end{align*}
            By choosing $ \kappa_k=\lambda_{C} $ for all $ k $ in \eqref{eq:Lambda Matrix}, we define $\mathcal{Q}_{f_t} (C_{\bullet,t}, \tilde{C}_{\bullet,t})$ as a surrogate function for $f_t$ according to Lemma \ref{lem:Quadratic Upper Bound Principle}. 
            The update rule in \eqref{eq:QUBP Update Rule} gives then
            \begin{equation*}
                C_{\bullet,t} = \tilde{C}_{\bullet,t} - \Lambda^{-1}(\tilde{C}_{\bullet,t}) \nabla f_t (\tilde{C}_{\bullet,t}),
            \end{equation*}
            which leads to
            \begin{equation*}
                C_{\bullet,t} \leftarrow C_{\bullet,t} \circ \dfrac{B^\intercal A_t^\intercal Y_{\bullet,t}}{B^\intercal A_t^\intercal A_t (BC)_{\bullet,t} + \mu_{C} C_{\bullet,t} + \lambda_{C} \vec{1}_{K\times 1} }.
            \end{equation*}
            
\section{Parameter Choice} \label{app:sec:parameter Choice}
    \begin{table}[H]
        \ra{1.1}
        \resizebox{\columnwidth}{!}{%
        \pgfplotstabletypeset[
            col sep=&, row sep=\\,
            string type,
            every head row/.style={
                before row={
                    \toprule
                    & \multicolumn{2}{c}{\ref{eq:NMF model BC}} & \phantom{i} & \multicolumn{2}{c}{\ref{eq:NMF model BC-X}} & \phantom{i} & \multicolumn{2}{c}{\texttt{gradTV}}\\
                    \cmidrule{2-3} \cmidrule{5-6} \cmidrule{8-9}
                },
                after row=\midrule,
            },
            every last row/.style={
                after row=\bottomrule},
            every even row/.style={
                before row={\rowcolor[gray]{0.9}}
            },
            columns/para/.style     ={column name=\textbf{Parameter}},
            columns/BC1/.style      ={column name=1\% noise},
            columns/BC3/.style      ={column name=3\% noise},
            columns/phantom1/.style ={column name={}},
            columns/BCX1/.style     ={column name=1\% noise},
            columns/BCX3/.style     ={column name=3\% noise},
            columns/phantom2/.style ={column name={}},
            columns/PCA1/.style     ={column name=1\% noise},
            columns/PCA3/.style     ={column name=3\% noise},
            ]{
                para & BC1 & BC3 & phantom1 & BCX1 & BCX3 & phantom2 & PCA1 & PCA3 \\
    			{$\alpha$} & -- &	-- && 70 & 70 && -- & --\\
    			{$\mu_C$} & 0.1 & 0.1 && 0.1 & 0.1 && -- & --\\
    			{$\tau$} & 10 & 50 && 6 & 20 && -- & --\\
    			{$\rho_{\text{grad}}$} & -- & -- && -- & -- && {$1\cdot 10^{-3}$} & {$8\cdot 10^{-4}$} \\
    			{$\rho_{\text{thr}}$} & -- & -- && -- & -- && {$7\cdot 10^{-4}$} & {$1\cdot 10^{-3}$}\\
    			{$\rho_{\text{TV}}$} & -- & -- && -- & -- && {$1\cdot 10^{-2}$} & {$2.5\cdot 10^{-2}$}\\
    			{$\tilde{\mu}_{C}$} & -- & -- && -- & -- && {$0.1$} & {$0.1$}\\
            }
        }
		\caption{Parameter choice of the considered methods for the dynamic Shepp-Logan Phantom for 1\% and 3\% Gaussian noise.}
		\label{tab:dynShepp:parameterChoice}
    \end{table}
    
    \begin{table}[H]
        \ra{1.1}
        \resizebox{\columnwidth}{!}{%
        \pgfplotstabletypeset[
            col sep=&, row sep=\\,
            string type,
            every head row/.style={
                before row={
                    \toprule
                    & \multicolumn{2}{c}{\ref{eq:NMF model BC}} & \phantom{i} & \multicolumn{2}{c}{\ref{eq:NMF model BC-X}} & \phantom{i} & \multicolumn{2}{c}{\texttt{gradTV}}\\
                    \cmidrule{2-3} \cmidrule{5-6} \cmidrule{8-9}
                },
                after row=\midrule,
            },
            every last row/.style={
                after row=\bottomrule},
            every even row/.style={
                before row={\rowcolor[gray]{0.9}}
            },
            columns/para/.style     ={column name=\textbf{Parameter}},
            columns/BC1/.style      ={column name=1\% noise},
            columns/BC3/.style      ={column name=3\% noise},
            columns/phantom1/.style ={column name={}},
            columns/BCX1/.style     ={column name=1\% noise},
            columns/BCX3/.style     ={column name=3\% noise},
            columns/phantom2/.style ={column name={}},
            columns/PCA1/.style     ={column name=1\% noise},
            columns/PCA3/.style     ={column name=3\% noise},
            ]{
                para & BC1 & BC3 & phantom1 & BCX1 & BCX3 & phantom2 & PCA1 & PCA3 \\
    			{$\alpha$} & -- &	-- && {$3\cdot 10^{2}$} & {$3\cdot 10^{2}$} && -- & --\\
    			{$\mu_C$} & 1 & 1 && 1 & 1 && -- & --\\
    			{$\tau$} & {$1.3\cdot 10^{2}$} & {$4.3\cdot 10^{2}$} && {$90$} & {$3\cdot 10^{2}$} && -- & --\\
    			{$\rho_{\text{grad}}$} & -- & -- && -- & -- && {$2\cdot 10^{-4}$} & {$8\cdot 10^{-5}$} \\
    			{$\rho_{\text{thr}}$} & -- & -- && -- & -- && {$2\cdot 10^{-4}$} & {$2.5\cdot 10^{-4}$}\\
    			{$\rho_{\text{TV}}$} & -- & -- && -- & -- && {$2\cdot 10^{-2}$} & {$4\cdot 10^{-2}$}\\
    			{$\tilde{\mu}_{C}$} & -- & -- && -- & -- && {$0.1$} & {$0.1$}\\
            }
        }
		\caption{Parameter choice of the considered methods for the vessel phantom for 1\% and 3\% Gaussian noise.}
		\label{tab:phantomVessel:parameterChoice}
    \end{table}

\end{document}